\documentclass[12pt,a4paper]{article}
\usepackage[margin=2.2cm]{geometry}

\usepackage{graphicx}
\usepackage{enumerate}
\usepackage{algpseudocode}
\usepackage{xspace}
\usepackage{boxedminipage}
\usepackage[ruled]{algorithm}
\usepackage{subfigure}
\usepackage{booktabs}
\usepackage{tabularx}
\usepackage{hyperref} 
\usepackage{times}
\usepackage{cite}
\usepackage{amsmath,amssymb,amsxtra,amsthm}
\usepackage{subfigure}
\usepackage{geometry}                

\usepackage[T1]{fontenc}
\usepackage[utf8]{inputenc}
\usepackage{authblk}

\newcommand{\Proofname}{Proof}
 
\newtheoremstyle
    {note}
    {}
    {}
    {}
    {}
    {\bfseries}
    {.}
    {1.0ex}
    {}
    
\swapnumbers
\numberwithin{equation}{section}

\theoremstyle{note}

\newtheorem{Rem}[subsection]{Remark}

\newtheorem{Pbm}[subsection]{Problem}

\providecommand{\highlight}[1]{{\color{red}#1}}
\usepackage{ifthen}
         \provideboolean{showchanges}
         \setboolean{showchanges}{true}
         \provideboolean{newchanges}
         \setboolean{newchanges}{true}
         \providecommand{\changes}[1]{
           \ifthenelse{\boolean{showchanges}}{{\highlight{#1}}}{#1}
         }
 \providecommand{\newchanges}[1]{
           \ifthenelse{\boolean{newchanges}}{{\highlight{#1}}}{#1}
         }
         \providecommand{\changefromto}[3][replace with]{
           \ifthenelse{\boolean{showchanges}}
           {{\sout{#2}\margnote{#1}}{\highlight{#3}}}
           {#3\xspace}
         }
         \providecommand{\ChangePar}[2]{
           \ifthenelse{\boolean{showchanges}}
           {{\par$\mapsfrom$ \textcolor{red!20}{#1}}{\par$\mapsto$ \textcolor{blue}{#2}}}
           {\par #2}
         }
         \providecommand{\InsertPar}[1]{
           \ifthenelse{\boolean{showchanges}}
           {{\par$\mapsto$ \textcolor{blue}{#1}}}
           {\par #1}
         }

\provideboolean{showdelete}
\setboolean{showdelete}{false}
\newcommand{\delete}[1]{
  \ifthenelse{\boolean{showdelete}} {{\color{red}{#1}}}{}
}

\renewcommand{\vec}[1]{\ensuremath{\boldsymbol{#1}}}

\newcommand{\myall}{\ensuremath{\quad \forall}}
\newcommand{\vp}{\ensuremath{\varphi}}
\newcommand{\vpo}{\ensuremath{\varphi_{obs}}}
\newcommand{\Lag}{\ensuremath{\mathcal{L}}}
\newcommand{\fv}{\ensuremath{\delta}}

\newcommand{\pdt}{\ensuremath{\partial_t}}

\newcommand{\mdt}[1]{\ensuremath{\partial^{\bullet}_{#1}}}

\newcommand{\Reals}{\ensuremath{{\mathbb{R}}}}

\newcommand{\diff}{\ensuremath{{\operatorname{d}}}}
\renewcommand{\O}{\ensuremath{{\Omega}}}

\newcommand{\normal}{\ensuremath{{\vec{\nu}}}}

\newcommand{\G}{\ensuremath{{\Gamma}}}
\newcommand{\Gt}{\ensuremath{{\G(t)}}}

\newcommand{\Go}{\ensuremath{{\G_{obs}}}}

\newcommand{\Lp}[1]{\ensuremath{\operatorname{L}_{#1}}}

\newcommand{\Hil}[1]{\ensuremath{\operatorname{H}^{#1}}}

\newcommand{\V}{\ensuremath{{\mathbb{V}}}}

\newcommand{\ltwon}[2]{\ensuremath{\left\|#1\right\|}_{\Lp{2}\left(#2\right)}}

\usepackage{upgreek}
\newcommand{\lap}{\ensuremath{\Updelta}}

\newcommand{\T}{\ensuremath{{\mathcal{T}}}}

\newcommand{\Lagrange}{\ensuremath{{\Lambda^h}}}

\usepackage{color}

\definecolor{MyGreen}{rgb} {0.05,0.4,0.05}
\definecolor{RedViolet}{rgb} {0.1,0.1,0.75}

          \providecommand{\highlight}[1]{{\color{blue}#1}}
         \provideboolean{shownotes}
          \setboolean{shownotes}{true}
          \newcommand{\standout}[1]{\colorbox{red}{\textcolor{white}{#1}}}
          \newcounter{margnote}[page]
          \newcommand{\margnotemark}{{\standout{\footnotesize\upshape\texttt{\arabic{margnote}}}}}
          \newcommand{\margnote}[2][]{
            \ifthenelse{
              \boolean{shownotes}
            }{\stepcounter{margnote}\margnotemark\marginpar{
                \texttt{
                  \begin{minipage}{2cm}
                    \raggedright\tiny
                    \margnotemark{#1}: 
                    #2
                  \end{minipage}
            }}}{}
          }
          \provideboolean{showchanges}
          \setboolean{showchanges}{true}
          \providecommand{\chcolor}{\color{blue}}
          \providecommand{\changes}[1]{
            \ifthenelse{\boolean{showchanges}}
          	     {{\chcolor #1}}
          	     {#1
          }
          }
          \providecommand{\changefromto}[3][replace with]{
            \ifthenelse{\boolean{showchanges}}
            {{\sout{#2}\margnote{#1}}{\highlight{#3}}}
            {#3\xspace}
          }
          \providecommand{\ChangePar}[2]{
            \ifthenelse{\boolean{showchanges}}
            {{\par$\mapsfrom$ \textcolor{red!20}{#1}}{\par$\mapsto$ \textcolor{blue}{#2}}}
            {\par #2}
          }
          \providecommand{\InsertPar}[1]{
            \ifthenelse{\boolean{showchanges}}
            {{\par$\mapsto$ \textcolor{blue}{#1}}}
            {\par #1}
         }
\usepackage[normalem]{ulem}
\usepackage{hyperref}
\usepackage[british]{babel}
\renewenvironment{itemize}{\begin{list}{\labelitemi}{\leftmargin=1em} }
{\end{list}
}
\newcommand \beq{\begin{equation}}
\newcommand \eeq{\end{equation}}

\def\eps{\varepsilon}

\def\bpmat{\begin{pmatrix}}
\def\epmat{\end{pmatrix}}
\def\mathref#1{\ifmmode\mathrm{(\ref{#1})}\else{\rm(\ref{#1})}\fi} 
\def\nref#1{\ifmmode\mathrm{\ref{#1}}\else{\rm\ref{#1}}\fi}

\title{Whole cell tracking through the optimal control of geometric evolution laws}
\author[a]{Konstantinos N. Blazakis}
\author[a]{Anotida  Madzvamuse}
\author[b]{Constantino-Carlos Reyes-Aldasoro}
\author[a]{Vanessa Styles}
\author[a]{Chandrasekhar Venkataraman\thanks{Corresponding author. Email : C.Venkataraman@sussex.ac.uk}}
\affil[a]{University of Sussex, Department of Mathematics, Brighton, UK, BN1 9QH}
\affil[b]{City University London, School of Mathematical Sciences, Computer Science and Engineering, Northampton Square, London, UK, EC1V OHB}

\date{}

\begin{document}

\maketitle
\abstract{
Cell tracking algorithms which automate and systematise the analysis of time lapse image data sets of cells are an indispensable tool in the modelling and understanding of cellular phenomena. In this study we present a theoretical framework and an algorithm for whole cell tracking. Within this work we consider that ``tracking'' is equivalent to a dynamic reconstruction of the whole cell data  (morphologies) from static image datasets. The novelty of our work is that the tracking algorithm is driven by a model for the motion of the cell. This model may be regarded as a simplification of a recently developed physically meaningful model for cell motility. The resulting problem is  the optimal control of a geometric evolution law and we discuss the formulation and numerical approximation of the optimal control problem. The overall goal of this work is to design a framework for cell tracking  within which the recovered data reflects the physics of the forward model.  A number of numerical simulations are presented that illustrate the applicability of our approach.\\

}

\section{Introduction}
\label{sec:Intro}
     Cell migration is a fundamental process in cell biology and is tightly linked to many important physiological and pathological events such as the immune response, wound healing, tissue differentiation, metastasis, embryogenesis, inflammation and tumour invasion \cite{bray2001cell}. Experimental advances provide techniques to observe migrating cells both  \emph{in vivo} and \emph{in vitro}. Inferring dynamic quantities from this static data is an important task that has many applications in biology and related fields.   The field of cell tracking arose from this need and is concerned with the  development of methods to track and analyse dynamic cell shape changes from a series of still images captured within a time frame (see for example \cite{xiong2010tools,meijering2012methods} for reviews). 

On the other hand, a major focus of current research is the derivation of  mathematical models for cell migration based on physical principles, e.g.,  \cite{venk11chemotaxis}. Furthermore, such models appear to show good qualitative and quantitative agreement with experimental observations of migrating cells. Despite this, very little research has focused on incorporating these mathematical modelling advances into appropriate cell tracking algorithms. In a related work, we investigated fitting parameters in models for cell motility to experimental image data sets of migrating cells where observations of both the position of the cells and the concentrations of cell-resident proteins related to motility were available \cite{2013arXiv1311.7602C}.  

In this study  we present a first step towards the development of a framework for cell tracking based on novel  models of cell motility. Specifically, we propose a cell tracking algorithm which can be thought of as fitting a simplified, yet physically meaningful, model for cell migration to experimental observations and data. 
We focus on the setting, prevalent in cell tracking problems, where only the position of the cell at a series of discrete times is available and no further biological information is given. We present a mathematical model based on physical principles for the cell movement that consists of a geometric evolution equation. We then formulate an inverse problem, which takes the form of a PDE constrained optimisation problem, for fitting the model to the static experimental observations. To solve the optimisation problem we propose an algorithm  based on the optimal control of geometric evolution laws \cite{hausser2010influence,hausser2012control}. 

The objective of this study is to serve as a useful first step in the development of cell tracking algorithms in which the underlying model for the evolution is based on physical principles, rather than purely geometric considerations. In this setting, one hopes to attain estimates of motility-related features such as trajectories, velocities, persistence lengths, circularity, etc., which reflect the physics underlying the model. We illustrate the fact that the tracking procedure we propose allows us to incorporate physically important aspects of cell migration by including volume conservation in the model for the evolution. This is motivated by the observation that, for many cells, while the surface area of the cell membrane may change significantly during migration the volume enclosed by the cell remains roughly constant \cite{shao2010computational}.  Of course other physical aspects of the migration could be included in the model, such as a spontaneous curvature of the membrane which is relevant for more complex models of cell motility involving the Helfrich model \cite{marth2013signaling}. 
     
The remainder of our discussion proceeds as follows. In \S \ref{sec:prob} we briefly describe the problem of cell tracking and introduce our approach to cell tracking, which may be regarded as fitting a mathematical model to experimental image data sets. We present the geometric evolution law model we seek to fit, which is a simplification of recently developed models in the literature that show good agreement with experiments \cite{shao2010computational,ziebert2011model,neilson2011modelling,neilson2011chemotaxis,venk11chemotaxis,shao2012coupling,marth2013signaling}. We finish \S \ref{sec:prob} by reformulating our model into the phase field framework, which appears more suitable for the problem in hand, and we formulate the cell tracking problem as a PDE constrained optimisation problem. In \S \ref{sec:implementation} we propose an algorithm for the resolution of the PDE	constrained optimisation problem and we discuss some practical aspects related to the implementation. In particular we note that the theoretical and computational framework may be applied directly to multi-cell image data sets and raw image data sets (of sufficient quality) without segmentation. In \S \ref{sec:application} we present some numerical examples for the case of $2d$ single and multi-cell image data sets. Finally in \S \ref{sec:conc} we present some conclusions of our study and discuss future extensions and applications of the work.

\section{Problem Formulation}
\label{sec:prob}
\subsection{Approaches to cell tracking}
\label{subsec:cell_tracking}
Our focus is on developing whole cell tracking algorithms which track the morphology of the cell rather than particle tracking in which particles such as the cell centroid or cell resident proteins or (macro-)molecules are tracked \cite{henry2013phagosight}. A number of approaches have proved successful in cell tracking with level-set \cite{bosgraaf2009analysis} or electrostatic based methods among the  most widely used \cite{tyson2010high}. One feature of such  methods is that the trajectories they generate are not physical in nature rather they are designed with the goal of achieving nice geometric properties, e.g., equidistribution of vertices, smoothness of the trajectories and so on.   Our approach differs to these purely geometric approaches  in that we start with a model derived from physical principles and it is this model for the evolution that drives the tracking algorithm. In this sense our approach is similar in spirit to the parameter identification procedure described in \cite{2013arXiv1311.7602C} as in both studies the goal may be regarded as fitting a mathematical model to experimental image data sets.

We now summarise the main problems that must be addressed by a cell tracking algorithm. In general cell tracking consists of three main steps
\begin{enumerate}
\item Segmentation: In this step the raw image data set is processed and the cells are separated from the background in each frame.
\item  Matching: The cells segmented in the first step must then be associated from frame to frame (note this is only relevant in the case of multiple cell image data sets)  such that where possible (in practice cells may disappear or spontaneously appear in images) there is a one-to-one map that uniquely associates individual cells from one frame to the next.
\item Linking: Finally the linking step consists of estimating dynamic data  from the associated segmented static cells. 
\end{enumerate}
In this work we will largely neglect the segmentation step. We assume either that we have segmented image data to work with or that the image data is of sufficient quality that the contrast between the cell and the background is clear and a simple thresholding step is sufficient to label the different cells. In the case of segmented image data, we assume this data consists of closed surfaces (or curves in $2d$)  that describe the boundaries of each individual cell.  We briefly comment on the role of segmentation in our approach and the fact that it may be unnecessary to extract the sharp contours of the cell boundaries in Remark \ref{Rem:segmentation}. For ease of presentation we will also describe the algorithm in the case of single cell image data and thus the matching step is redundant. However the theoretical aspects of our approach apply equally to multiple cell data and one potentially attractive feature of our cell tracking procedure is that the matching problem may be resolved 
implicitly without the need to associate cells form one frame to the next, c.f., Remark \ref{Rem:multi-cells}.
Our main focus is on the tracking step and in the remainder of this study we outline a theoretical framework, a practical implementation and numerical examples, for linking data at a series of discrete times which allows the  recovery of the whole cell morphologies in time.  
\subsection{Model}
\label{subsec:models}
As mentioned above in contrast to many of the existing approaches for cell tracking, the framework we propose in this study is based on fitting a model, derived from physical principles, for the motion of the cell to experimental image data. The general class of models to which our approach is applicable are PDE based models for the motion, where the movement of the cell membrane is described by a geometric evolution law. We introduce some notation for the formulation of the model.

We denote by $\G$ the cell membrane, which is assumed to be a closed smooth oriented $d-1$ dimensional hypersurface in $\Reals^d$,  $d=2,3$, with outward pointing unit normal $\normal$. Given a function $\eta$ defined in a neighbourhood of $\G$, the tangential or surface gradient of $\eta$ denoted by $\nabla_\G$ is defined as
\begin{equation}
\nabla_\G\eta:=\nabla\eta-\nabla\eta\cdot\normal\normal,
\end{equation}
where $\nabla$ denotes the Cartesian gradient in $\Reals^d$. The Laplace-Beltrami operator $\lap_\G$ is defined as the tangential divergence of the tangential gradient, i.e.,
\begin{equation}
\lap_\G\eta:=\nabla_\G\cdot\left(\nabla_\G\eta\right).
\end{equation}
The mean curvature $H$ of $\G$ with respect to the normal $\normal$ is defined as
\begin{equation}
H:=\nabla_\G\cdot\normal.
\end{equation}

In this study  we model the evolution of the cell membrane as being governed by volume conserved mean curvature flow with forcing, given by
\beq\label{eqn:MCF_forcing}
\begin{cases}
\vec V(\vec x,t)=\left(-\sigma H(\vec x,t)+\eta(\vec x,t)+\lambda_V(t)\right)\normal(\vec x,t)
\quad \text{on }\Gt,t\in(0,T],\\
\G(0)=\G^0,
\end{cases}
\eeq
where $\G$ is the closed surface that represents the cell membrane, $\vec V$ is the material velocity of $\G$, $\sigma$ is the surface tension and $\lambda_V(t)$ is a spatially uniform force accounting for volume conservation, physically this may be thought of as an interior pressure. The forcing function $\eta$ is the main driver of the directed migration.  The model we present is phenomenological and hence it is difficult to directly relate $\eta$ to biophysical processes. However, as positive values of  $\eta$ correspond to protrusive forces and negative values of  $\eta$  correspond to contractile forces one interpretation of the forcing function $\eta$ is that it accounts for both protrusive forces generated by polymerisation of actin at the leading edge of the cell and contractile forces generated by the action of myosin motors at the rear of the cell.

The evolution law (\ref{eqn:MCF_forcing}) is a simplification of a large class of models that arise in the modelling of cell motility which take the following form
\beq\label{eqn:evolution_law}
\begin{cases}
\vec V(\vec x,t)=\Bigg(g_1(H(\vec x,t))+g_2(\vec a(\vec x,t))+\lambda_V(t)\Bigg)\normal(\vec x,t)
\quad \text{on }\Gt,t\in(0,T],\\
\G(0)=\G^0,
\end{cases}
\eeq
where $g_1$ models the dependence of the evolution on geometric quantities, such as resistance of the membrane to stretching which could be modelled by mean curvature terms as in \eqref{eqn:MCF_forcing} or bending which can be modelled through the inclusion of Willmore or Helfrich flow type terms. The function $g_2$  appearing in (\ref{eqn:evolution_law}) captures the dependence of the evolution on a vector of bulk and/or surface resident species $\vec a$. The  surface resident species $\vec a$ could, for example, satisfy another PDE such as a surface reaction-diffusion system
\beq
\label{eqn:RDS}
\begin{cases}
\mdt{\vec V}\vec a+\vec a\nabla_{\Gt}\cdot\vec V-\vec D\lap_{\Gt}\vec a=\vec f(\vec a) \quad \text{on }\Gt,t\in(0,T],\\
\vec a(\cdot,0)=\vec a^0(\cdot)\quad \text{on }\G(0),
\end{cases}
\eeq
where $\vec{a} = (a_1, \dots, a_{n_a})^T$, $n_a$ is the number of chemical species involved, $a_i$ denotes the density of the $i$th chemical species, $\vec{V}$ is the material velocity of the surface, 
\begin{equation}\label{eqn:material_derivative}
\mdt{\vec V}{\vec{a}}:=\pdt{\vec{a}}+\vec{V}\cdot\nabla{\vec{a}},
\end{equation}
is the material derivative with respect to the velocity $\vec V$, $\vec{D}$ is a diagonal matrix of positive diffusion coefficients and $\vec{f}(\vec{a})$ is the nonlinear reaction. Models of the form (\ref{eqn:evolution_law})-(\ref{eqn:RDS}) have been used successfully to model cell motility in \cite{venk11chemotaxis,neilson2010use,neilson2011modelling,neilson2011chemotaxis} while models coupling evolution laws of the form (\ref{eqn:evolution_law}) to bulk PDEs (i.e., equations posed in the cell interior) have been considered in \cite{shao2012coupling,marth2013signaling,ziebert2011model}. 

Despite its simplicity the evolution law (\ref{eqn:MCF_forcing}) may be regarded as a prototype of the more complex models  for cell motility of the form (\ref{eqn:evolution_law})-(\ref{eqn:RDS}). The geometric evolution component (\ref{eqn:evolution_law}) is often the most challenging component of the model to solve numerically and developing an understanding of how to construct cell tracking algorithms assuming a  geometric evolution law based model for the motion is an important first step towards developing  tracking algorithms based on more realistic physical models. In many applications it is also the case that the only information available from the data is the position of the cell membrane and no adequate model for the biochemistry of the motility related species involved is available. Without any knowledge of the relevant biochemistry it is difficult to identify which motility related species should influence the evolution let alone propose how the evolution depends on their distribution (i.e., a $g_2$ in (\ref{eqn:evolution_law})) or a model for the species dynamics (i.e., an equation such as (\ref{eqn:RDS})).  Nevertheless one may still wish to extract dynamic quantities from static image data sets, in this setting it may be reasonable to consider the evolution law (\ref{eqn:MCF_forcing}) as a stand alone model for the motion as at least the mechanical aspects of the membrane evolution are accounted for through a physical model derived from basic physical principles. 

\subsection{An optimal control approach to cell tracking}
\label{subsec:optimal_control}
The cell tracking approach we consider in this study corresponds to the following problem.
\begin{Pbm}[Cell tracking]\label{Pbm:tracking}
Given an initial cell membrane position $\G^0$  and an observation of the position  $\Go$,  find a space-time distributed forcing $\eta$ such that the evolution of the cell membrane, $\Gt,t\in[0,T]$ satisfies (\ref{eqn:MCF_forcing}) with $\G(0)=\G^0$ and $\G(T)$ the position of the cell membrane at time $t=T$,  is close to $\Go$.
\end{Pbm}

As the volume enclosed by the cell may vary over the images it is inappropriate to enforce conservation of a constant volume. 
Instead we enforce, with the help of a Lagrange multiplier $\lambda_V(t)$,  that the volume enclosed by the cell is given by 
$\widetilde{V}(t)=V^0+\frac{t}{T}(V_{obs}-V^0)$, i.e. that the volume of the cell is a time-dependent linear interpolant of the volumes of the data. 


%

Problem \ref{Pbm:tracking} is an optimal control of a free boundary problem, where the free moving boundary problem is that of forced mean curvature flow and the control variable is the space time distributed forcing. The theory of optimal control of geometric evolution laws is in its infancy, in fact only recently has progress been made on the optimal control of parabolic equations on evolving surfaces even in the case of prescribed evolution \cite{vierling2014parabolic}. On the other hand the theory for the optimal control of semilinear parabolic equations is more mature (see, for example, \cite{troltzsch2010optimal}). We wish to exploit this fact and to this end we consider the phase field approximation of (\ref{eqn:MCF_forcing}) given by the Allen-Cahn equation;
 \begin{equation}\label{eqn:pf}
  \begin{cases}
 \partial _t \vp({\vec {{x}}},t)&=\Delta\vp({\vec {x}},t)-\frac{1}{\varepsilon^2}G^\prime(\vp({\vec {x}},t))-\frac{1}{\varepsilon}\left({c_{G}}
\eta({\vec {x}},t)-\lambda(t)\right)\text{ in }\O\times(0,T],\\
\nabla \vp\cdot\normal_\O&=0\text{ on }\partial\O\times(0,T],\\
\vp(\cdot,0)&=\vp^0(\cdot)\text{ in }\O,
\end{cases}
\eeq
where $\O\subset\Reals^d$ is a bulk domain, with normal $\normal_\O$, that contains $\G$,  $\vp^0$ is a diffuse interface representation of $\G^0$ and $\eps>0$ is a small parameter which governs the width of the diffuse interface. For details on the asymptotic analysis of (\ref{eqn:pf}) and the convergence (as $\eps\to0$) to a solution of \eqref{eqn:MCF_forcing} we refer the reader, for example, to \cite{chen1992generation,blowey1993curvature,bellettini1996anisotropic,brassel2011modified} and references therein.
The function $G$ appearing in (\ref{eqn:pf}) is a double well potential, for example the quartic potential
\beq\label{eqn:G_def}
G(\vp)=\frac{1}{4}\left(\vp^2-1\right)^2
\eeq 
which has minima at $\pm1$. The constant $c_G=\frac{1}{\sqrt{2}}\int_{-1}^1G(r)^{1/2}\diff r$ appearing in (\ref{eqn:pf}) is a scaling constant that depends on the double well potential. 
We enforce the time-dependent volume constraint following the approach of \cite{blowey1993curvature}. Specifically our diffuse interface formulation of the constraint on the enclosed  
volume is given by a constraint on $\int_\O[\varphi(\vec x,t)]_{+}\diff\vec x$, where $[a]_{+}=max(a,0)$.
We define $M_\vp$, the linear interpolant of $\int_\O[\varphi(\vec x,t)]_{+}\diff\vec x$ of the initial and target diffuse interface data  by
$$M_\vp(t):=\int_\O[\vp^0]_{+}+\frac{t}{T}\left([\vpo]_{+}-[\vp^0]_{+}\right)\diff \vec x,$$
 and determine $\lambda(t)$ in (\ref{eqn:pf}) such that $M_\vp(t)=\int_\O[\vp(\vec x,t)]_{+}\diff\vec x$.
 We have used $\lambda$ (rather than $\lambda_V$) for the Lagrange multiplier in (\ref{eqn:pf}) to reflect the fact that our constraint is on $\int_\O[\varphi(\vec x,t)]_{+}\diff\vec x$. However we shall refer to this constraint as a volume constraint  in order to highlight the physical feature the constraint is intended to model.  We also investigated an alternative approach to enforcing the volume constraint via penalising deviations from a target volume following  \cite{du2006simulating} (see also \cite{marth2013signaling}), in our numerical tests this strategy  proved less robust than the volume constraint proposed above.

To formulate the cell tracking problem as a PDE constrained optimal control problem we define the objective functional we shall seek to minimise as follows
\beq\label{eqn:objective}
J(\vp,\eta)=\frac{1}{2}\int_\O\left(\vp(\vec x,T)-\vpo(\vec x)\right)^2\diff \vec x+\frac{\theta}{2}\int_0^T\int_\O\eta(\vec x,t)^2\diff \vec x\diff t,
\eeq
where $\vpo$ is a diffuse interface representation of the observation $\Go$ and $\theta>0$ is a regularisation parameter. The first term on the right of \eqref{eqn:objective} is the so called fidelity term that measures the distance between the solution to the model and the target data and the second term is the regularisation which is necessary to ensure a well-posed problem (for example see \cite{troltzsch2010optimal}).

Our optimal control approach to the cell tracking problem may now be stated as the following minimisation problem.
\begin{Pbm}[Optimal control problem]\label{Pbm:control}
Given an initial diffuse interface representation of the cell membrane position $\vp^0$ and an observation of the position $\vpo$,  find a space-time distributed forcing $\eta^*:\O\times[0,T]\to\Reals$ such that with $\vp$ a solution of \eqref{eqn:pf} with initial condition $\vp(\cdot,0)=\vp^0$, the forcing $\eta^*$ solves the minimisation problem
\beq
\min_{\eta}J(\vp,\eta),\quad\text{with }J\text{ given by \eqref{eqn:objective}.}
\eeq
\end{Pbm}
\subsection{Optimality conditions}
\label{subsec:opt_con}

To apply the theory of optimal control of semilinear PDEs for the solution of the tracking problem, we briefly outline the derivation of the optimality conditions, for further details see for example \cite{hinze2009optimization,troltzsch2010optimal}. 
Introducing  the Lagrange multiplier  (adjoint state) $p$, we define the Lagrangian functional
\beq\label{eqn:Lagrangian}
\begin{split}
\Lag(\vp,\eta,p)=J(\vp,\eta)-&\int_0^T\int_\O\Bigg( \partial _t \vp({\vec {{x}}},t)-\Delta\vp({\vec {x}},t)
\\
&
+\frac{1}{\varepsilon^2}G^\prime(\vp({\vec {x}},t))+
\frac{1}{\varepsilon}\big(c_{G}\eta({\vec {x}},t)-\lambda(t)\big)\Bigg)p(\vec x,t)\diff \vec x\diff t.	
\end{split}
\eeq
Requiring stationarity of the Lagrangian with respect to the adjoint state yields the state equation (\ref{eqn:pf}) and requiring stationarity  of the Lagrangian, at the optimal control $\eta^*$ and associated optimal state $\vp^*$, with respect to the state and the control,  yields the (formal) first order optimality conditions 
\begin{align}
\label{eqn:oc1}
\fv_\vp \Lag(\vp^*,\eta^*,p)\vp&=0,\myall\vp:\vp(\vec x,0)=0,\\
\label{eqn:oc2}
\fv_\eta \Lag(\vp^*,\eta^*,p)\eta&=0,\myall\eta.
\end{align}
Condition (\ref{eqn:oc1}) yields the adjoint equation, which is the following  linear parabolic PDE for the adjoint state $p$,
 \begin{equation}\label{eqn:adjoint}
  \begin{cases}
 \partial _t p({\vec {{x}}},t)&=-\Delta p({\vec {x}},t)+\frac{1}{\varepsilon^2}G^{\prime\prime}(\vp({\vec {x}},t))p(\vec x,t)\quad\text{ in }\O\times(0,T],\\
\nabla p\cdot\normal_\O&=0\quad\text{ on }\partial\O\times(0,T],\\
p(\vec x,T)&=\vp(\vec x,T)-\vpo(\vec x)\quad\text{ in }\O.
\end{cases}
\eeq
Note that equation \eqref{eqn:adjoint} is posed backwards in time and hence is equipped with terminal conditions. 
Condition (\ref{eqn:oc2}) together with the  Riesz representation theorem yields the optimality condition  (c.f., \cite{troltzsch2010optimal})
\beq\label{eqn:control_opt_con}
\fv_\eta \Lag(\vp^*,\eta^*,p)=\theta\eta^*+\frac{c_G}{\eps}p=0.
\end{equation}


\begin{Rem}[Choice of the potential]
 We note that our approach to the optimal control problem involving the formulation of the adjoint problem appears to require a smooth potential $G$ (c.f., \eqref{eqn:G_def}). The formulation of the adjoint problem is to our best knowledge an open problem for other widely used, but non smooth or unbounded, potentials such as the obstacle or logarithmic potential.
\end{Rem}
\section{Practical considerations, implementation and algorithm}
\label{sec:implementation}
As is standard, we use the optimality conditions to construct an iterative optimisation loop to solve the optimal control problem, Problem \ref{Pbm:control}. The basic idea is that in each step of the loop we first solve the state equation \eqref{eqn:pf} with a given control, then solve the adjoint equation (\ref{eqn:adjoint}) with the computed states and then update the control using the optimality condition \eqref{eqn:control_opt_con}. For this initial study to ensure robustness of the algorithm and to aid in the clarity of the exposition we employ a simple  gradient based update of the control \cite{snyman2005practical}. 
Given $\eta^k$ and $p$ we compute the updated control $\eta^{k+1}$ via steepest descent.
That is we choose as an update direction the negative gradient, the formula for the update of the control is 
 \begin{equation}\label{eqn:update}
\eta^{k+1}(\vec x,t)=\eta^k(\vec x,t)-\alpha\Big(\theta\eta^k(\vec x,t)+\frac{c_G}{\varepsilon}p(\vec x,t)\Big),\quad(\vec x,t)\in\O\times[0,T),
\end{equation}
where $\alpha$ is a step size. For simplicity in this study we take a constant step size of $\alpha=0.01$.

For the termination criteria for the algorithm we stop if the objective functional $J$ is less than a given tolerance $tol_J$, 
the update in the control is less than a given tolerance, i.e., 
if $\ltwon{\alpha(\theta\eta^k+\frac{c_G}{\varepsilon}p)}{\O\times[0,T)}<tol_\eta$ or if a maximum number of iterations $K_{max}$ is reached. In practice the forward (\ref{eqn:pf}) and adjoint equations (\ref{eqn:adjoint}) must be approximated and to do this we utilise a finite element method, using continuous piecewise linear elements. The details of the numerical method for the approximation of the forward and adjoint equations are given in  \ref{app:FEM}.

 The cell tracking algorithm we propose may now be stated in pseudocode as follows (the notation employed is as in \ref{app:FEM});
\newcommand{\Algoname}[1]{\ensuremath{\text{\textsf{#1}\xspace}}}
\subsection{$\Algoname{Optimal control based cell tracking algorithm}$} 
\begin{algorithmic}\label{alg:control}
\Require

Data: $\vp^0_h,(\vpo)_h$ the initial and target (discrete) diffuse interface data.
 
 Numerical Parameters: $T>0$ end-time and $M>0$ number of timesteps.
 
 Optimisation Parameters: tolerances ${tol_J}$, ${tol_\eta}$, and $K_{max}$. 
 
 Initial guess for the control: $(\eta_h)^0:=(\eta_h^i)^0\in\V,i=0,\dotsc,M.$
\State $k:=0$
\While{$\Big(\ltwon{\alpha(\theta(\eta_h)^k+\frac{c_G}{\varepsilon}p_h^{k+1})}{\O\times[0,T)}>tol_\eta$, $J>tol_J$ and $k<K_{max}\Big)$}\\

\State Solve state equation for $\{(\vp_h^i)^{k+1},(\lambda^i)^{k+1}\},i=1,\dots,M$ with $(\eta_h^i)^k$ and initial data
$(\vp_h^0)^{k+1}=\vp^0_h.$
 \Comment{c.f., \S \ref{app:FEM}}\\
\State Solve adjoint equation for $(p_h^i)^{k+1},i=M-1,\dots,0$ with computed $(\vp_h^i)^{k+1}$ and with terminal data 
$(p_h^M)^{k+1}=(\vp_h^M)^{k+1}-(\vpo)_h.$\\
\State Update control $(\eta_h^i)^{k+1}=(\eta_h^i)^k-\alpha(\theta(\eta_h^i)^k+\frac{c_{G}}{\varepsilon}(p_h^i)^{k+1})\quad i=0,\dots,M.$\\
\State Compute $J$ according to (\ref{eqn:objective}).\\
\State $k:=k+1$.\\
\EndWhile 
\end{algorithmic}

\begin{Rem}[Segmentation and image data]\label{Rem:segmentation}
An important aspect of any cell tracking algorithm is its ability to extract data suitable for the tracking algorithm from the experimental image data set. In many cases the experimental image data set is grayscale data with the intensity (brightness) indicating whether a point is in the interior of a cell, i.e., points inside the cell appear bright for example and points outside appear dark. For many tracking algorithms this intensity data is then post processed via a segmentation algorithm (e.g., active contour methods \cite{chan2001active,dormann2002simultaneous}) to yield sharp interface representations of the cell membrane.  Assuming a sharp interface representation of the cell membrane is available, diffuse interface representations may be easily initialised (see for example \cite{2013arXiv1311.7602C}).  We note however that the raw intensity data produced by many imaging procedures may already be close to a diffuse interface representation of the cell, this is typically the case when the data is relatively free of noise and the contrast between the cell and the background is high. In this case one may wish to exploit this fact in the algorithm and work with the raw image data set itself (or a post processed e.g., thresholded version), thus circumventing the extra error induced by segmentation.  
\end{Rem}

\begin{Rem}[Observations at multiple points in time]\label{Rem:multi-obs}
For clarity of exposition we focus on the case of fitting to a single observation. The approach generalises straightforwardly to multiple observations with the first term in (\ref{eqn:objective}) simply replaced by a sum over the distinct times at which the observations are taken of the difference between the solution (at the appropriate time) and the target data.  
\end{Rem}

\begin{Rem}[Multiple cells and matching problems]\label{Rem:multi-cells}
As mentioned above a major focus of many cell tracking algorithms is to track multiple cells in the same image and the resolution of the so called matching problem. Our approach can be applied to multi-cell image data. Here $\vp^0$ and $\vpo$ would be diffuse interface representations of the multi-cell image data set and the diffuse interfaces would consist of multiple disjoint phases. The remaining aspects of the  approach remain unchanged and the matching problem is solved implicitly in the computation of the optimal control. 

There are however multiple practical issues which arise in this setting related to the separation between distinct cells, which affects the choice of $\eps$, and the fact that the evolution law (\ref{eqn:pf}) allows changes in the topology of the phases which may lead to cell splitting, the annihilation of a phase (which would correspond to the disappearance of a cell) or the nucleation of a phase (i.e., the spontaneous appearance of a cell). We intend to comment on practical approaches to multi-cell tracking elsewhere.
\end{Rem}

\section{Numerical examples}\label{sec:application}

We now present some benchmark numerical examples illustrating the application of the algorithm to artificial image data sets. For all the simulations we report on in this section, in the state equation (\ref{eqn:pf}) we set $\eps=0.1$, and we took the end-time $T=0.4$. As mentioned previously the parameter $\eps$ governs the width of the diffuse interface and should in general be taken as small as is computationally feasible, smaller values of $\eps$ necessitate a finer grid, in this initial study with uniform grids we set $\eps=0.1$ as the CPU times become prohibitive for smaller values of $\eps$.  The end time $T$ corresponds to the nondimensionalised time between snapshots and could in principle be related to an acquisition time between images given real biological data. For each of the experiments, apart from those of \S \ref{subsec:initial-guess},  we set the initial guess for the control  to be  constant in space and time (zero in the single cell case  and one for the multi-cell examples).
For the approximation of the forward and adjoint equations we used triangulations with $8321$ DOFs in all the simulations, apart from those of \S \ref{subsec:convergence-test},  and selected a uniform timestep $\tau=1\times10^{-3}$. The same numerical parameters for the optimisation algorithm were used for all the experiments and  are given in Table \ref{table1}. In every example we report on the algorithm terminated due to the update of the control being less than the prescribed tolerance.
\begin{table}[H]
\centering
\begin{tabular}{c c c c c}
\hline
$\alpha$ &$\theta$ &${tol_J}$&${tol_\eta}$&$K_{max}$\\
\hline
0.01&0.01&$1\times10^{-4}$&$1\times10^{-4}$&3500\\
 \hline
\end{tabular}
\caption{Parameter values used for the numerical simulations.}
\label{table1}
\end{table}
The technical details of the hardware used to carry out the simulations is given in Remark \ref{rem:hardware}.

\begin{Rem}[Hardware details]\label{rem:hardware}
 All the numerical experiments have been performed on the high performance cluster (HPC) at the University of Sussex. Each of the simulations was carried out in serial using a single core of the cluster.  The HPC cluster currently consists of 3140 cores with an even mixture of Intel and AMD CPUs. The majority of the cluster are 64 core AMD nodes with 256GB RAM per node, and a smaller number of 512GB RAM nodes. The cluster uses the high-performance Lustre clustered-filesystem for I/O, and currently stands at 298TB of storage for research use.
\end{Rem}

\subsection{Application to synthetic data}
\label{subsec:results-single}
Here we apply the algorithm to a single synthetic cell data set  taken from the {\it{PhagoSight}} website \url{http://www.phagosight.org/synData.php}. 
The synthetic cell was generated as a mixture of Gaussians with Poisson noise that varied over time to simulate the displacement and change of shape of a neutrophil as observed in a Zebrafish embryo. The data for analysis consisted of points on the synthetic cell membrane at a series of times (for simplicity we used $2d$ data, i.e., the cell membrane was a $1d$ curve embedded in $\Reals^2$). The initial and target curves we took as test data for the algorithm  are shown in Figure \ref{fig:synthetic_data}.   To apply our algorithm, based on diffuse interface representations, we define the domain $\O$:=$[0,8]\times[0,6]$ which was such that both the initial and target curves were contained in the domain. We then constructed  diffuse interface representations of the target data following the procedure described in  \cite{2013arXiv1311.7602C}. Figures \ref{fig:synthetic_data_PF_initial} and \ref{fig:synthetic_data_PF_target2} show the diffuse interface representations of the initial and target data respectively. 
\begin{figure}[h!]
\begin{center}
\subfigure[Initial (red curve) and target (green curve) synthetic data. The cell centroids are shown together with the trajectory of the linear interpolant of the cell centroids (black line).]{\label{fig:synthetic_data}
{\includegraphics[ width=0.6\linewidth]{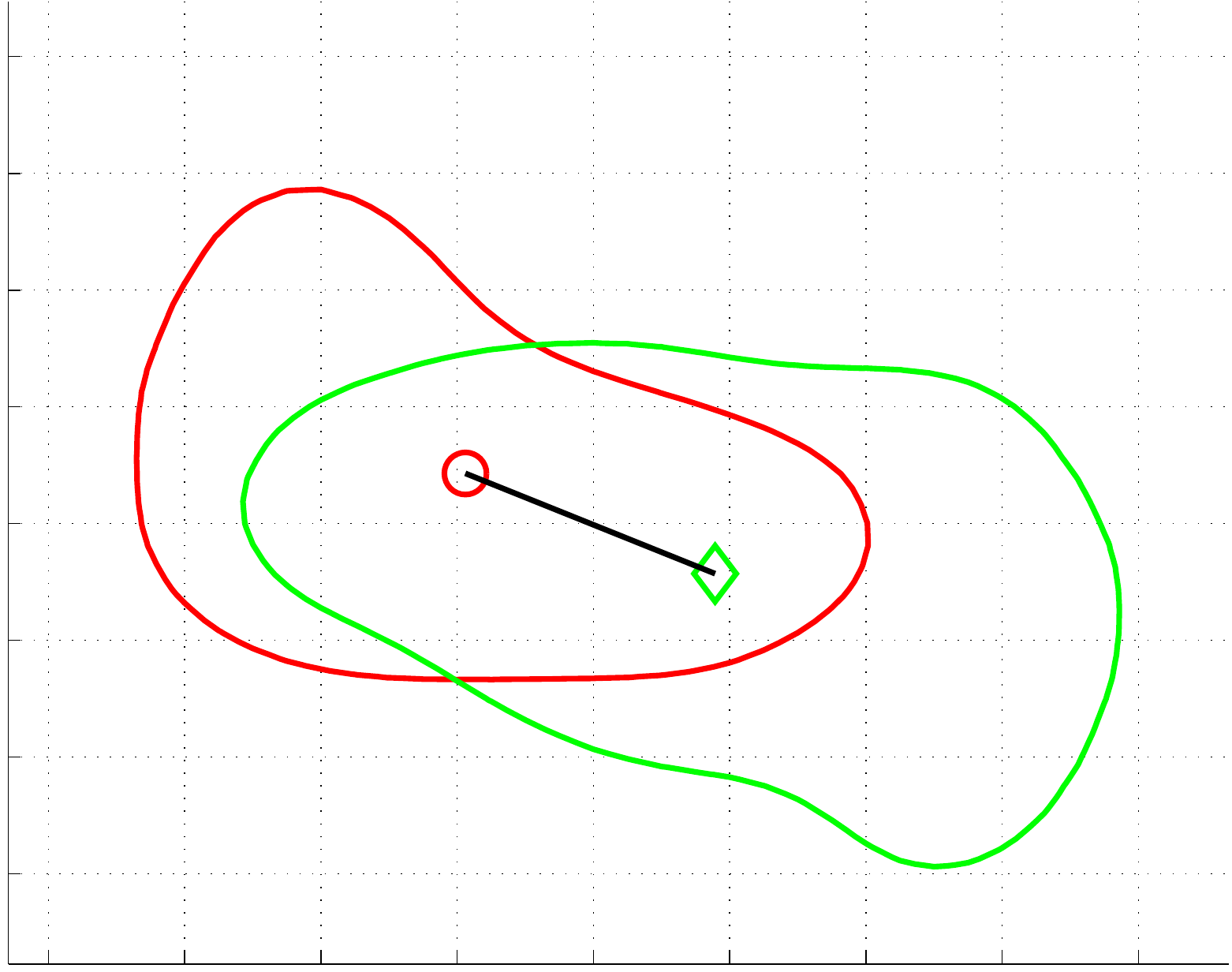}}}\\
\subfigure[Initial data ($\vp^0$).]{\label{fig:synthetic_data_PF_initial}
{\includegraphics[trim = 10mm 140mm 20mm 20mm, clip,width=0.45\linewidth]{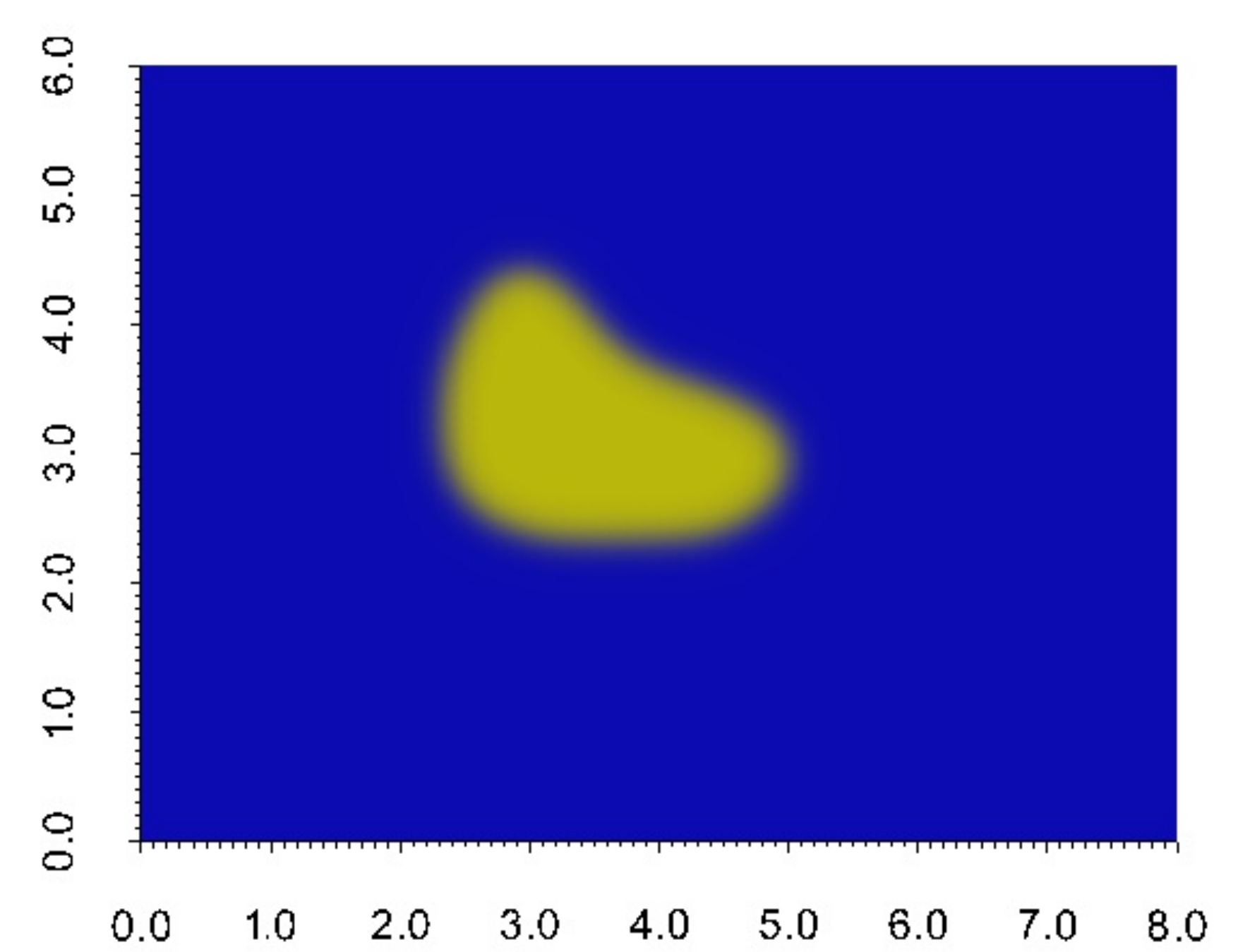}}}
\subfigure[Target data ($\vpo$).]{\label{fig:synthetic_data_PF_target2}
{\includegraphics[trim = 10mm 140mm 20mm 20mm,clip, width=0.45\linewidth]{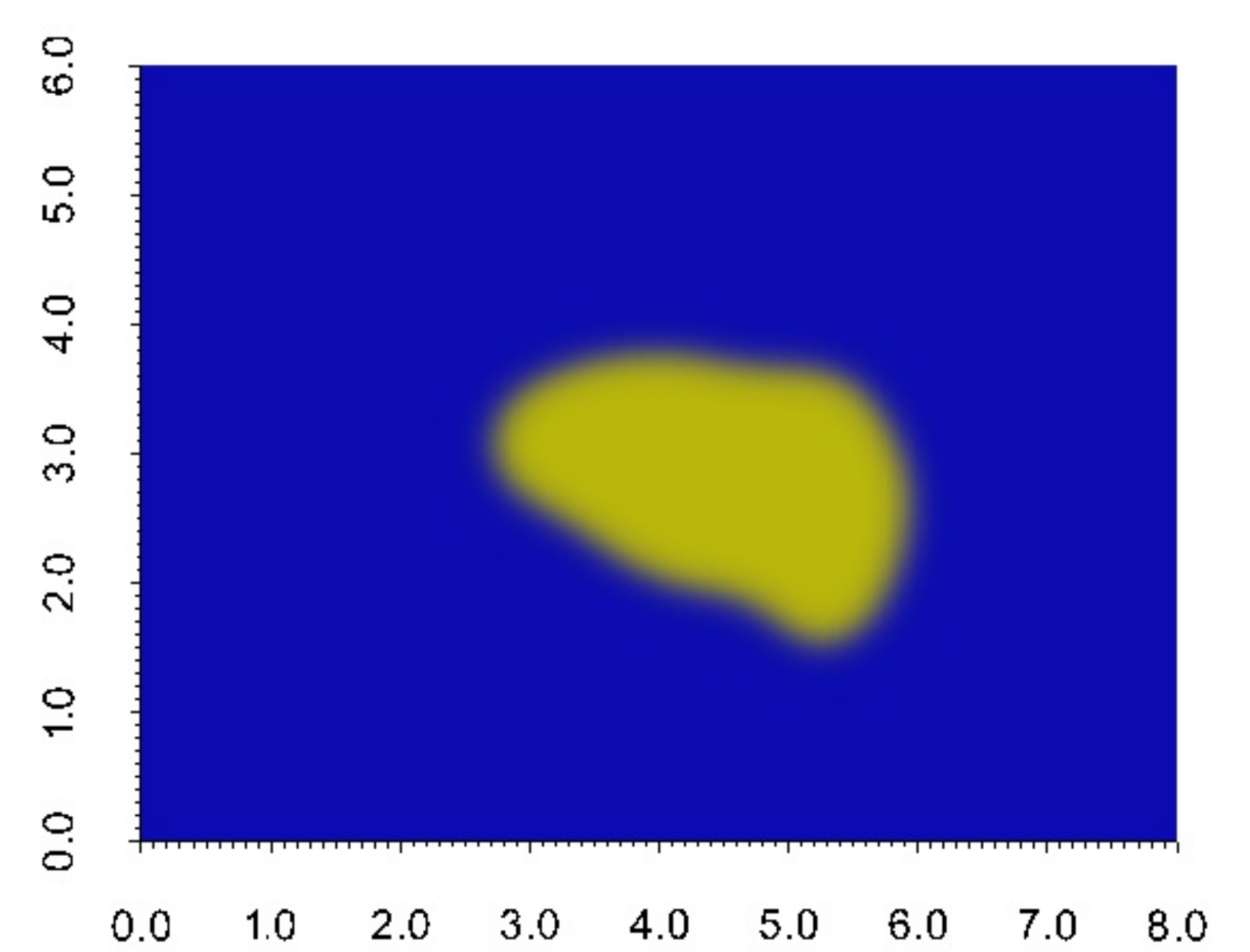}}}
\end{center}
 \caption{Initial and target data for the example of \S \ref{subsec:results-single}.}
   \end{figure}
 In order to investigate the influence of the volume constraint on the computed cell morphologies we performed two experiments, in the first we simply considered the forced Allen-Cahn
 model for the evolution with no volume constraint,  i.e., \eqref{eqn:pf} with $\lambda =0$, and in the second we included the volume constraint as described in \S \ref{sec:prob}.
 The algorithm took 1996 iterations to meet the stopping criteria with no volume constraint  and 2479 iterations with the volume constraint, corresponding to CPU times of  25238 and 119216 seconds respectively.
 
 Figure \ref{fig:cost_single} shows the value of the objective functional against the number of iterations of the optimisation algorithm with and without the volume constraint. We  observe similar behaviour in both cases with an initial rapid decrease in the objective functional followed by a more gradual reduction with each iteration as we approach the minimum. Figure \ref{fig:single_optimal_target} shows the zero level-set of the computed solution using the optimal control at the final time with and without the volume constraint. The curve corresponding to the zero level set is shaded by the value of the computed optimal control. The background shading corresponds to the target data. In both cases the position of the zero level-set of the computed solution shows good agreement with the target data. Qualitatively we observe cells with a clearly  defined ``front'' and ``rear'', with the computed control corresponding to protrusive forces at the front and contractive forces at the rear. 
     \begin{figure}[h!]
\begin{center}
\subfigure[Without the volume constraint]
{{\includegraphics[trim = 10mm 0mm 20mm 0mm,clip,width=0.49\linewidth]{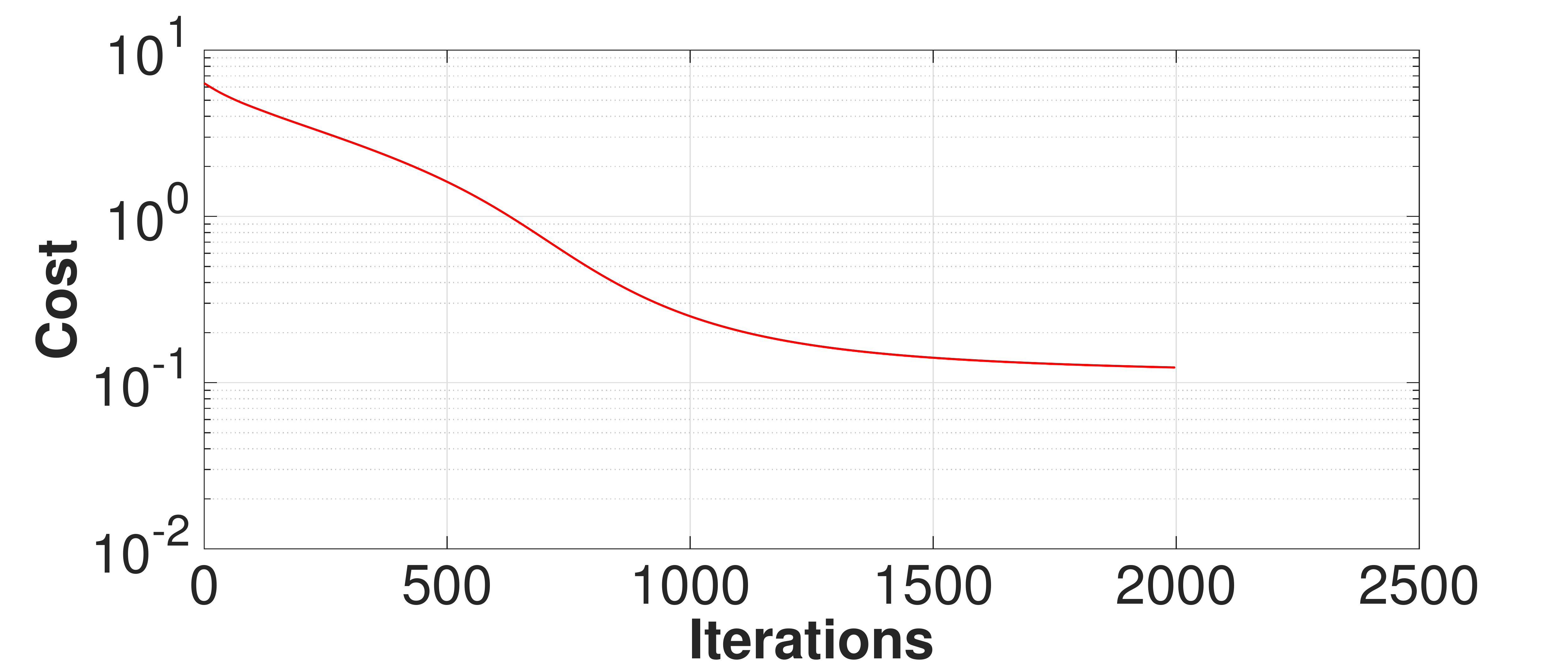}}}
\subfigure[With the volume constraint]
{{\includegraphics[trim = 10mm 0mm 20mm 0mm,clip,width=0.49\linewidth]{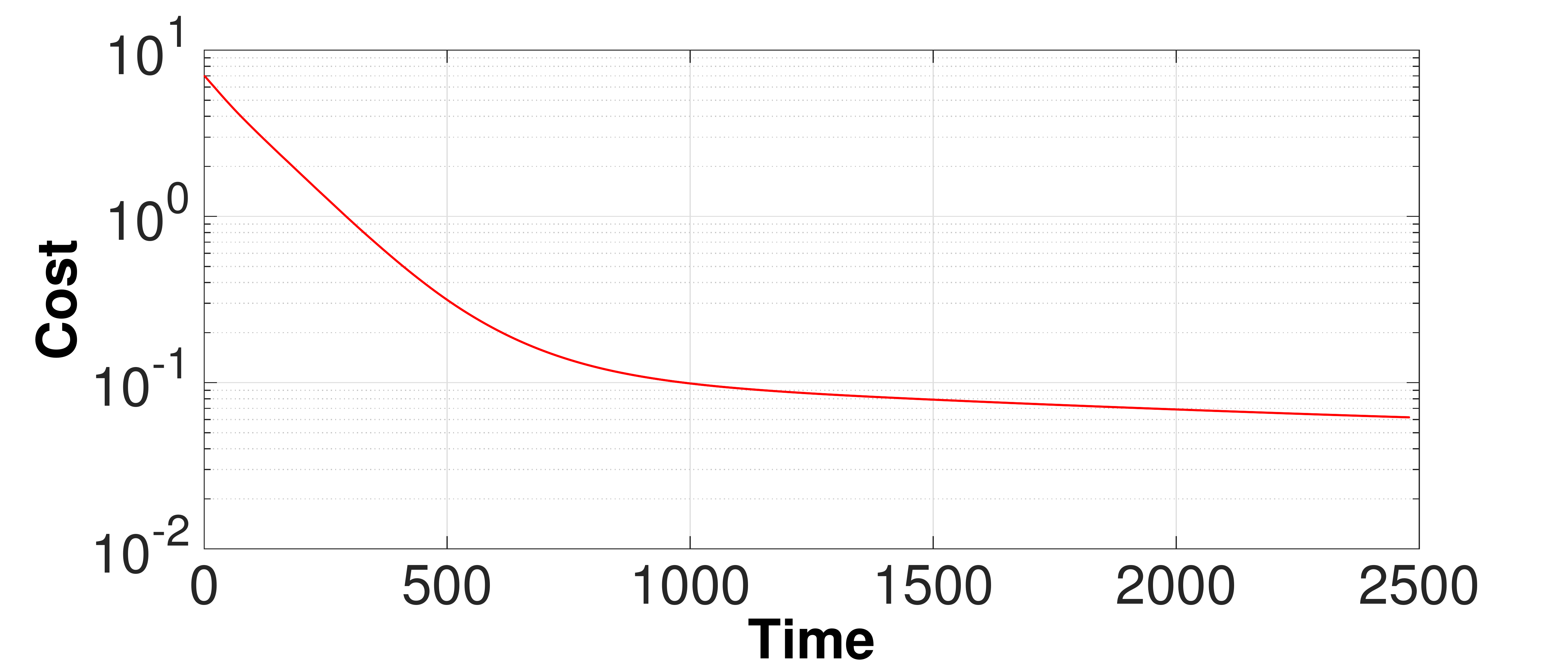}}}
\end{center}
 \caption{The value of the cost functional versus the number of iterations for the experiments of \S \ref{subsec:results-single} with and without the volume constraint. We observe a rapid decrease in the cost initially followed by a much more gradual  decrease as we approach the minimum, this is as expected since the steepest descent algorithm is used for the update of the control.}\label{fig:cost_single}
  \end{figure}
  \begin{figure}[h!]
\begin{center}
\subfigure[Without the volume constraint]
{{\includegraphics[trim = 10mm 140mm 20mm 20mm,clip,  width=0.49\linewidth]{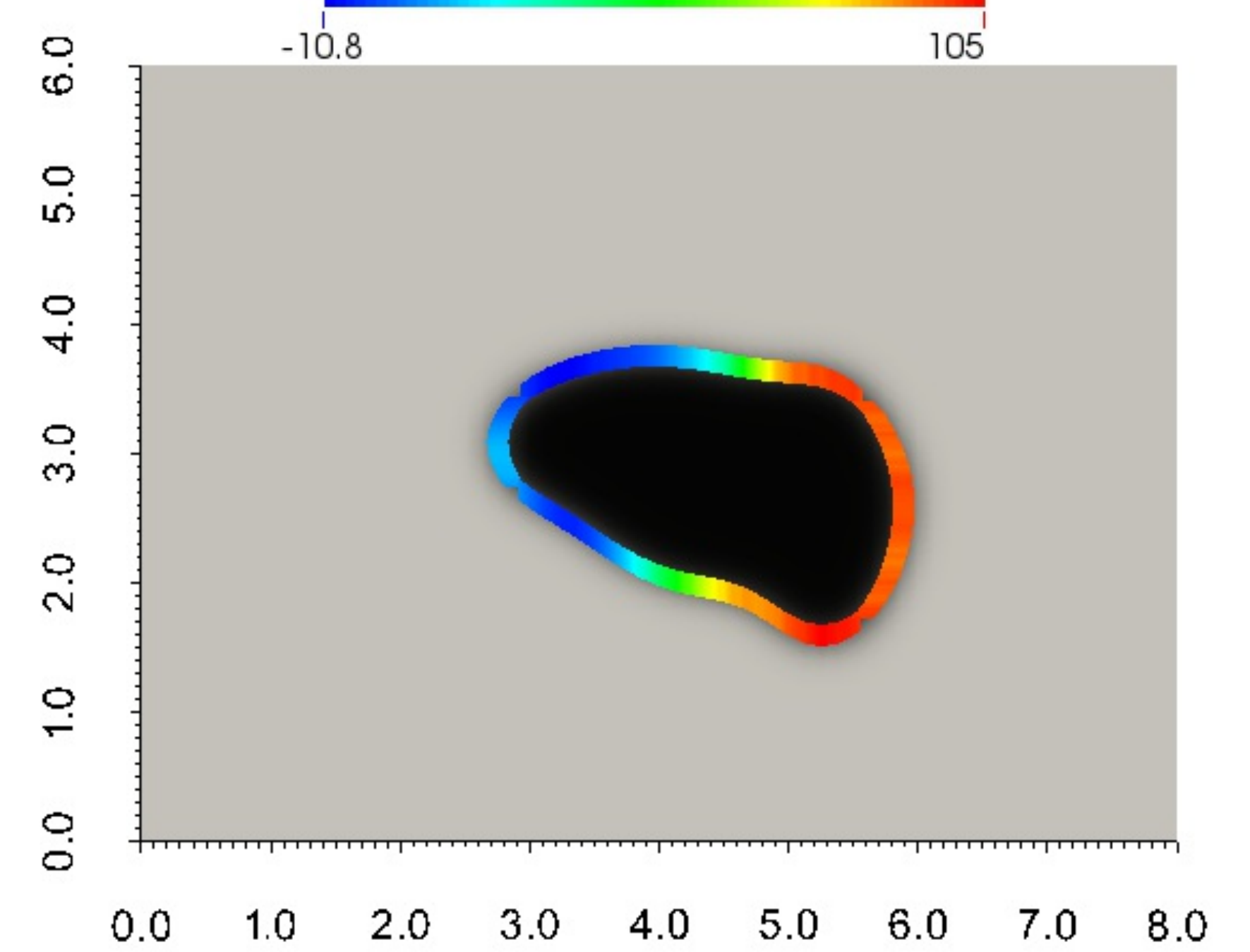}}}
\subfigure[With the volume constraint]
{{\includegraphics[trim = 10mm 140mm 20mm 20mm,clip,  width=0.49\linewidth]{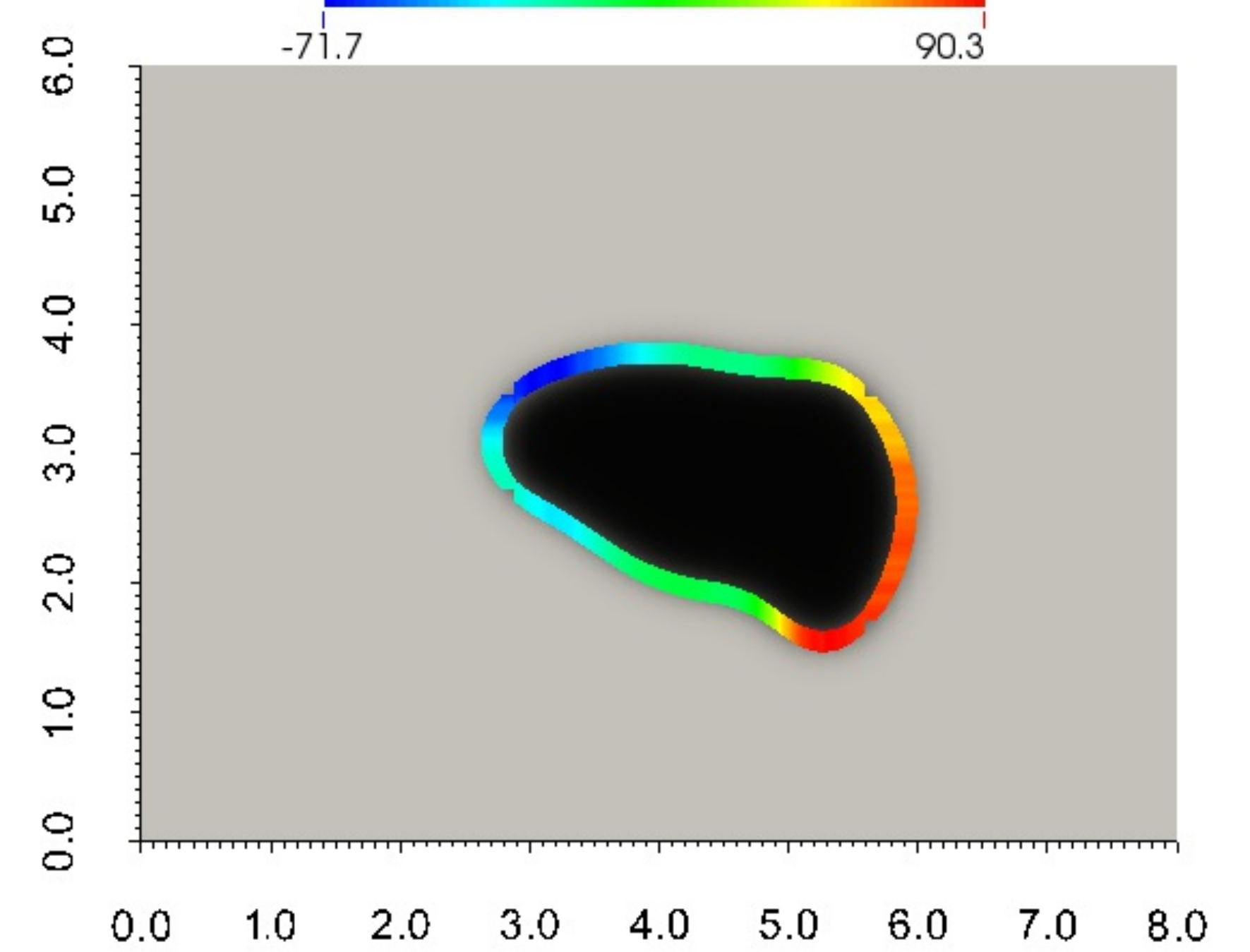}}}
\end{center}
 \caption{
 Zero level-set of the solutions ($\vp(\vec x,T)$) computed using the approximated optimal control ($\eta^*$) with and without the volume constraint for the experiments of \S\ref{subsec:results-single}.
 The curve (zero level-set of $\vp(\vec x,T)$)  is shaded by the approximated optimal control ($\eta^*(\vec x,T)$) and the background by the target data ($\vpo(\vec x)$). The color-bar corresponds to the scale for $\eta^*(\vec x,T)$.  We see good agreement between the zero level-set of the data computed with the optimal control and the target data in both cases.
 }\label{fig:single_optimal_target}
  \end{figure}
  
 Figure \ref{fig:area_single} shows the area enclosed by the zero level-set of the solution with the optimal control with and without the volume constraint together with the linear interpolant of the areas of the data. We see that without the volume constraint the area initially decreases then rapidly increases as we approach the final time whilst with the volume constraint (note the constraint is actually on the mass rather than the volume) the area is close to the  linear interpolant of the areas of the data. In terms of the computed  cell morphologies, Figure \ref{fig:single_optimal_curves} shows snapshots of the computed cell membranes (zero level-sets) for the two different cases. We clearly observe that the intermediate snapshot (blue curve) encloses a much smaller area if the volume constraint is not included in the algorithm. In Figure \ref{fig:centroids_single} we report on the trajectory and speed (magnitude of the velocity) of the centroid (center of mass) of the zero level-set of the computed solution with the optimal control, with and without the volume constraint. We observe similar trajectories with and without the volume constraint and in both cases we observe an increase in the speed as we approach the end-time. However in the case of no volume constraint this increase is more marked with a sharp spike in the centroid velocity observed close to the final time. The increasing centroid speed we observe may be unphysical and if a (roughly) constant centroid velocity is desired one strategy may be to impose pointwise constraints on the control, this would prevent the large increase in the maximum and minimum values of the control observed during the simulations as we approach the final time as shown in Figure \ref{fig:control_single}. Another possible strategy would be to modify the regularisation in (\ref{eqn:objective}).
    \begin{figure}[h!]
\begin{center}
{{\includegraphics[trim = 40mm 0mm 30mm 0mm,clip,width=\linewidth]{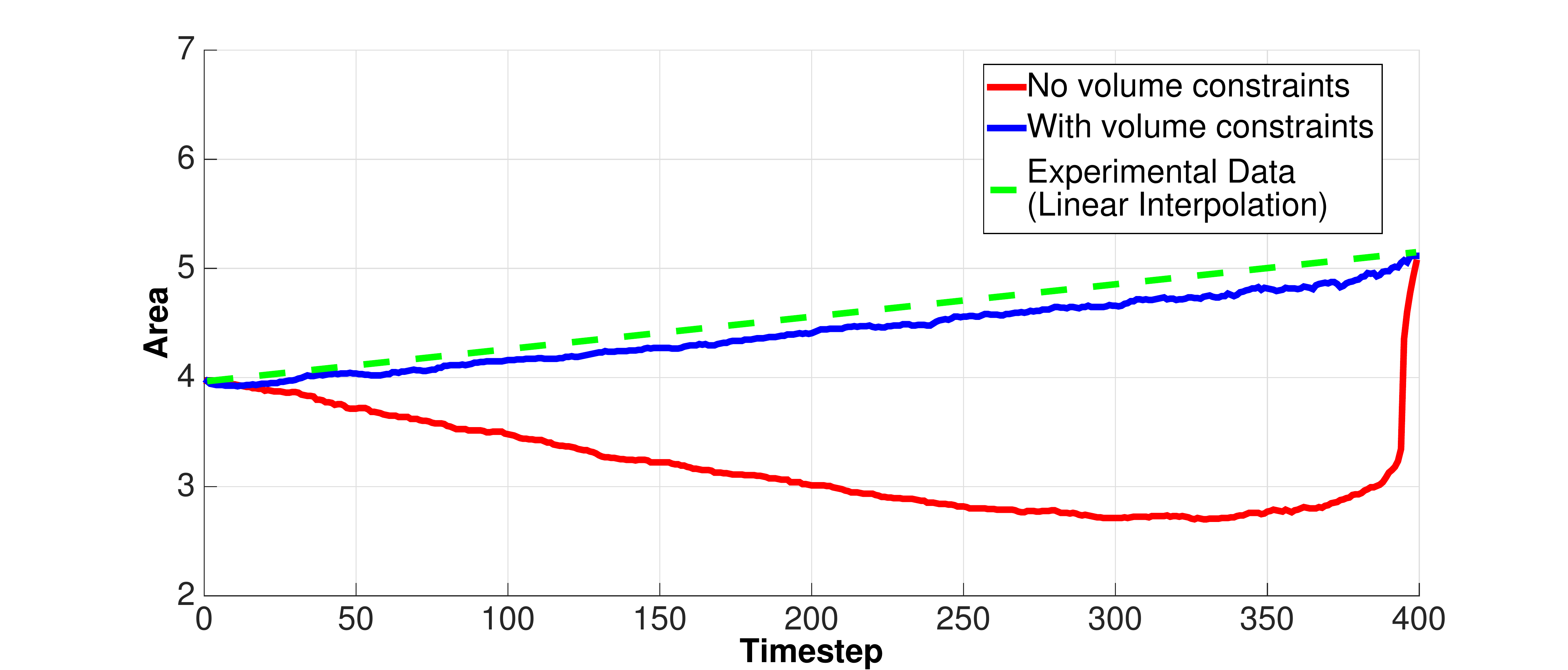}}}
\end{center}
 \caption{Area enclosed by the cell for the experiments of \S \ref{subsec:results-single} with and without the volume constraint. The cell shrinks considerably during the evolution without the volume constraint whilst a good fit to the linear interpolant of the area enclosed by the data is observed with the volume constraint.}\label{fig:area_single}
  \end{figure}
\begin{figure}[h!]
\begin{center}
\subfigure[Without the volume constraint]
{{\includegraphics[trim = 40mm 160mm 40mm 30mm,clip,  width=0.49\linewidth]{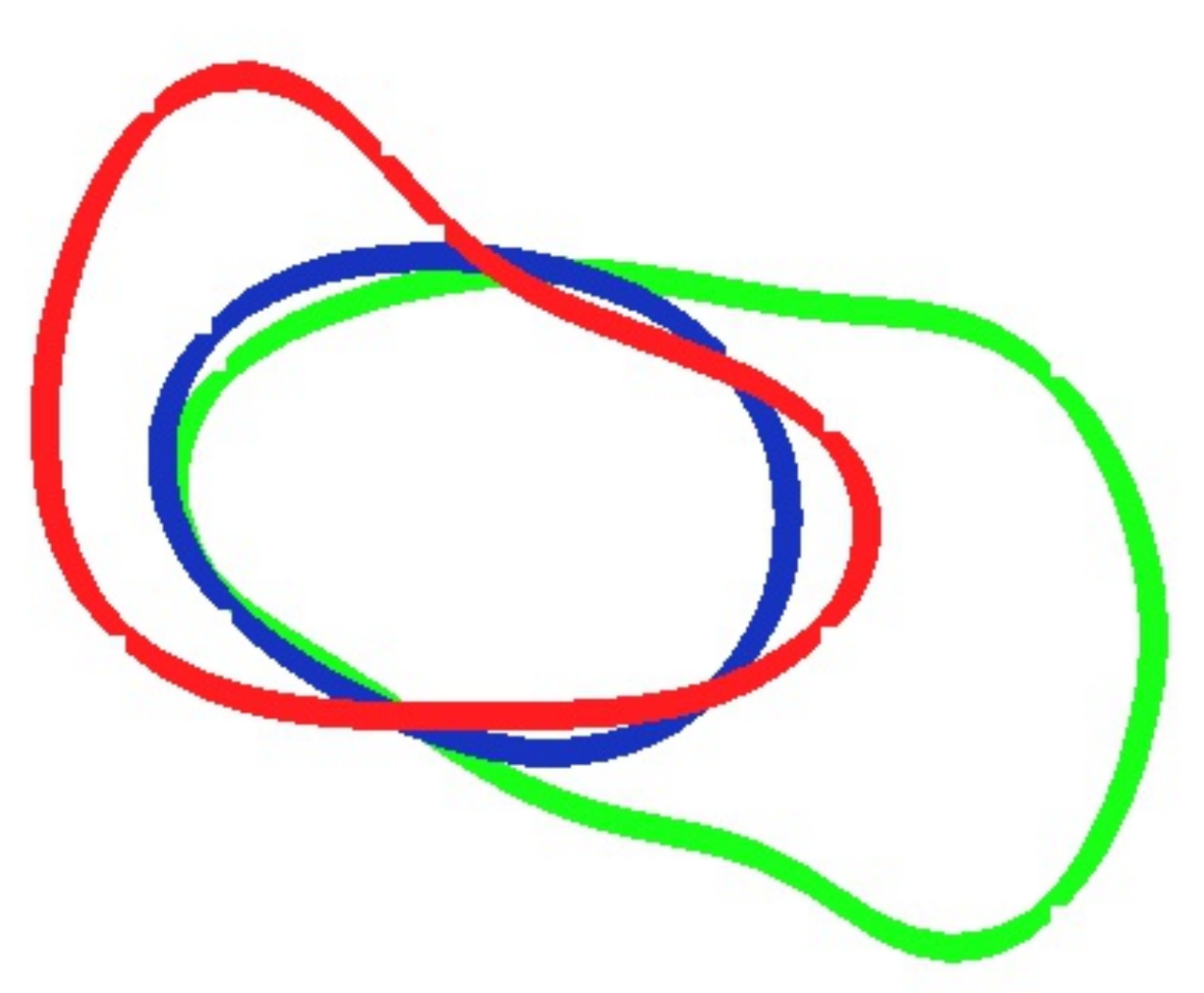}}}
\subfigure[With the volume constraint]
{{\includegraphics[trim = 40mm 160mm 40mm 30mm,clip,  width=0.49\linewidth]{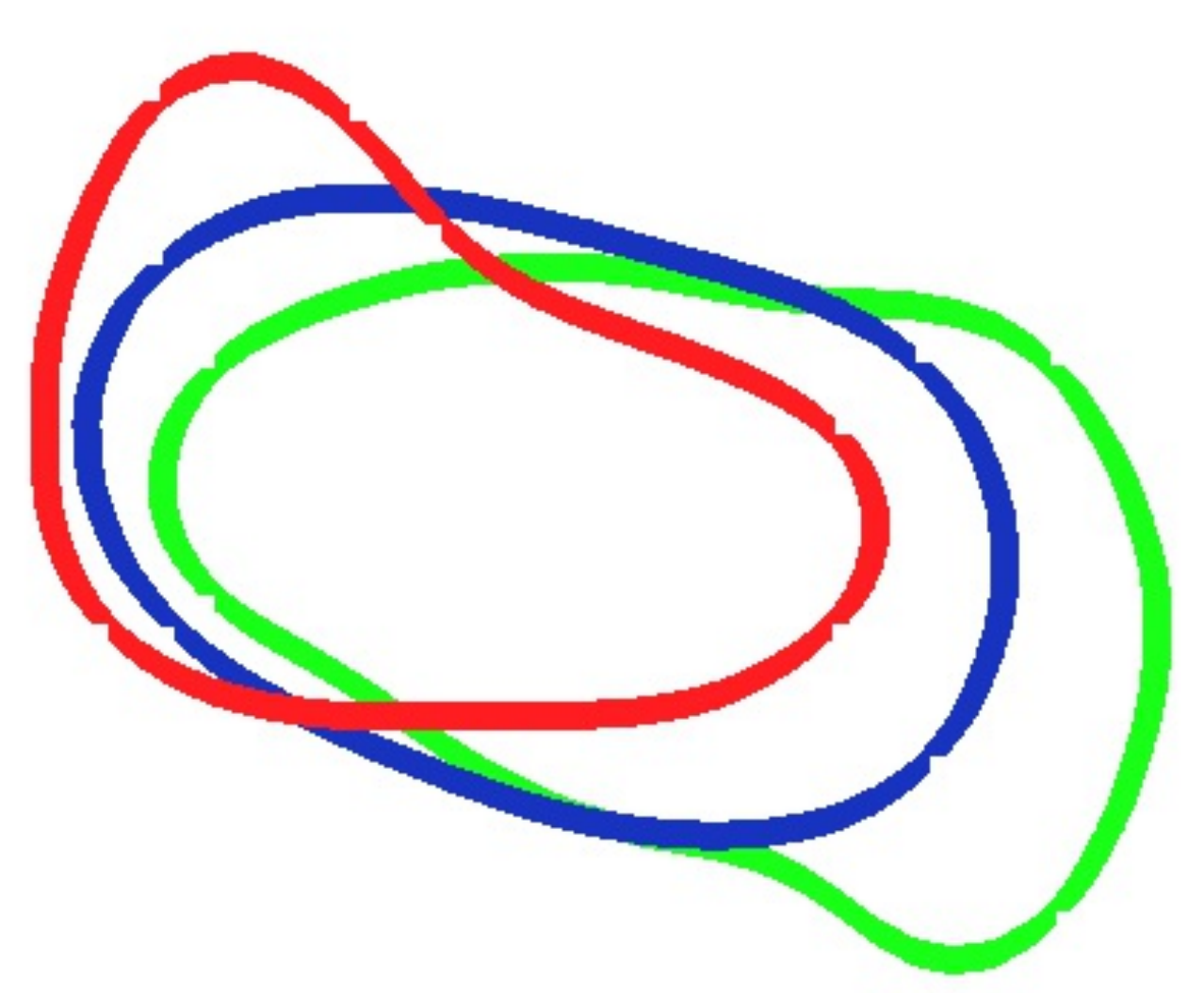}}}
\end{center}
 \caption{Zero level-sets of the solutions computed ($\vp(\vec x,t)$) with the optimal control ($\eta^*(\vec x,t)$) for the experiments of \S \ref{subsec:results-single} with and without the volume constraint at $t=0$ (red), $t=0.35$ (blue) and $t=0.4$ (green). We observe that the volume enclosed by the blue curve is significantly smaller than the volumes enclosed by the red and green curves without the volume constraint whilst this is not observed if the volume constraint is included.}\label{fig:single_optimal_curves}
  \end{figure}
\begin{figure}[h!]
\begin{center}
\subfigure[Centroid trajectories together with centroids of the data.]{{\includegraphics[ width=\linewidth]{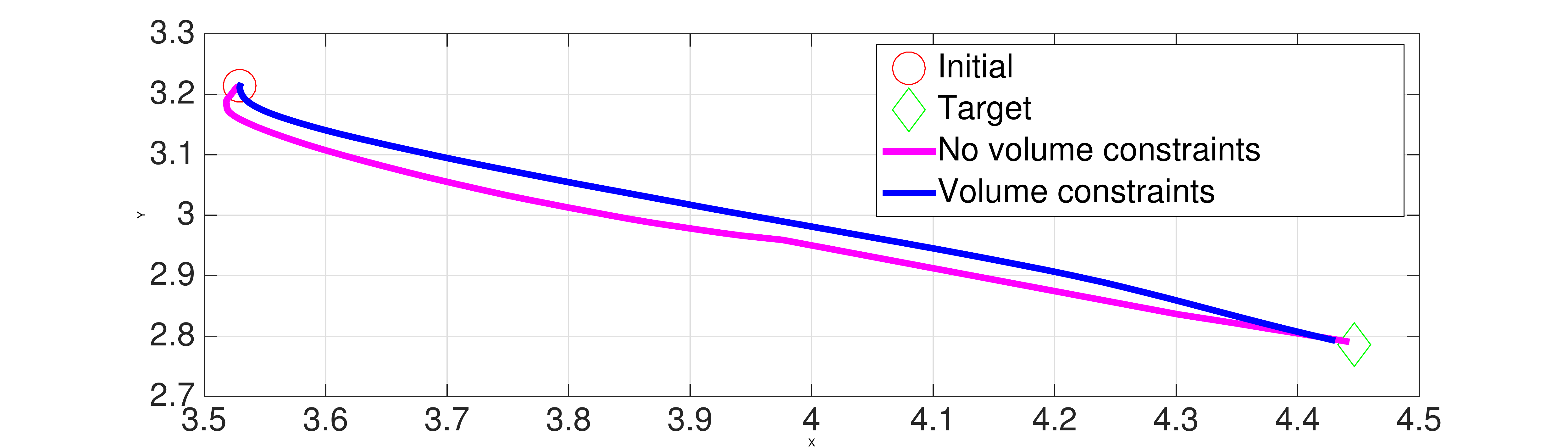}}}\\
\subfigure[Centroid speed (log scaling on the $y-$axis).]{{\includegraphics[ width=\linewidth]{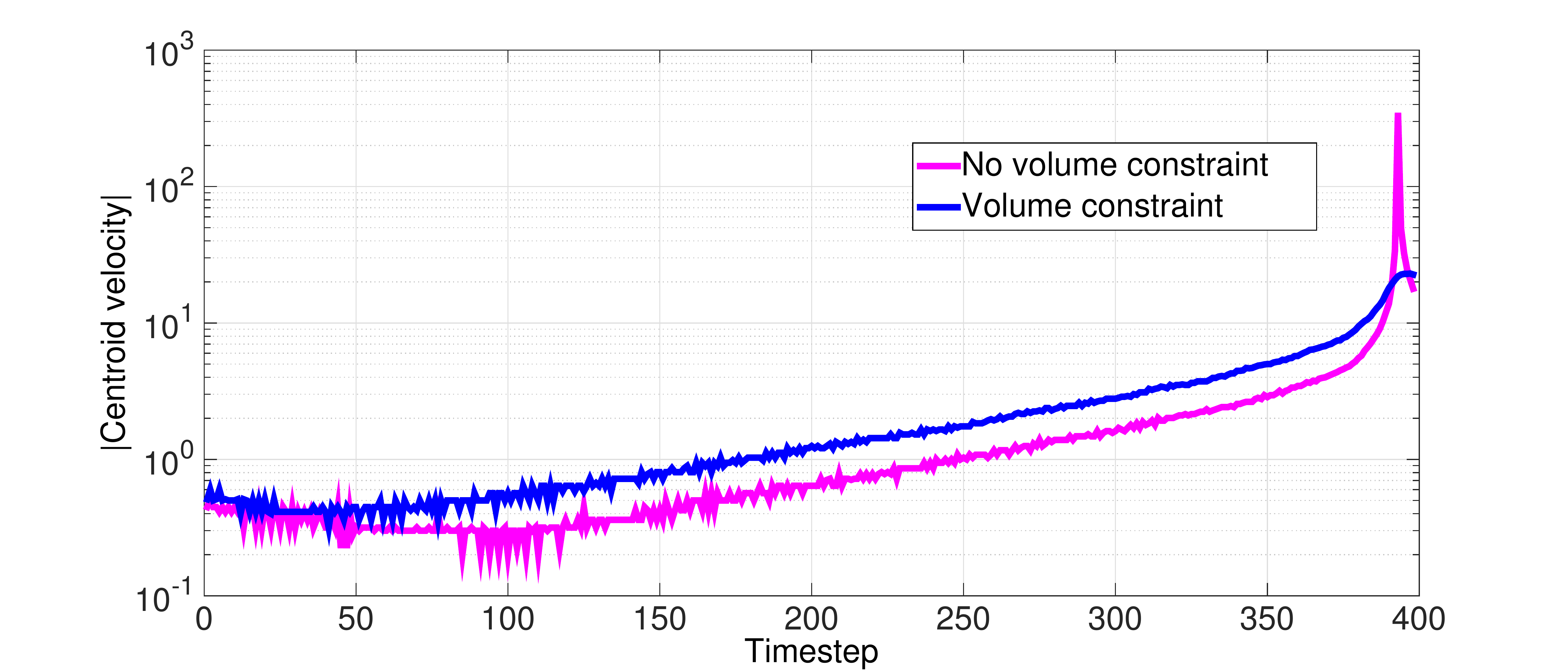}}}
\end{center}
 \caption{ Trajectories and speeds of the centroid of the zero level-sets of the solution with the optimal control with and without the volume constraint for the example of  \S \ref{subsec:results-single}.}
 \label{fig:centroids_single}
   \end{figure}
\begin{figure}[hp!]
\begin{center}
\subfigure[Minimum value of the control for a series of iterations against time for the algorithm without volume constraints.]{{\includegraphics[ width=0.55\linewidth]{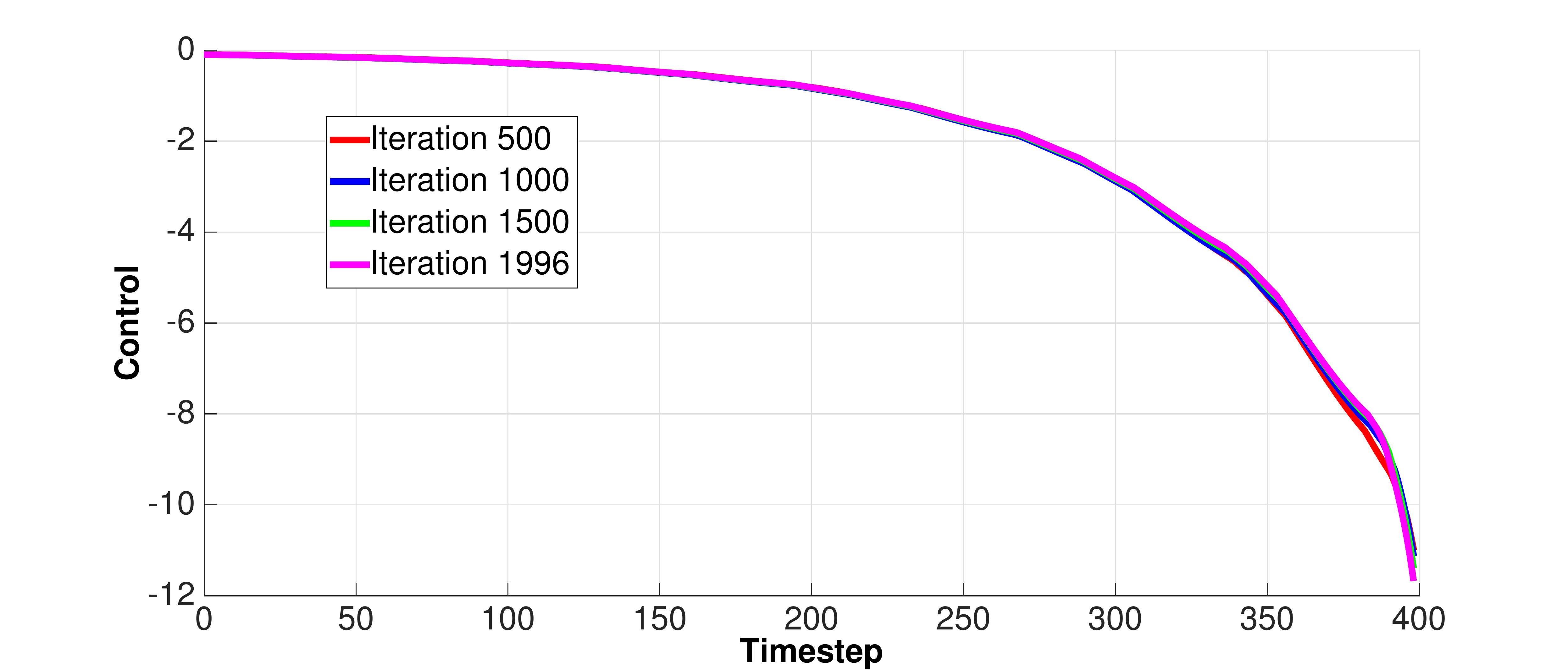}}}
\subfigure[Maximum value of the control for a series of iterations against time for the algorithm without volume constraints.]{{\includegraphics[ width=0.55\linewidth]{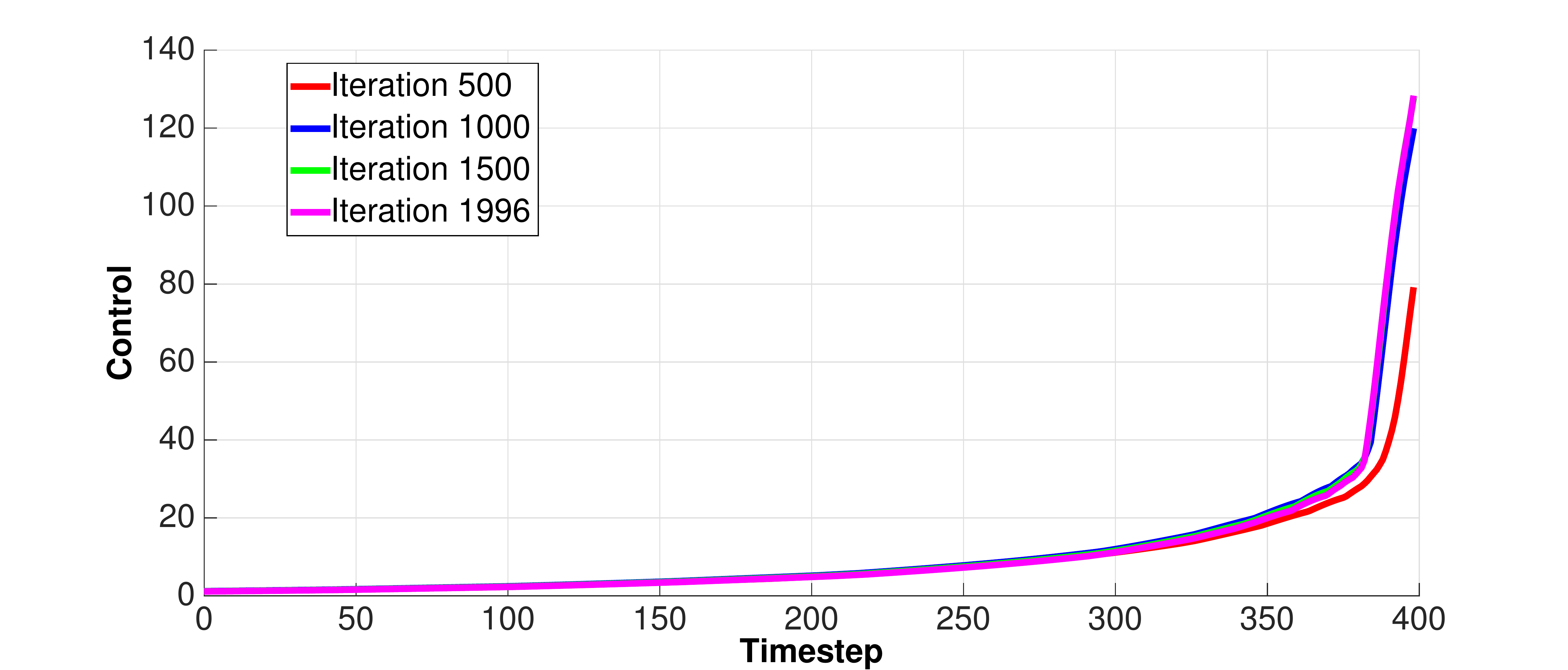}}}
\subfigure[Minimum value of the control for a series of iterations against time for the algorithm with volume constraints.]{{\includegraphics[ width=0.55\linewidth]{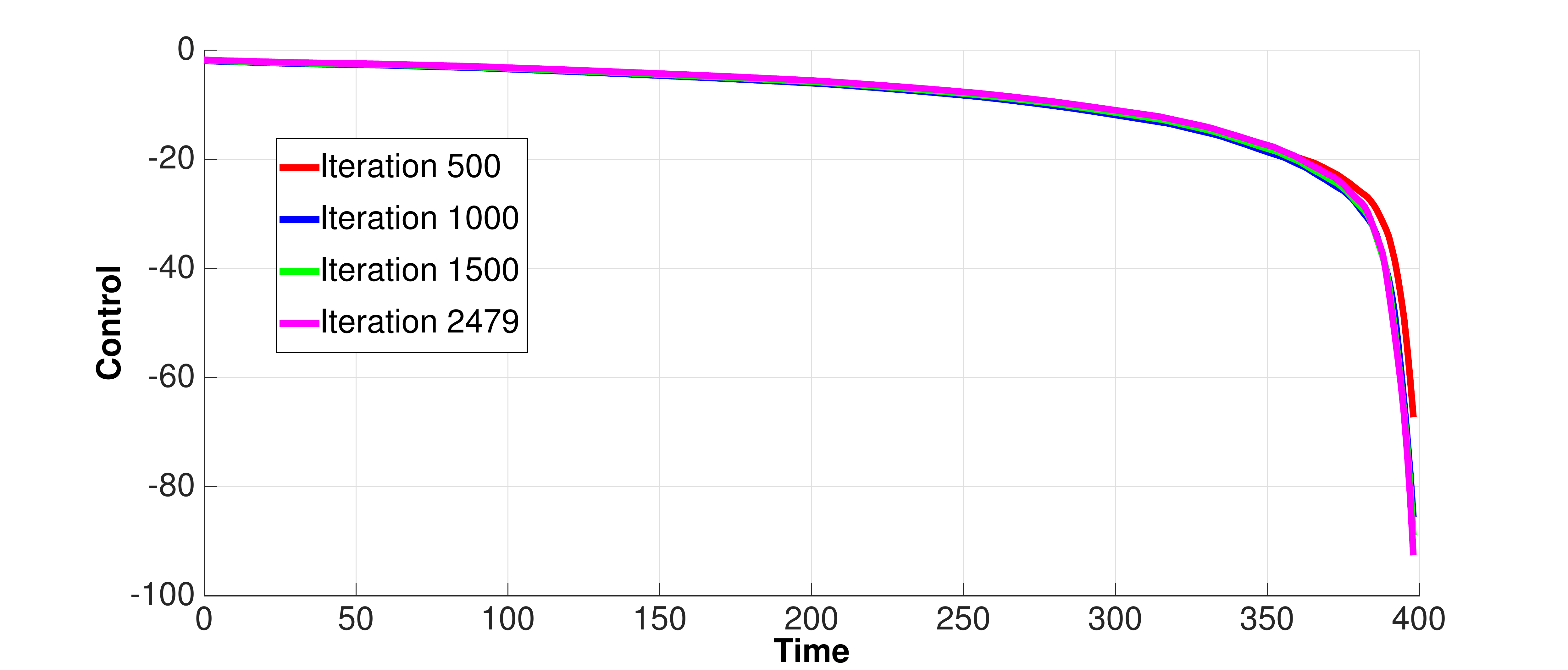}}}
\hskip1em
\subfigure[Maximum value of the control for a series of iterations against time for the algorithm with volume constraints.]{{\includegraphics[ width=0.55\linewidth]{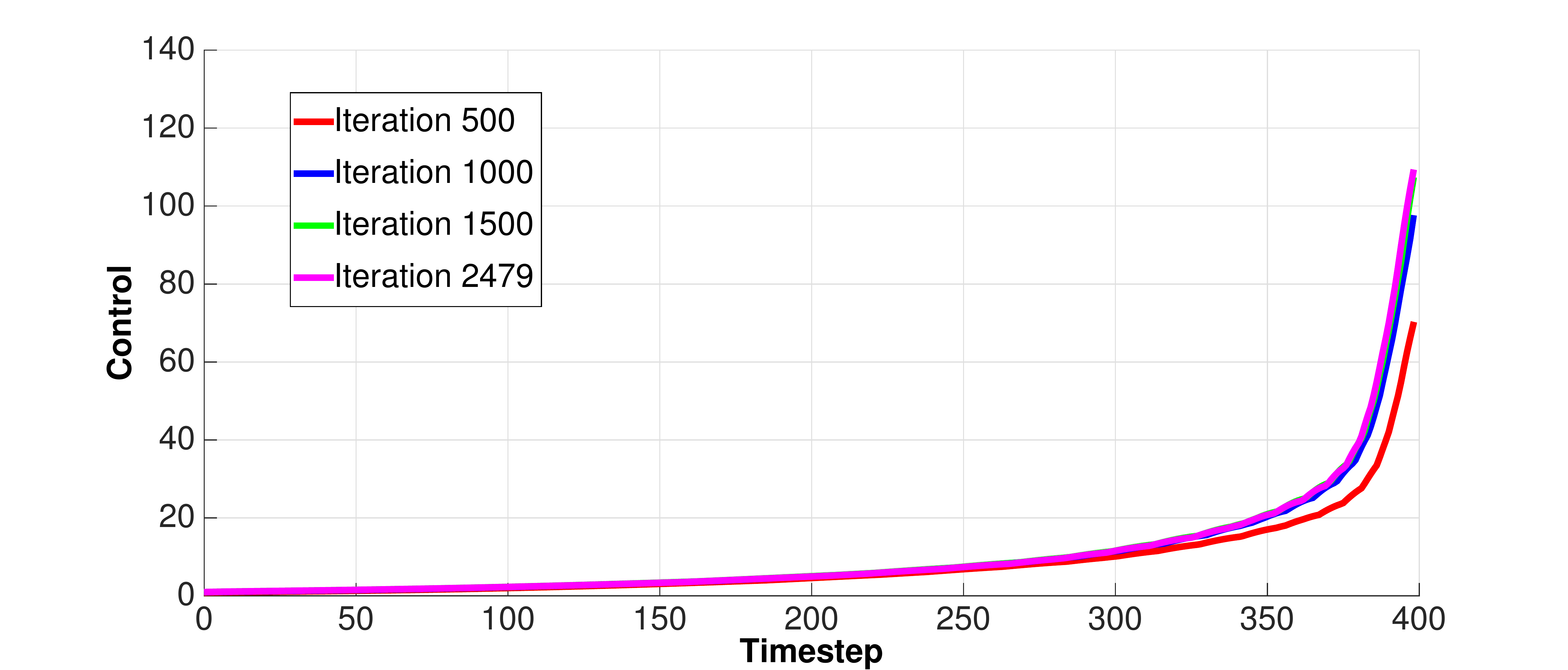}}}
\end{center}
 \caption{ Minimum and maximum values of the control $\eta$ with and without the volume constraint for the example of  \S \ref{subsec:results-single}.}
 \label{fig:control_single}
   \end{figure}
 
Figure \ref{fig:fid_term} shows the fidelity term $\ltwon{\vpo(\vec x)-\vp(\vec x,T)^k}{\O}$ with and without the volume constraint versus the number of iterations (where $k$ corresponds to the optimisation iteration number). The fidelity term may be considered as a quantitative measure for the ``goodness of fit" of the computed data to the observations.
We observe a steady decay in the fidelity term  as we approach the optimal control in both cases. 

     \begin{figure}[h!]
\begin{center}
\subfigure[Without the volume constraint]
{{\includegraphics[trim = 10mm 0mm 20mm 0mm,clip,width=0.49\linewidth]{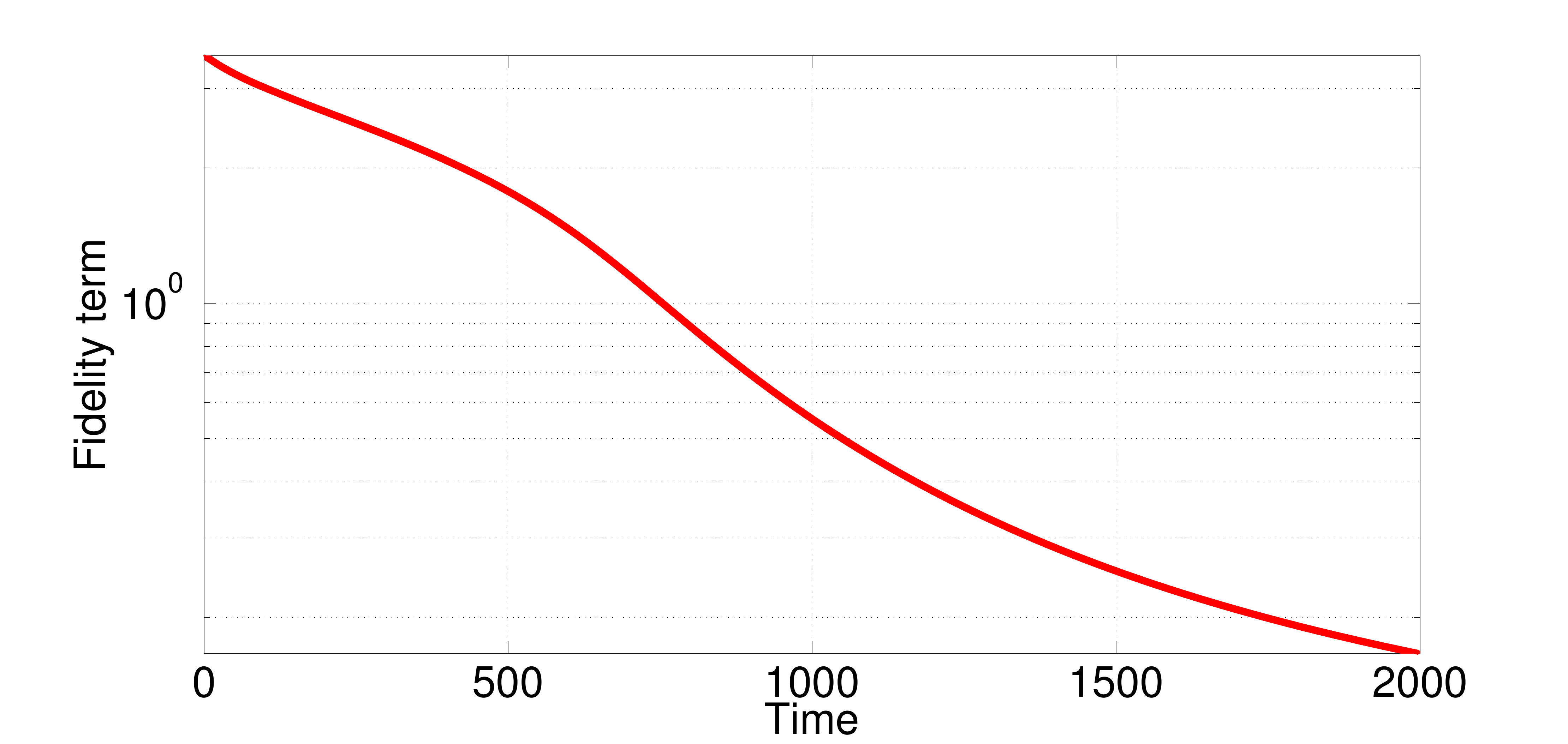}}}
\subfigure[With the volume constraint]
{{\includegraphics[trim = 10mm 0mm 20mm 0mm,clip,width=0.49\linewidth]{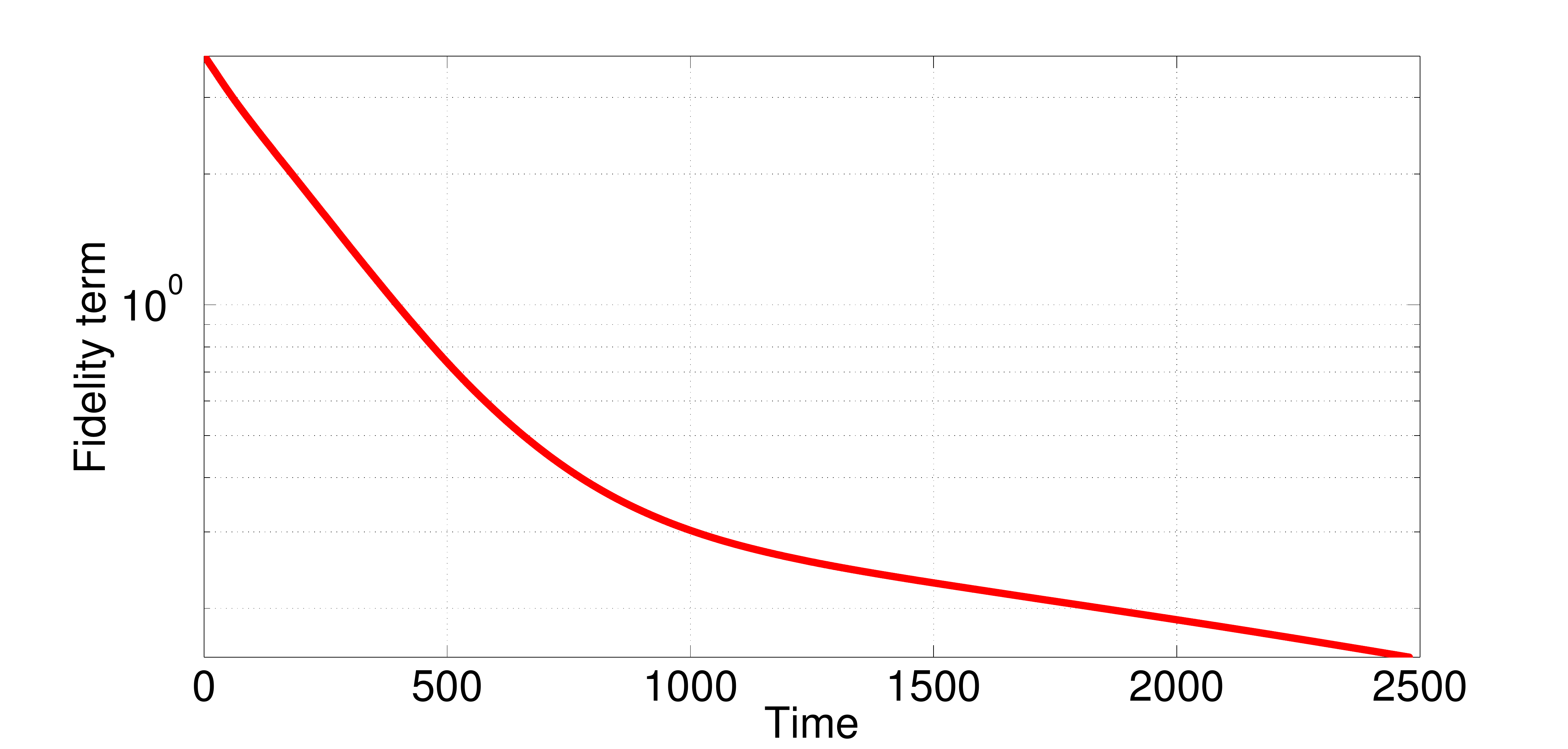}}}
\end{center}

 \caption
 {The value of the fidelity term versus the number of iterations for the experiments of \S \ref{subsec:results-single} with and without the volume constraint.}\label{fig:fid_term}
  \end{figure}
 
 \subsection{The effect of mesh refinement}
\label{subsec:convergence-test}
In this section we investigate the effect of the mesh-size on the results by refining the mesh whilst keeping the time step $\tau$ constant. We report on the value of the fidelity term computed using $\phi^*$, i.e., the forward state computed with the optimal control. The initial and target data are taken to be the same as in \S \ref{subsec:results-single} and we employ the algorithm with the volume constraint. The initial guess for the control is taken to be zero in each case.
\begin{table}[h!]
\centering
\begin{tabular}{c| c  c}
\hline
$DOFs$ &$|| \varphi(\vec x,t)-\varphi_{obs}||_{L_2(\Omega)}$&CPU time (secs)\\
\hline
545&0.348165&3142\\
2113&0.176357&26050\\
8321&0.154305&119216\\
33025&0.135178&291597\\
 \hline
\end{tabular}
\caption
   {Mesh refinement.}
\label{converg-test-table}
\end{table}
The results from the mesh refinements are presented in Table \ref{converg-test-table}. We observe a reduction in the fidelity term as we refine the mesh which implies an improved fit to the observed data. Although in principle it would be interesting to investigate the influence of refining both the timestep and mesh-size on the computed results, our tests indicate that the algorithm breaks down for time steps significantly larger than $0.01$ and hence, refinement of  the timestep and mesh-size together becomes computationally prohibitive.

\subsection{The influence of the initial guess for the control}
\label{subsec:initial-guess}

Here we apply the algorithm with the volume constraint on the simple example of a translated circle to illustrate the effect that the choice of the initial guess for the control $\eta$ has on the solution of the problem. 
To apply our algorithm we define the domain $\Omega$ to be $[-3,6]\times[-3,3]$ with a triangulation of $8321$ grid points. We selected a uniform timestep $\tau=1\times10^{-3}$ and set the interfacial thickness $\varepsilon=0.1$. We took the end-time $T=0.8$. The remaining numerical parameters for the optimisation algorithm are as given in Table \ref{table1}.  The initial data was taken to be a smoothed (by running a few steps of the Allen-Cahn solver) version of the function taking the value 1 inside $B_1(0,0)$ (a circle of radius 1  centred at the origin) and -1 in $\O/B_1(0,0)$.  The target data was taken to be a smoothed (by running a few steps of the Allen-Cahn solver) version of the function taking the value 1 inside $B_1(3,0)$  and -1 in $\O/B_1(0,0)$. Figure \ref{fig:circle3_dif} shows the initial and target diffuse interface data. 
\begin{figure}[h!]
\begin{center}
\subfigure[Initial data ($\vp^0$).]{\label{fig:circle3_initial}
{\includegraphics[trim = 40mm 95mm 30mm 95mm, clip,width=0.48\linewidth]{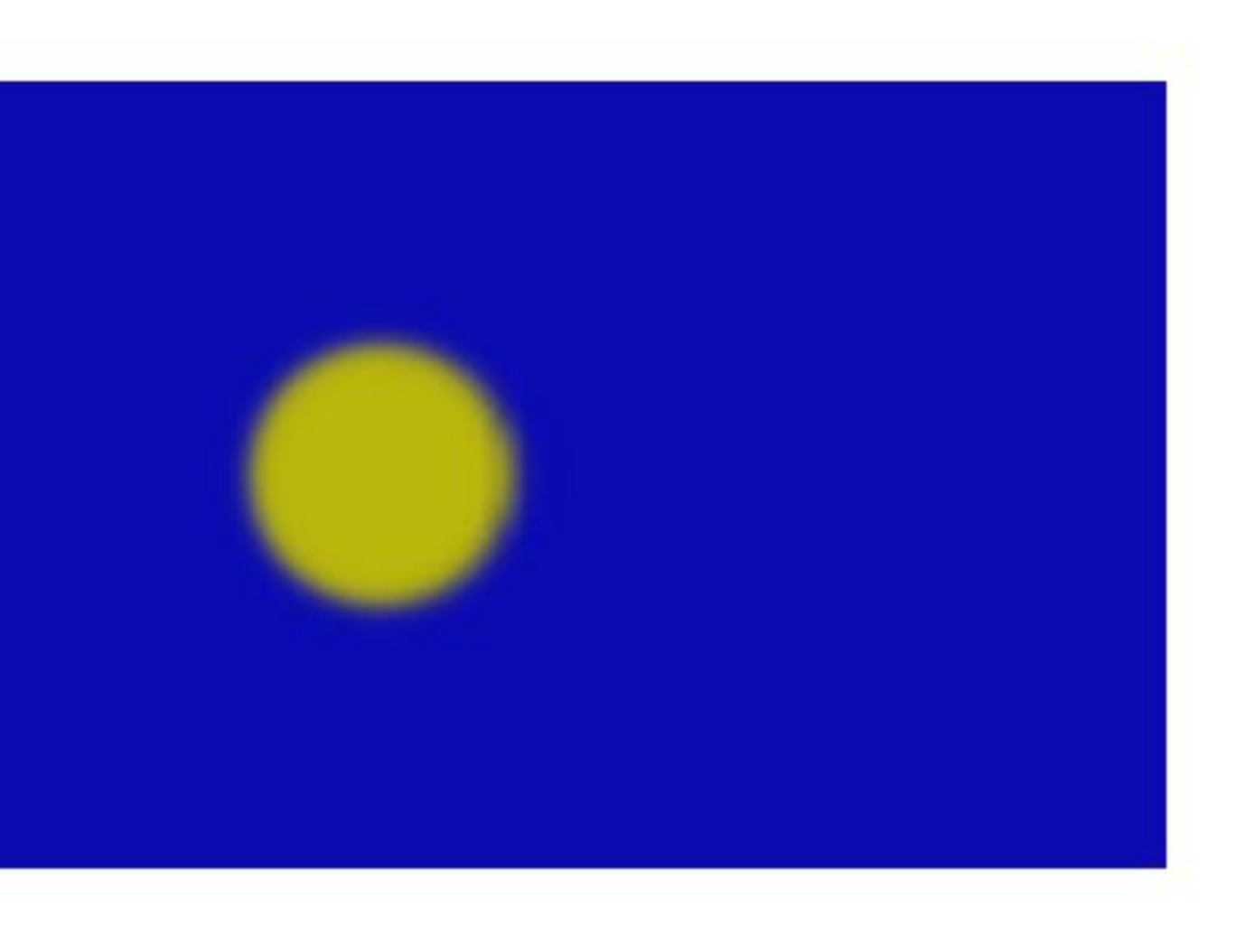}}}
\subfigure[Target data ($\vpo$).]{\label{fig:circle3_target}
{\includegraphics[trim = 40mm 95mm 30mm 95mm,clip, width=0.48\linewidth]{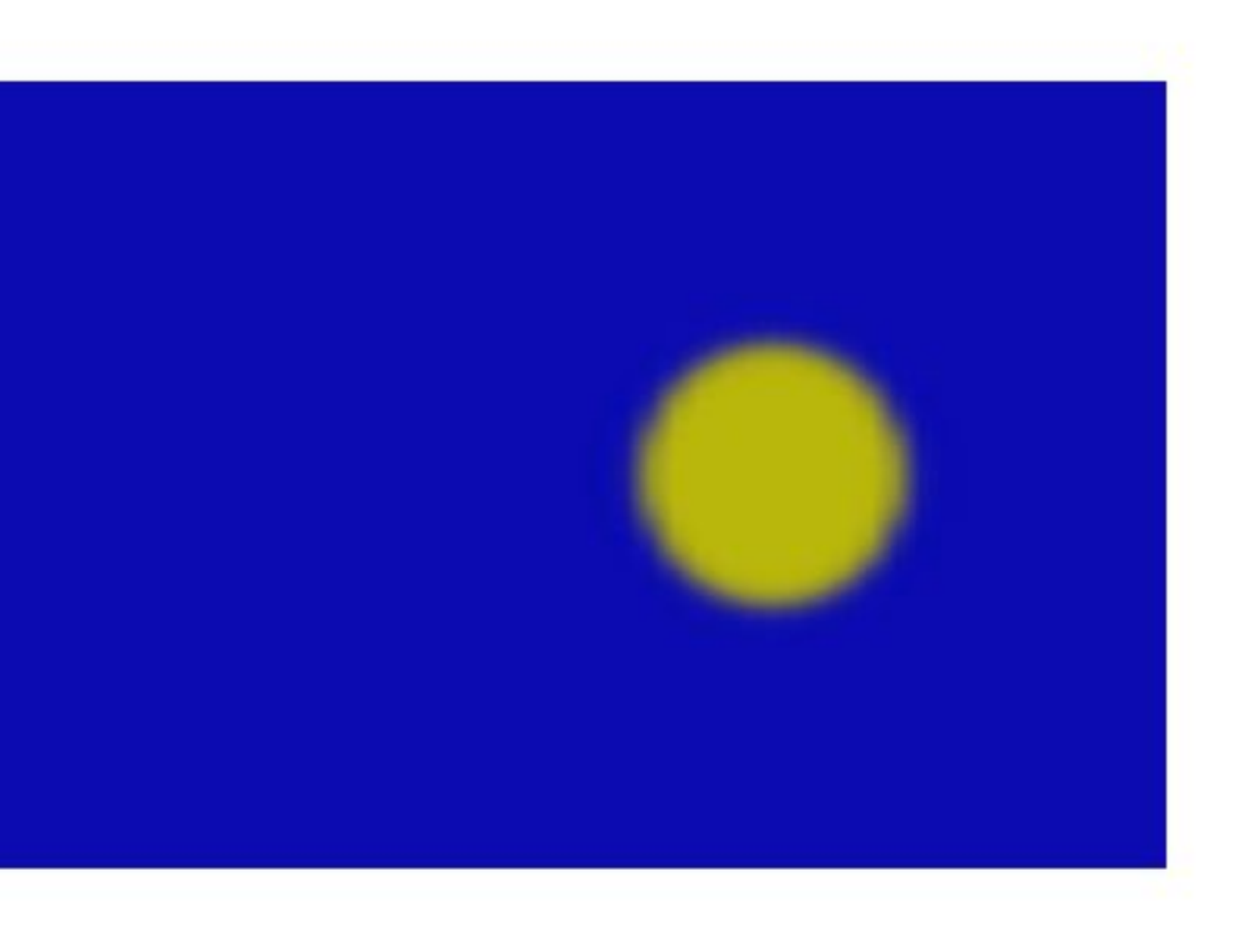}}}
\end{center}
 \caption
{Initial and target data for the examples of \S \ref{subsec:initial-guess}.}
 \label{fig:circle3_dif}
   \end{figure}
To illustrate the effect of the choice of initial guess on the algorithm, we consider two different values for the initial guess, firstly we set $\eta=0$ and secondly we set $\eta=\vec c\cdot\nabla \varphi$, where $\vec c=(2.5, 0)$, i.e., in the latter case the initial guess depends on the solution to the Allen-Cahn equation.  In both cases we used the algorithm with the volume constraints. With the zero initial guess the algorithm took 3262 iterations to meet the stopping criteria corresponding to a CPU time of 320433 seconds. With the second choice of initial guess the algorithm took 2056 iterations to meet the stopping criteria corresponding to a CPU time of 228173 seconds respectively.

Figure \ref{fig:circle3_optimal_target} shows the zero level-set of the computed solution using the optimal control at the final time. The curve corresponding to the zero level-set is shaded by the value of the control with the background shading corresponding to the target data. In both cases the position of the computed curve (zero level-set) with the optimal control shows good agreement with the target data. 
Figure \ref{fig:circle3_optimal_curves} shows snapshots of the computed zero level-sets with the two different initial guesses.  For the case with the initial value of $\eta=0$, we observe in Figure \ref{fig:circle3-sharp-eta0} that the interface remains close to the initial position for most of the time of the simulation, and at the very last moment it shrinks to a point with a new phase nucleated at the position of the target data corresponding to a change in topology. With the second choice of initial guess $(\eta=\vec c \cdot \nabla \varphi)$ we observe in Figure \ref{fig:circle3-sharp-etag} that there is a gradual motion towards the target position with no changes in topology.  Figure \ref{fig:area_circle3} shows the area enclosed by the zero level-set of the computed solution with the optimal control with the two different initial guesses together with the linear interpolant of the areas of the data. We observe a sharp increase in area towards the end of the time interval with the zero initial guess as the new phase is nucleated. With the second choice of initial guess, the area of the computed curve exhibits a good fit to the linear interpolant of the areas of the data. 
  \begin{figure}[h!]
\begin{center}
\subfigure[With initial guess $\eta=0$]
{{\includegraphics[trim = 40mm 95mm 30mm 80mm, clip,width=0.48\linewidth]{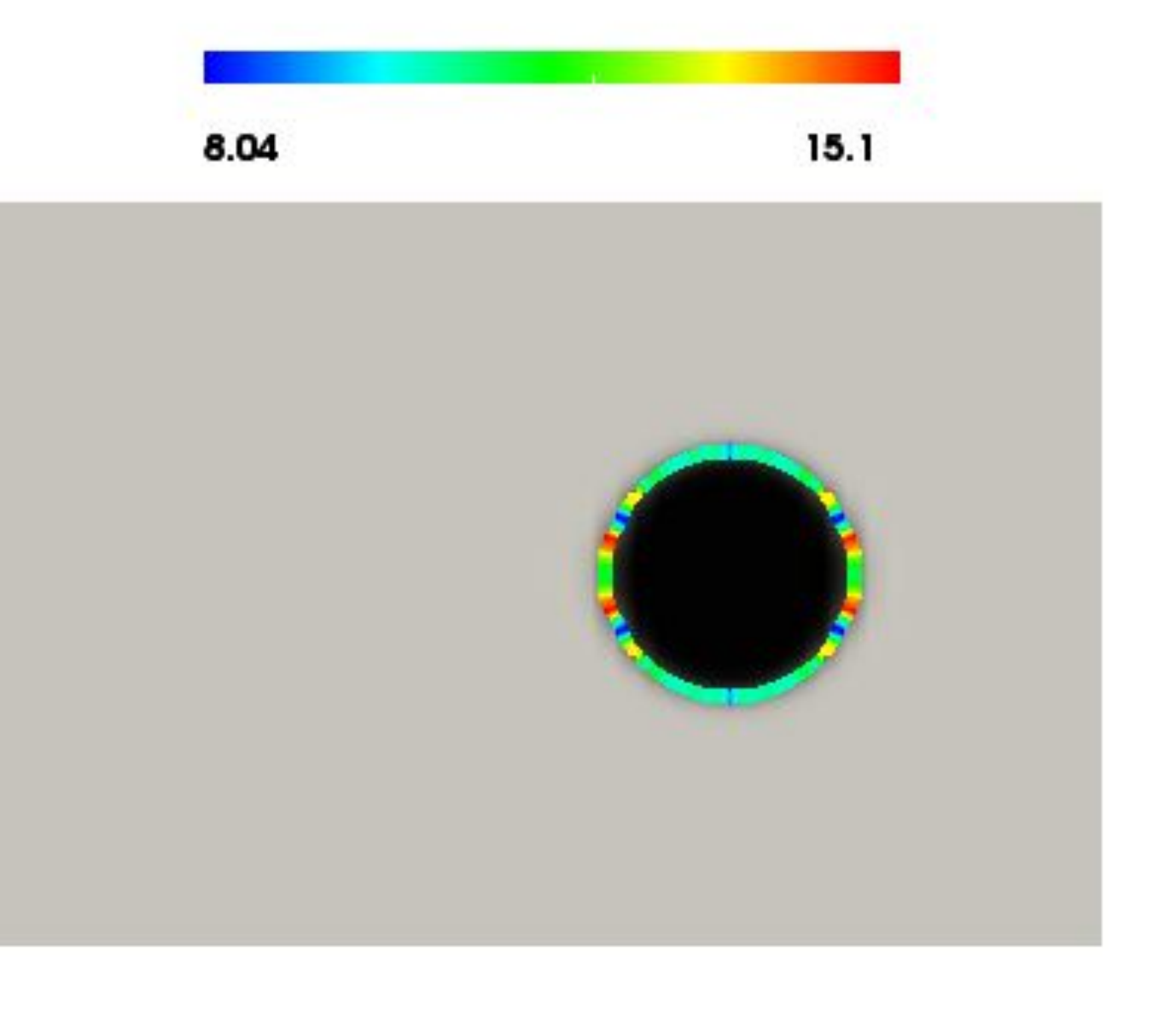}}}
\subfigure[With initial guess $\eta=\vec c\cdot \nabla \varphi$]
{{\includegraphics[trim = 40mm 95mm 30mm 80mm, clip,width=0.48\linewidth]{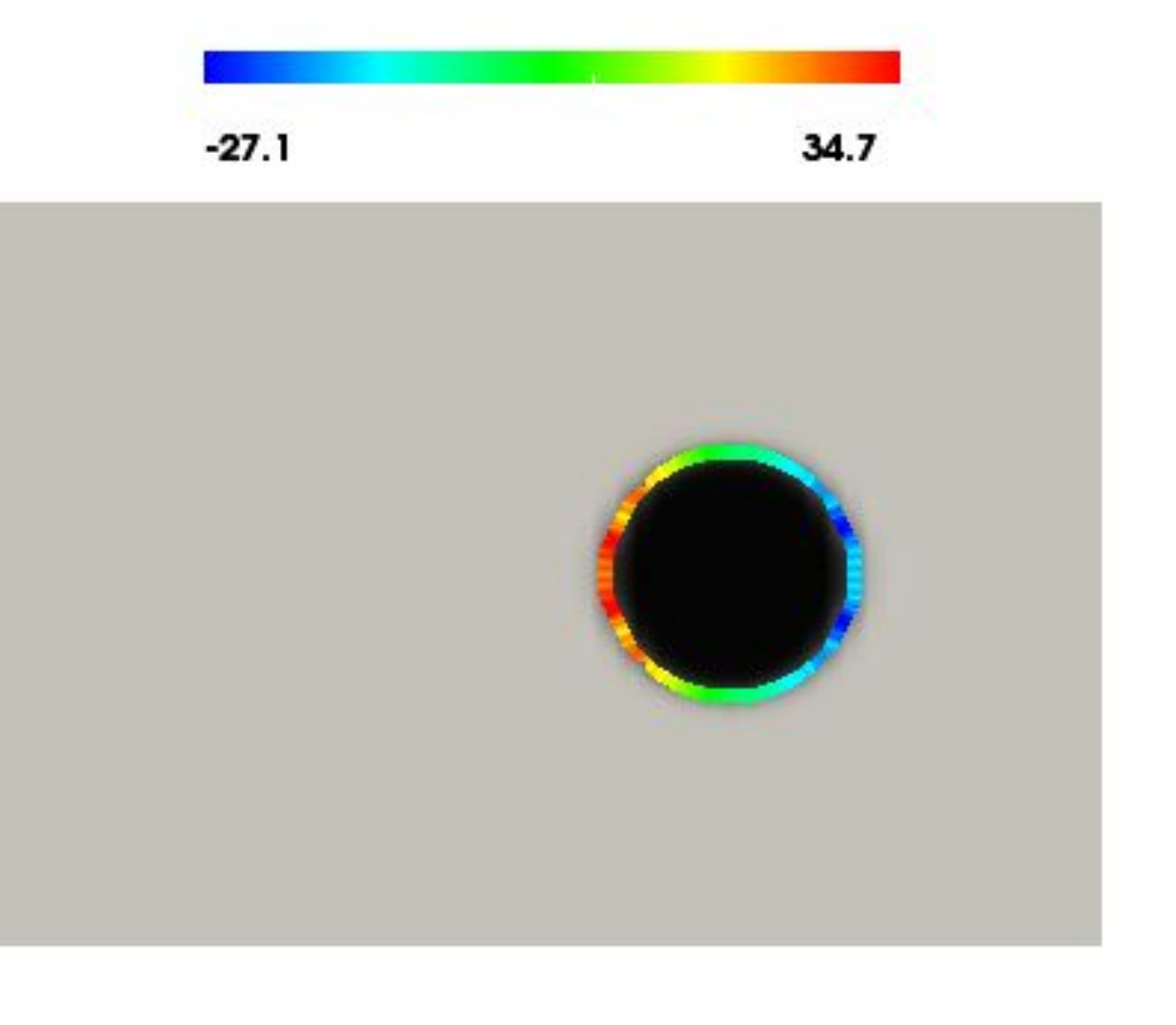}}}
\end{center}
 \caption{
 Zero level-set of the solutions ($\vp(\vec x,T)$) computed using the approximated optimal control ($\eta^*(\vec x,t)$) for the experiments of \S\ref{subsec:initial-guess}.
 The curve (zero level-set of $\vp(\vec x,T)$)  is shaded by the approximated optimal control ($\eta^*(\vec x,T)$) and the background by the target data ($\vpo(\vec x)$). The color-bar corresponds to the scale for $\eta^*(\vec x,T)$.  We see good agreement between the zero level-set of the data computed with the optimal control and the target data in both cases.
 }\label{fig:circle3_optimal_target}
  \end{figure}
       \begin{figure}[h!]
\begin{center}
\subfigure[With initial guess $\eta=0$]
{\label{fig:circle3-sharp-eta0}
{{\includegraphics[trim = 20mm 95mm 20mm 95mm, clip,width=0.48\linewidth]{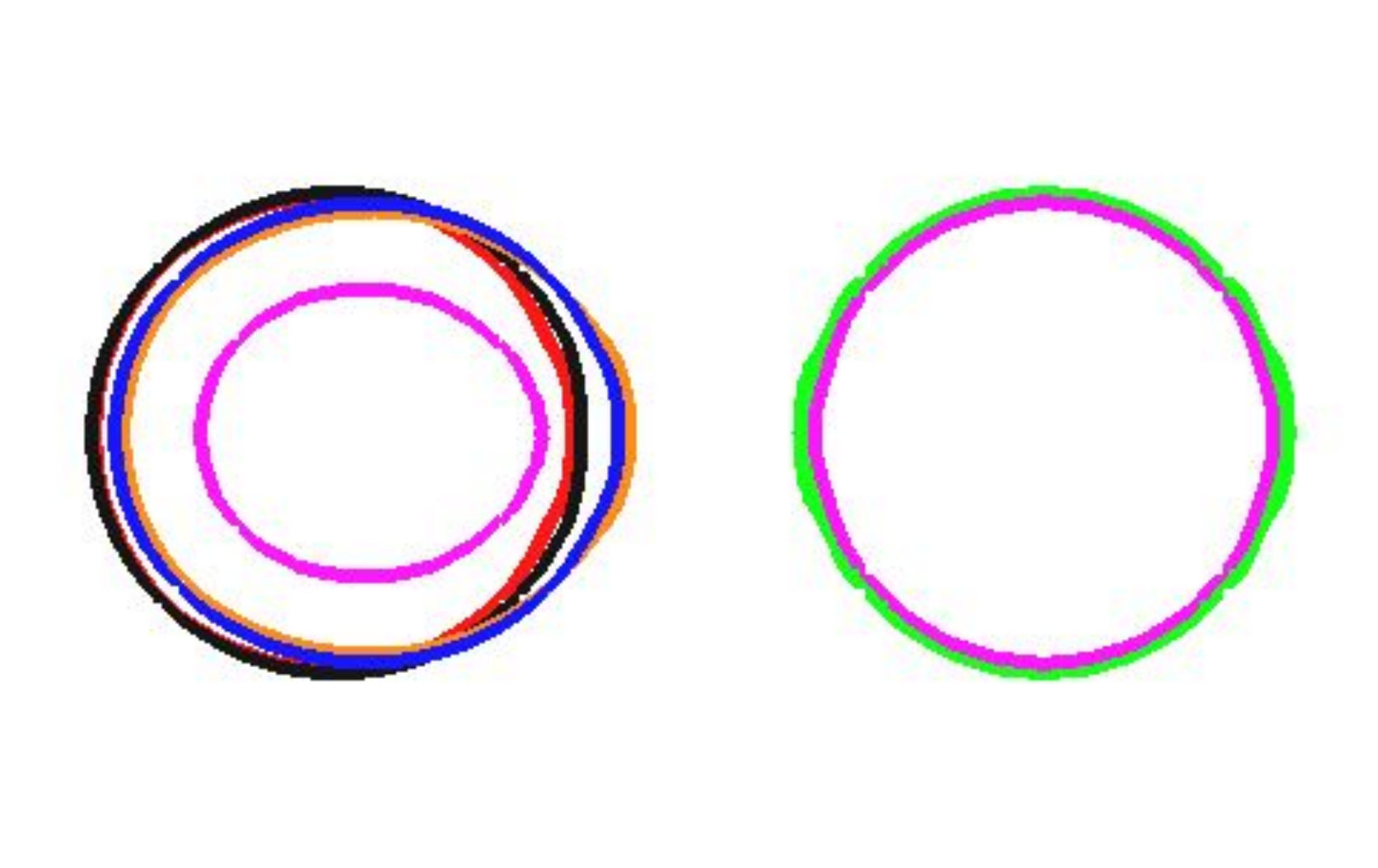}}}}
\subfigure[With initial guess $\eta=\vec c\cdot \nabla \varphi$]
{\label{fig:circle3-sharp-etag}
{{\includegraphics[trim = 20mm 95mm 30mm 95mm, clip,width=0.48\linewidth]{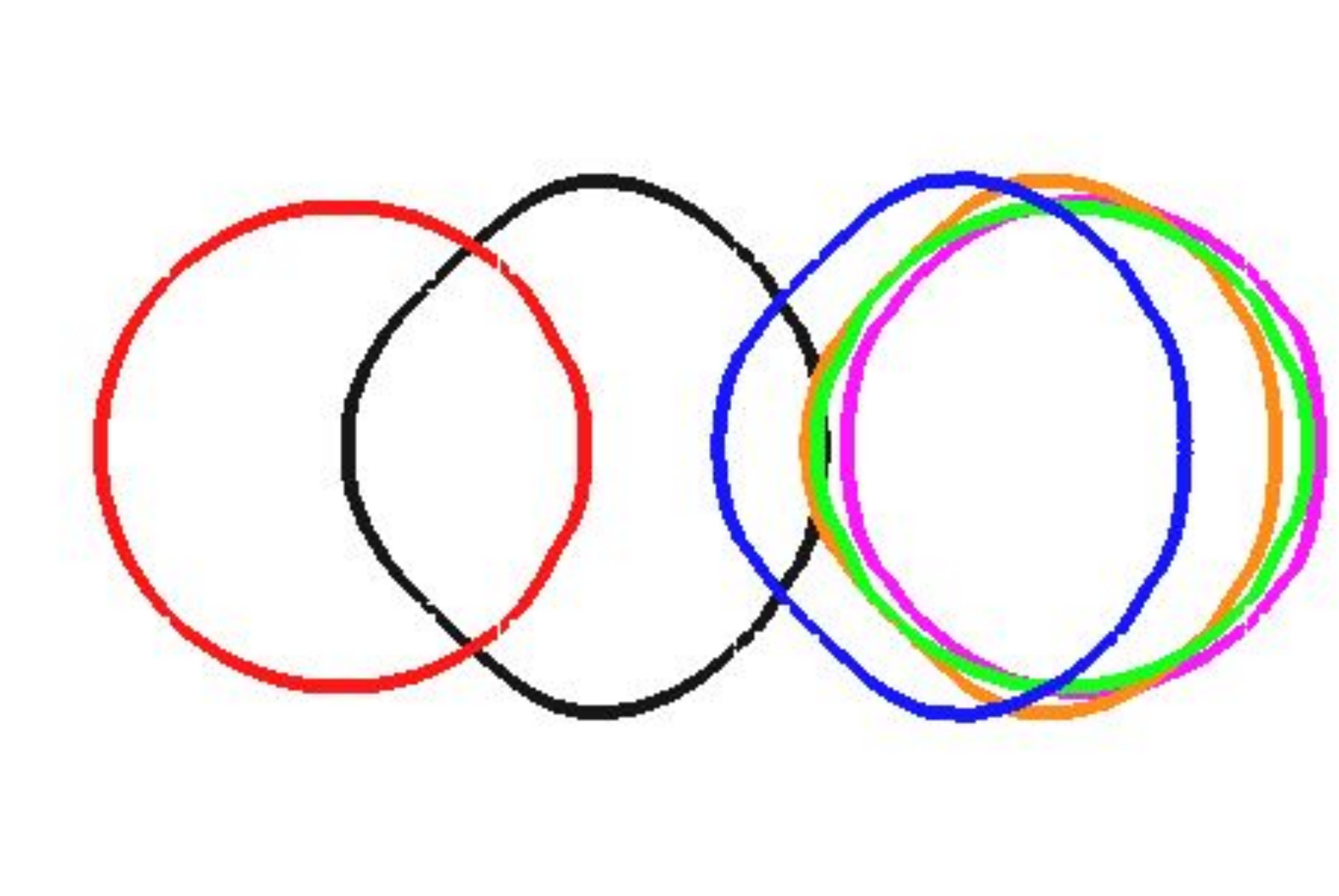}}}}
\end{center}
  \caption
{Zero level-sets of the solutions computed ($\vp(\vec x,t)$) with the optimal control ($\eta^*(\vec x,t)$) for the experiments of \S \ref{subsec:initial-guess} at $t=0$ (red), $t=0.2$ (black), $t=0.6$ (blue), $t=0.7$ (orange), $t=0.789$ (pink) and $t=0.8$ (green). We observe the nucleation of a phase and a change in topology with the zero initial guess whilst there are no evident changes in topology and the zero level-set maintains a fixed  topology in the case of the nonzero initial guess.}
 \label{fig:circle3_optimal_curves}
  \end{figure}
  \begin{figure}[h!]
\begin{center}
{{\includegraphics[trim = 20mm 0mm 30mm 0mm,clip,width=\linewidth]{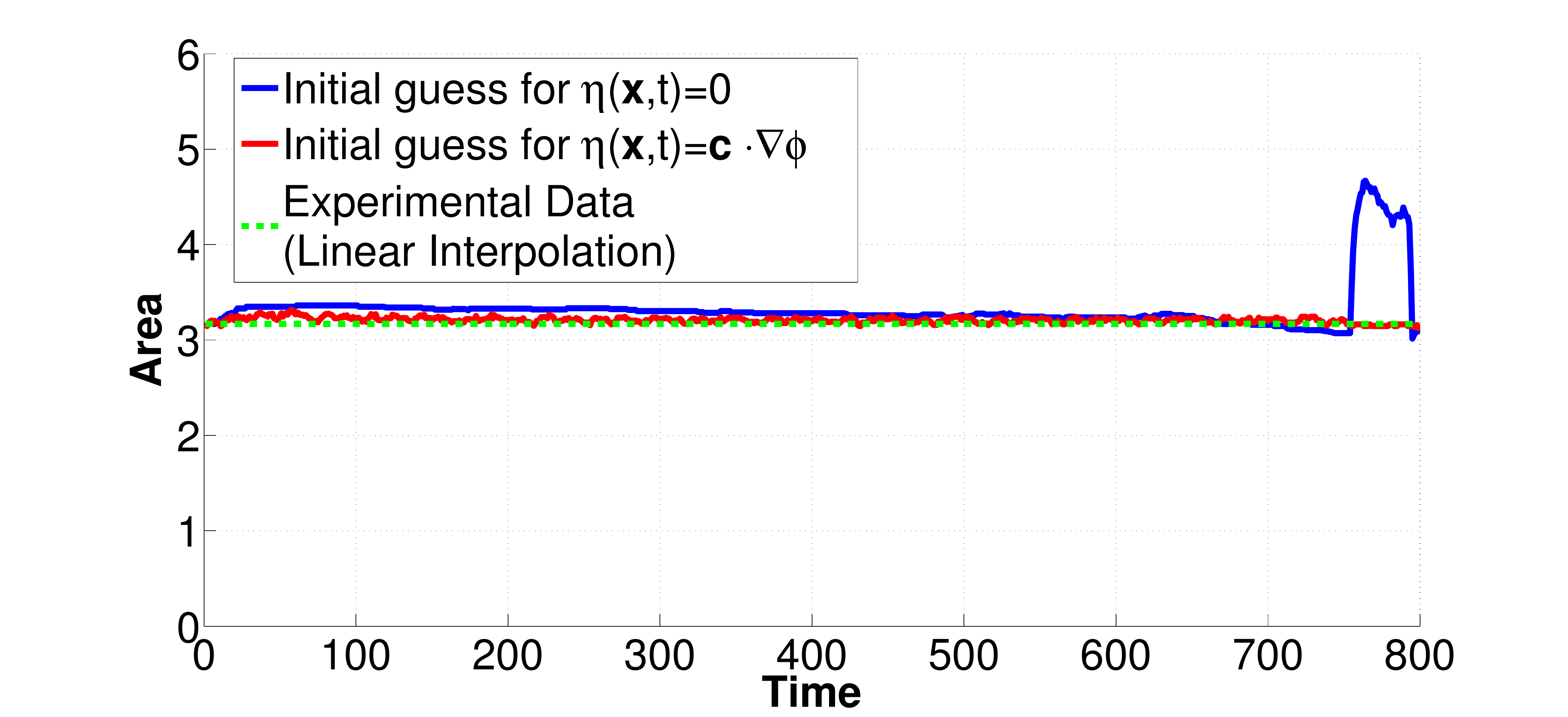}}}
\end{center}
 \caption{
Area enclosed by the curve for the experiments of \S \ref{subsec:initial-guess}. A good fit to the linear interpolant of the areas is only observed with the nonzero initial guess. We observe a rapid increase in the area near the end time for the zero initial guess, this corresponds to the time at which a new phase is nucleated, c.f., Figure \ref{fig:circle3-sharp-eta0}.}\label{fig:area_circle3}
  \end{figure}

\subsection{Application to multi-cell image data sets}
\label{subsec:results-multi}
We now apply the algorithm to the case of multi-cell image data sets. As a proof-of-concept we consider the simplest possible scenario where we have an initial and desired dataset both consisting of two cells that are well separated. 

For the first experiment we defined the initial data and target data as follows. Defining the domain $\O$ to be $[-2,8]\times[-2,2]$ we defined the subdomains $\O_1,$  $\O_2,$ $\O_3$ and $\O_4$ to be the simply connected bounded domains with boundary curves $\G_1,\G_2,\G_3$ and $\G_4$ defined by (the curves $\G_1,\G_2$ and $\G_3$ and $\G_4$ are the zero level-sets of the diffuse interfaces shown in Figures \ref{fig:multi_1_data_PF_intial} and \ref{fig:multi_1_data_PF_target} respectively).
\begin{align*}
&\G_1:=\left\{x\in\O\left|\right.x_{1}^2+x_{2}^2-0.8^2+0.1\sin(4x_{1})+0.1\sin(3x_{2})=0\right\},\\
&\G_2:=\left\{x\in\O\left|\right.\left(\frac{x_{1}}{2}-2\right)^2+\left(x_{2}-0.6\right)^2-0.7^2+0.1\sin\left(\frac{5x_{1}}{2}\right)+0.3\sin(2x_{2})=0\right\},\\
&\G_3:=\left\{x\in\O\left|\right.(x_{1}-0.4)^2+(x_{2}-0.5)^2-0.8^2+0.1\sin(6x_{1})+0.1\sin(7x_{2})=0\right\},\\
&\G_4:=\left\{x\in\O\left|\right.\left(\frac{x_{1}}{2}-2.5\right)^2+(x_{2}-1)^2-0.7^2+0.1\sin\left(\frac{7x_{1}}{2}\right)+0.1\sin(1.5x_{2})=0\right\}.
\end{align*}
We then set the initial and target data to be a smoothed (by running a few steps of the Allen-Cahn solver) version of the function 
\beq
\notag
\vp^0=
\begin{cases}
1\text{ for }\vec x\in\O_1\cup\O_2,\\
-1\text{ for }\vec x\in\O/\left(\O_1\cup\O_2\right),
\end{cases}
\text{ and }
\quad
\vpo
=
\begin{cases}
1\text{ for }\vec x\in\O_{3}\cup\O_{4},\\
-1\text{ for }\vec x\in\O/\left(\O_{3}\cup\O_{4}\right).
\end{cases}
\eeq
Figure \ref{fig:multi_initial_target} shows the initial and target diffuse interface data.

 As previously, we compare the results of the algorithm with and without the volume constraint.
For this experiment, the algorithm took 2035 iterations to meet the stopping criteria with no volume constraint  and 2199 iterations with the volume constraint, corresponding to CPU times of  28608 and 105750 seconds respectively.

 Figure \ref{fig:cost_multi} shows the value of the objective functional against the number of iterations of the optimisation algorithm with and without the volume constraint.  Figure \ref{fig:multi_optimal_target} shows the zero level-set of the computed solution using the optimal control at the final time with and without the volume constraint shaded by the value of the control with the background shading corresponding to the target data. The results are similar to the single cell simulations of \S \ref{subsec:results-single} with an initial rapid decrease in the cost followed by a subsequent gradual decrease. The cells (zero level-sets) computed with the optimal control show good agreement with the target data for both versions of the algorithm and for both cells.   For each of the versions of the algorithm, both of the computed cells again posses a clearly  defined ``front'' and ``rear'' similar to the single cell case. 
 
Figure \ref{fig:area_multi} shows the area enclosed by the zero-level set of the computed solution with the optimal control with and without the volume constraint together with the linear interpolant of the areas of the data. We observe analogous behaviour to the single cell. In terms of the computed  cell morphologies, Figure \ref{fig:multi_optimal_curves} shows snapshots of the computed  zero level-sets for the two different versions of the algorithm. We see that in this multi-cell setting the algorithm has implicitly solved the matching problem by generating two disjoint cells whose topology remains fixed throughout the evolution.  We  observe that the loss of volume in the case of no volume constraint corresponds to one of the cells in the intermediate snapshot (blue curve) enclosing a much smaller area.
\begin{figure}[h!]
\begin{center}
\subfigure[Initial data ($\vp^0$).]{\label{fig:multi_1_data_PF_intial}
{\includegraphics[trim = 10mm 195mm 20mm 20mm, clip,width=0.48\linewidth]{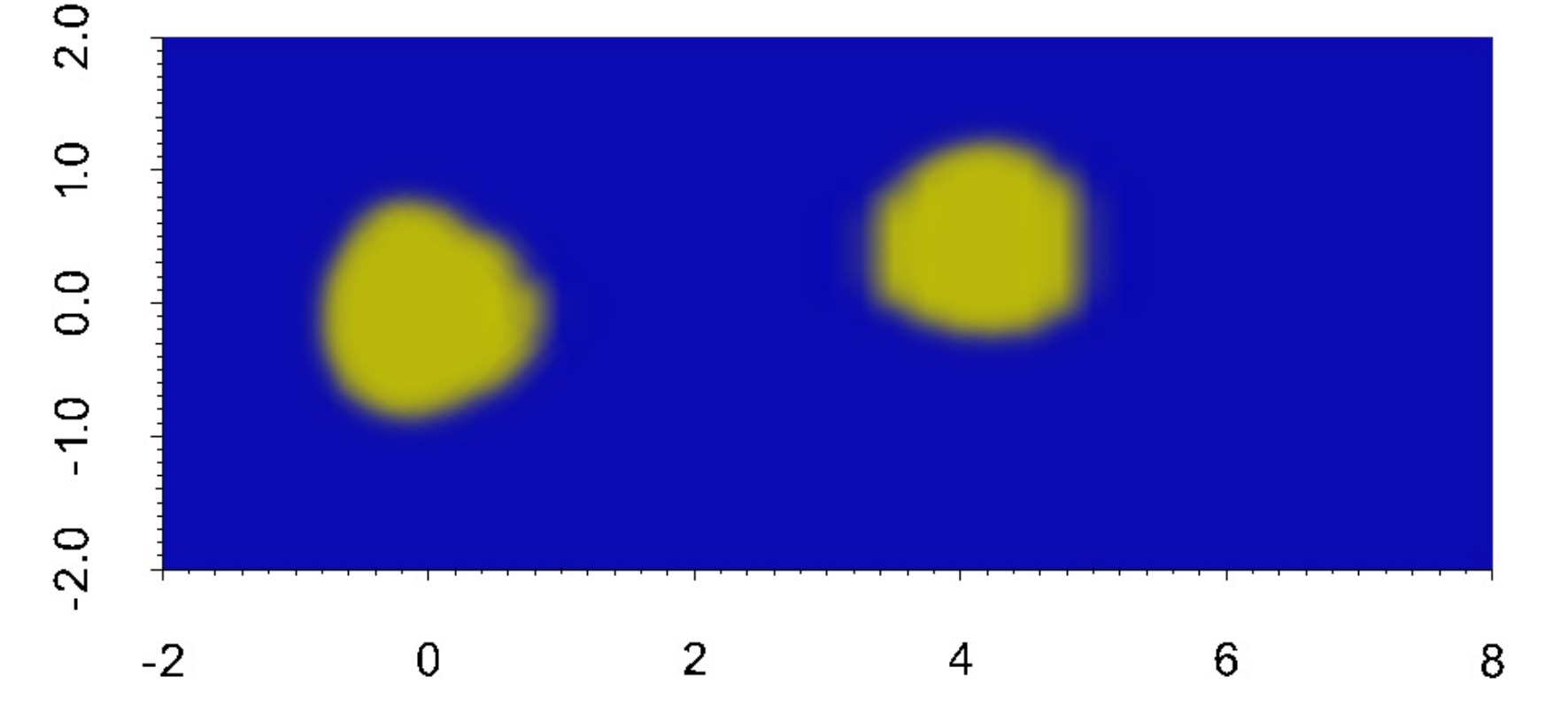}}}
\subfigure[Target data ($\vpo$).]{\label{fig:multi_1_data_PF_target}
{\includegraphics[trim = 10mm 195mm 20mm 20mm,clip, width=0.48\linewidth]{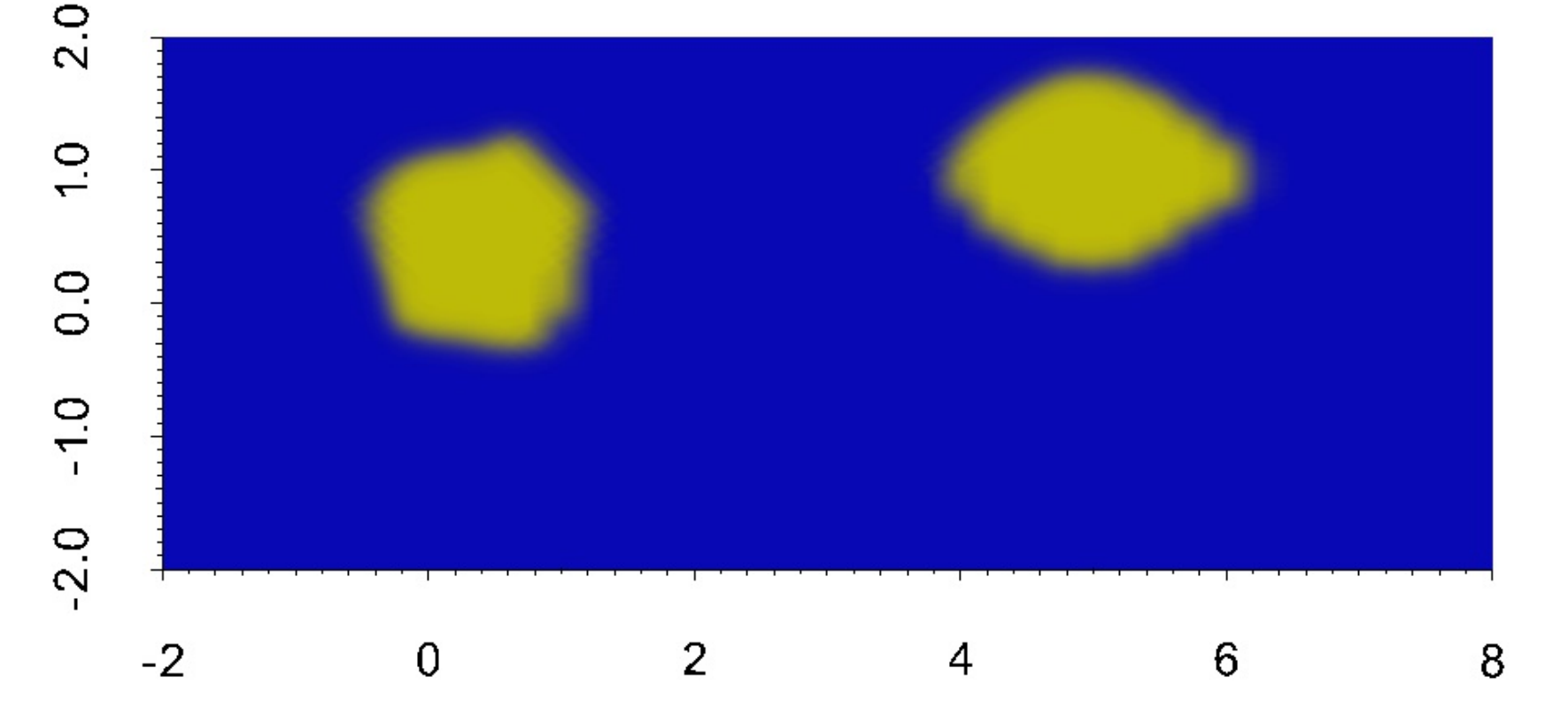}}}
\end{center}
 \caption{Initial and target data for the examples of \S \ref{subsec:results-multi}.}
 \label{fig:multi_initial_target}
   \end{figure}
     \begin{figure}[h!]
\begin{center}
\subfigure[Without the volume constraint]
{{\includegraphics[trim = 10mm 0mm 20mm 0mm,clip,width=0.49\linewidth]{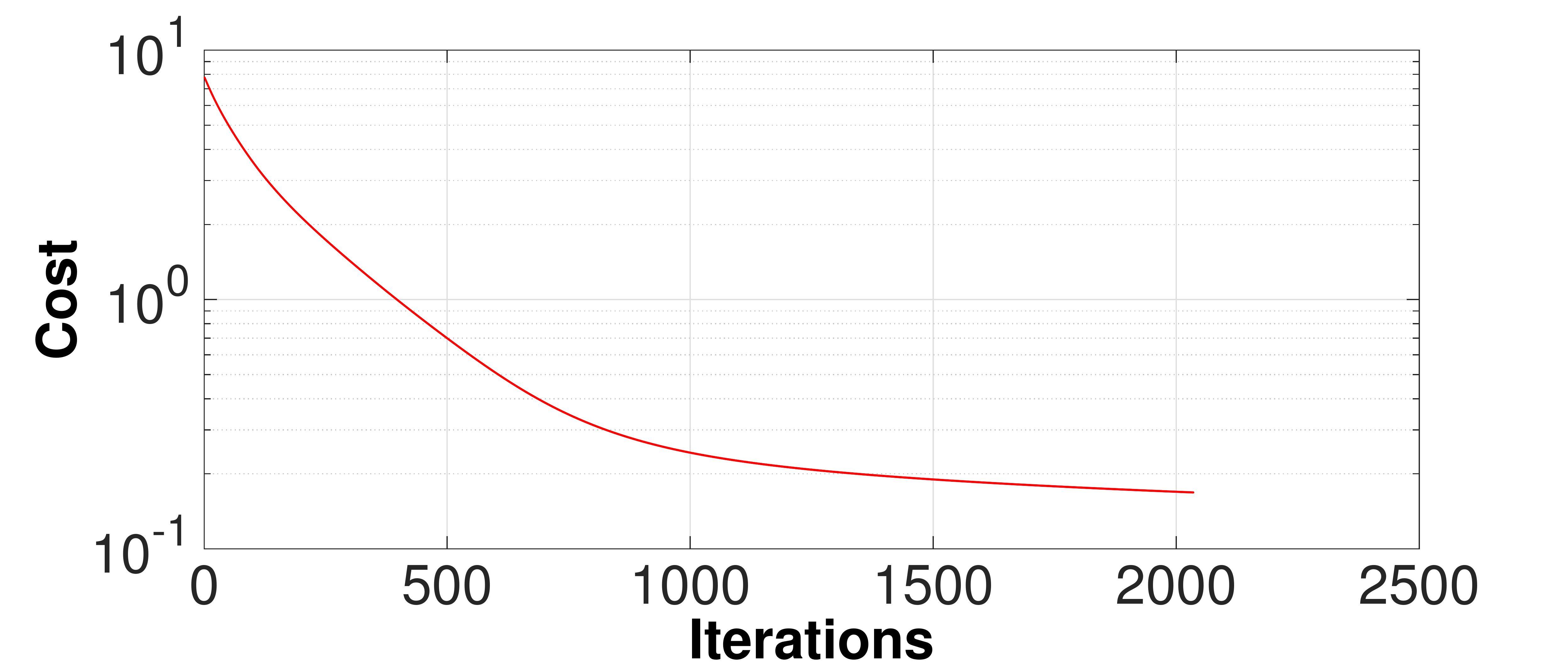}}}
\subfigure[With the volume constraint]
{{\includegraphics[trim = 10mm 0mm 20mm 0mm,clip,width=0.49\linewidth]{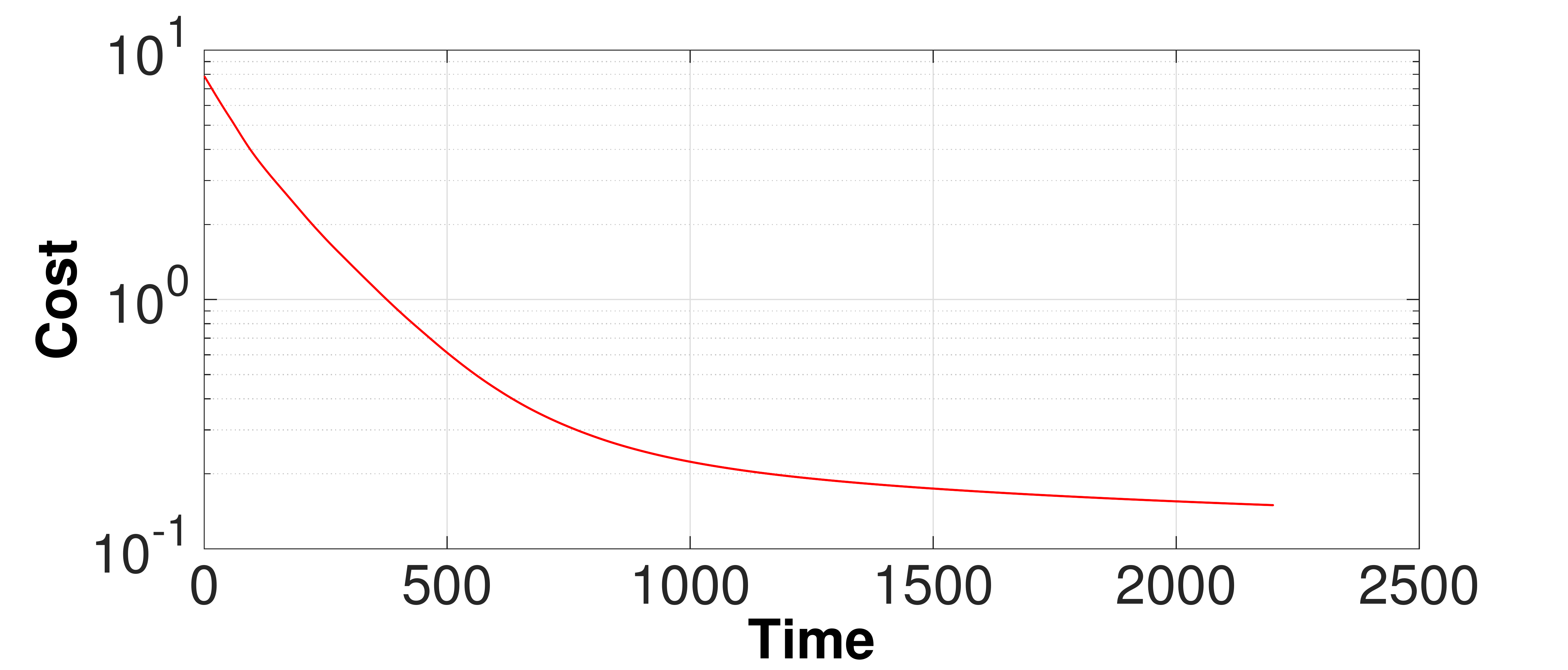}}}
\end{center}
 \caption{Cost functional versus the number of iterations for the examples of \S \ref{subsec:results-multi}. }\label{fig:cost_multi}
  \end{figure}
  \begin{figure}[h!]
\begin{center}
\subfigure[Without the volume constraint]
{{\includegraphics[trim = 10mm 190mm 20mm 10mm,clip,  width=0.49\linewidth]{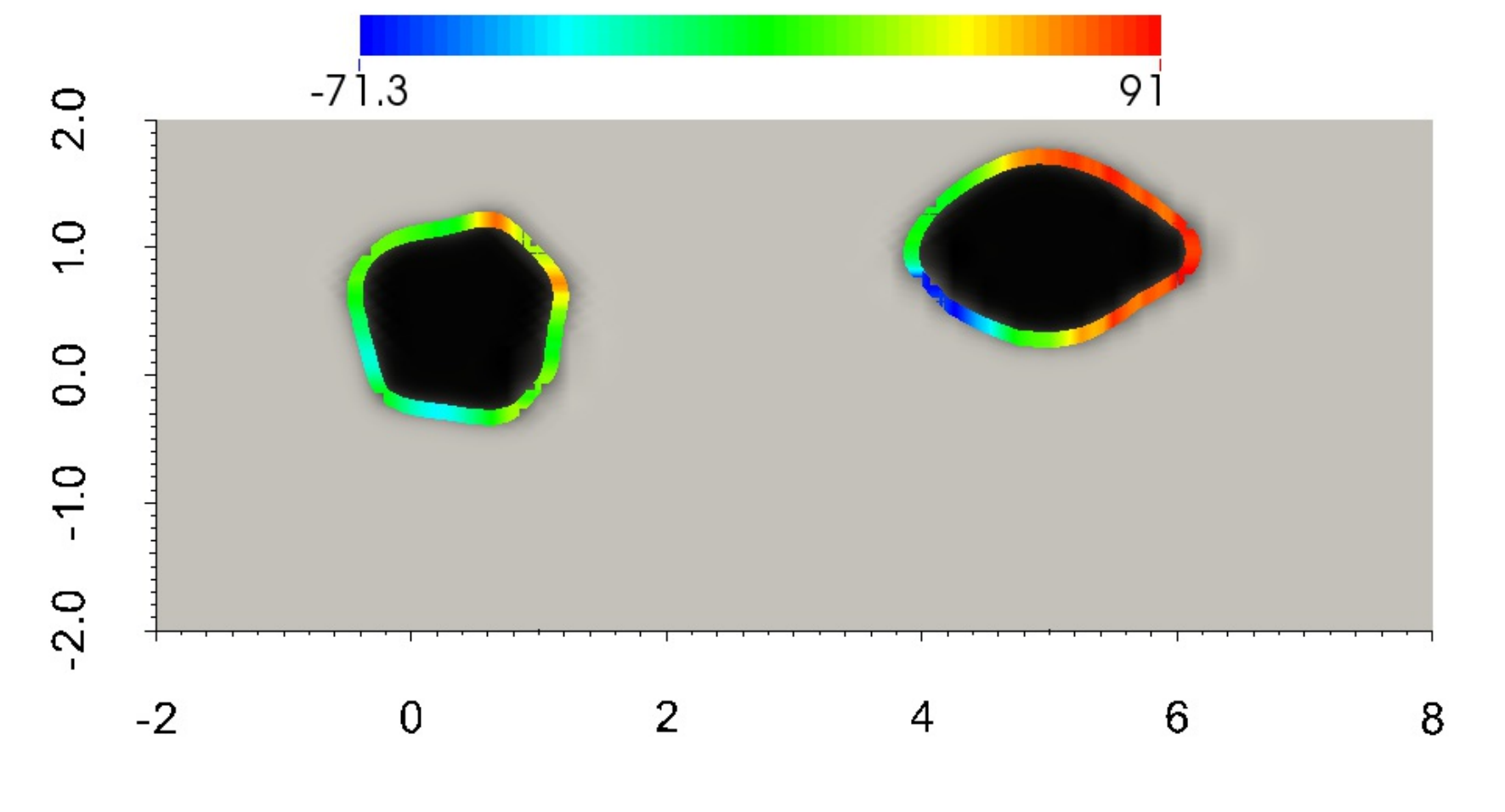}}}
\subfigure[With the volume constraint]
{{\includegraphics[trim = 10mm 190mm 20mm 10mm,clip,  width=0.49\linewidth]{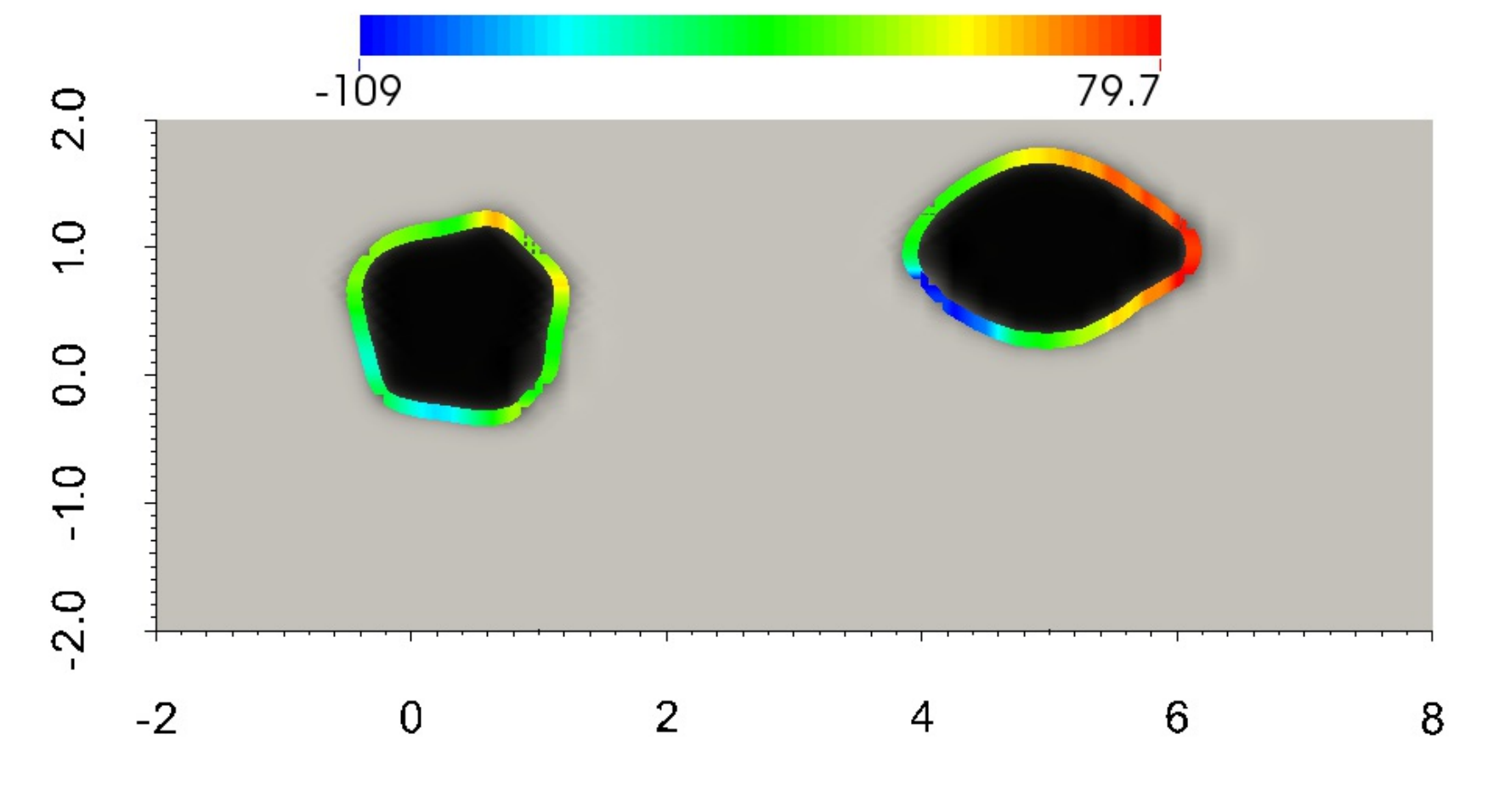}}}
\end{center}
 \caption{
 Zero level-set of the solutions ($\vp(\vec x,T)$) computed using the approximated optimal control ($\eta^*(\vec x,t)$) with and without the volume constraint for the examples of \S \ref{subsec:results-multi}.
 The curve (zero level-set of $\vp(\vec x,T)$)  is shaded by the approximated optimal control ($\eta^*(\vec x,T)$) and the background by the target data ($\vpo(\vec x)$). The color-bar corresponds to the scale for $\eta^*(\vec x,T)$.  We see good agreement between the zero level-set of the data computed with the optimal control and the target data in both cases.
 }\label{fig:multi_optimal_target}
  \end{figure}
    \begin{figure}[h!]
\begin{center}
{{\includegraphics[trim = 40mm 0mm 30mm 0mm,clip,width=\linewidth]{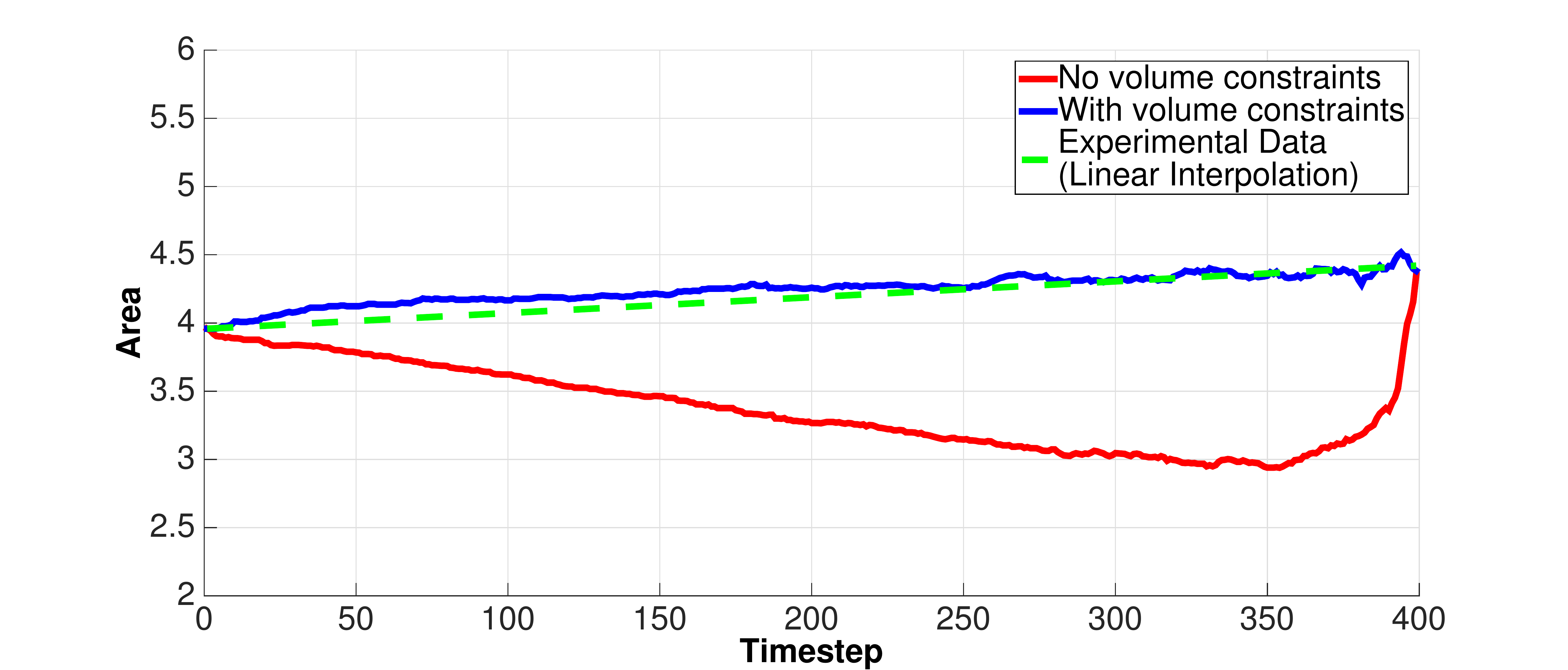}}}
\end{center}
 \caption{Area enclosed by the cell for the experiments of \S \ref{subsec:results-multi} with and without the volume constraint. As with the single cell data, the area (now the sum of the area of the two cells) shrinks considerably during the evolution without the volume constraint whilst a good fit to the linear interpolant of the area enclosed by the data is observed with the volume constraint.}\label{fig:area_multi}
  \end{figure}
\begin{figure}[h!]
\begin{center}
\subfigure[Without the volume constraint]
{\label{fig:single_ini_target_1_opt}{\includegraphics[trim = 40mm 220mm 40mm 10mm,clip,  width=0.45\linewidth]{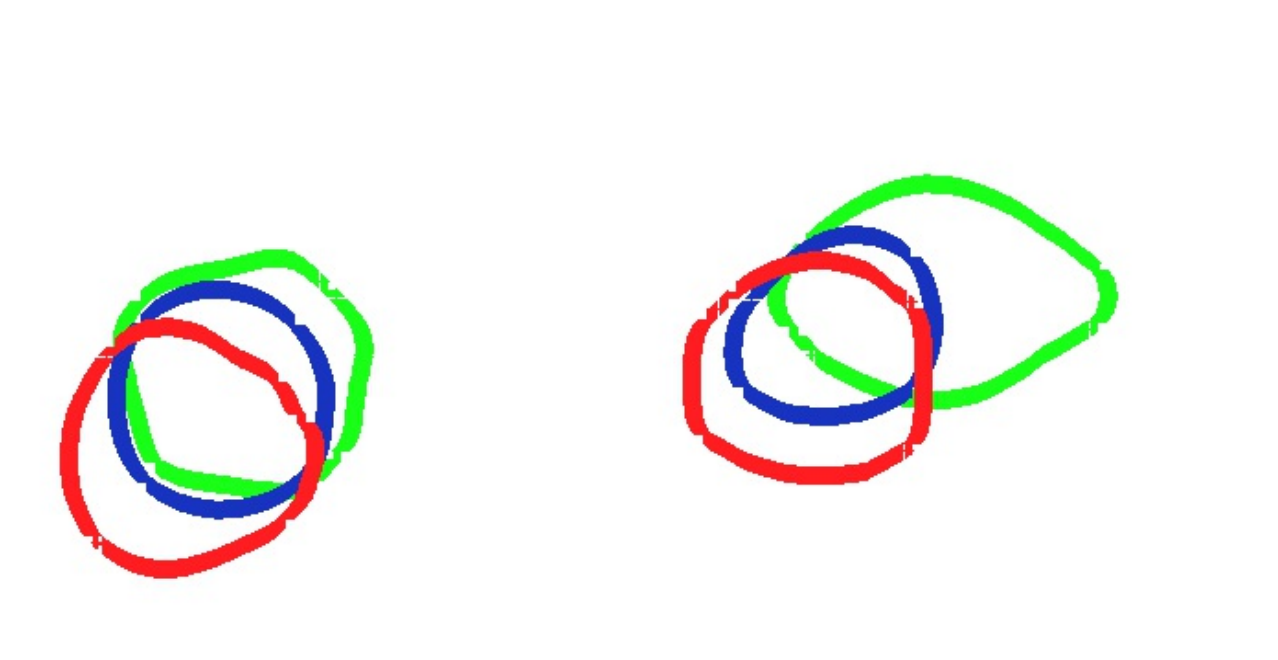}}}
\hskip1em
\subfigure[With the volume constraint]
{\label{fig:single_ini_target_1_opt}{\includegraphics[trim = 40mm 220mm 40mm 10mm,clip,  width=0.45\linewidth]{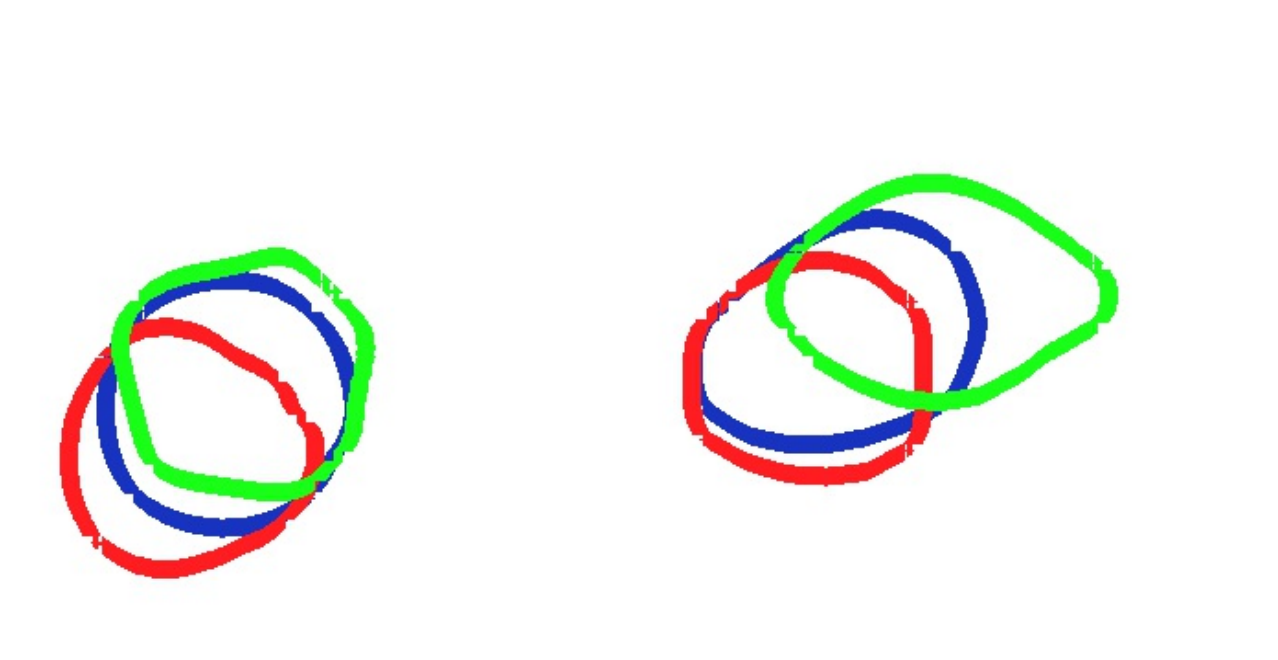}}}
\end{center}
 \caption{Zero level-sets of the solutions computed with the optimal control for the examples of \S \ref{subsec:results-multi} with and without the volume constraint at $t=0$ (red), $t=0.35$ (blue) and $t=0.4$ (green). The volume enclosed by both cells shrinks during the evolution without the volume constraint whilst this is not observed if the volume constraint is included.  Both with and without the volume constraint,  the implicit solution of the matching problem in this case generates two disjoint cells which do not change in topology.}\label{fig:multi_optimal_curves}
  \end{figure}

\subsection{An example with topological change}
\label{subsec:results-multi-splits}

Of course, in general our algorithm may generate cells whose topology is not fixed as in \S \ref{subsec:initial-guess}. To this end we report on another experiment.

Defining the domain $\O$ to be  $[-2,6.3]\times[-2.5,2.5]$ we defined the subdomains $\O_1,$  $\O_2,$ $\O_3$ and $\O_4$ to be the simply connected bounded domains with boundary curves $\G_1,\G_2,\G_3$ and $\G_4$ defined by (see Figure \ref{fig:splits_initial_target})
\begin{align*}
&\G_1:=\left\{x\in\O\left|\right.x_{1}^2+x_{2}^2-0.9^2+0.1\sin(4.5x_{1})+0.11\sin(3x_{2}))=0\right\},\\
&\G_2:=\left\{x\in\O\left|\right.\left(x_{1}-5\right)^2+x_{2}^2-0.7^2+0.1\sin\left(\frac{5x_{1}}{2}\right)+0.3\sin(2x_{2})=0\right\},\\
&\G_3:=\left\{x\in\O\left|\right.(x_{1}-0.35)^2+(x_{2}-0.7)^2-0.8^2+0.1\sin(6x_{1})+0.1\sin(7x_{2})=0\right\},\\
&\G_4:=\left\{x\in\O\left|\right.\left(x_{1}-0.3\right)^2+(x_{2}-1.1)^2-0.7^2-0.1\sin\left(\frac{7x_{1}}{2}\right)+0.1\sin(1.5x_{2})=0\right\}.
\end{align*}
We then set the initial and target data to be a smoothed (by running a few steps of the Allen-Cahn solver) version of the function 
\beq
\notag
\vp^0=
\begin{cases}
1\text{ for }\vec x\in\O_1\cup\O_2,\\
-1\text{ for }\vec x\in\O/\left(\O_1\cup\O_2\right),
\end{cases}
\text{ and }
\quad
\vpo
=
\begin{cases}
1\text{ for }\vec x\in\O_{3}\cup\O_{4},\\
-1\text{ for }\vec x\in\O/\left(\O_{3}\cup\O_{4}\right).
\end{cases}
\eeq
Figure \ref{fig:splits_initial_target} shows the initial and target diffuse interface data.

 As previously, we compare the results of the algorithm with and without the volume constraint.
For this experiment, the algorithm took 1960 iterations to meet the stopping criteria with no volume constraint  and 1937 iterations with the volume constraint, corresponding to CPU times of  27553 and 93150 seconds respectively.

 Figure \ref{fig:cost_splits} shows the value of the objective functional against the number of iterations of the optimisation algorithm with and without the volume constraint. Figure \ref{fig:splits_optimal_target} shows the zero level-set of the computed solution using the optimal control at the final time with and without the volume constraint shaded by the value of the control with the background shading corresponding to the target data. The results are similar to the previous simulations with an initial rapid decrease in the cost followed by a subsequent gradual decrease and good agreement with the target data for both versions of the algorithm and for both cells.   For each of the versions of the algorithm, both of the computed cells again posses a clearly  defined ``front'' and ``rear''. 
 
 Figure \ref{fig:area_splits} shows the area of the domain in which the computed solution is positive with and without the volume constraint together with the linear interpolant of the areas of the data. For this experiment we observe that the area enclosed by the computed cells differs significantly from the linear interpolant areas of the data both with and without the volume constraint. This may be due to the change in topology of the interface during the evolution.  Figure \ref{fig:splits_optimal_curves} shows snapshots of the computed  zero level-sets for the two different versions of the algorithm. Unlike the previous examples we see that for this particular choice of initial and target data, the algorithm yields cells which change in topology with one of the curves shrinking until it disappears whilst the other splits eventually becoming two disjoint curves. Thus our algorithm generates trajectories corresponding to the annihilation (via shrinking) of one cell whilst the other cell splits to form the two cells observed in the image data set. 
\begin{figure}[h!]
\begin{center}
\subfigure[Initial data ($\vp^0$).]{
{\includegraphics[trim = 10mm 140mm 20mm 20mm, clip,width=0.48\linewidth]{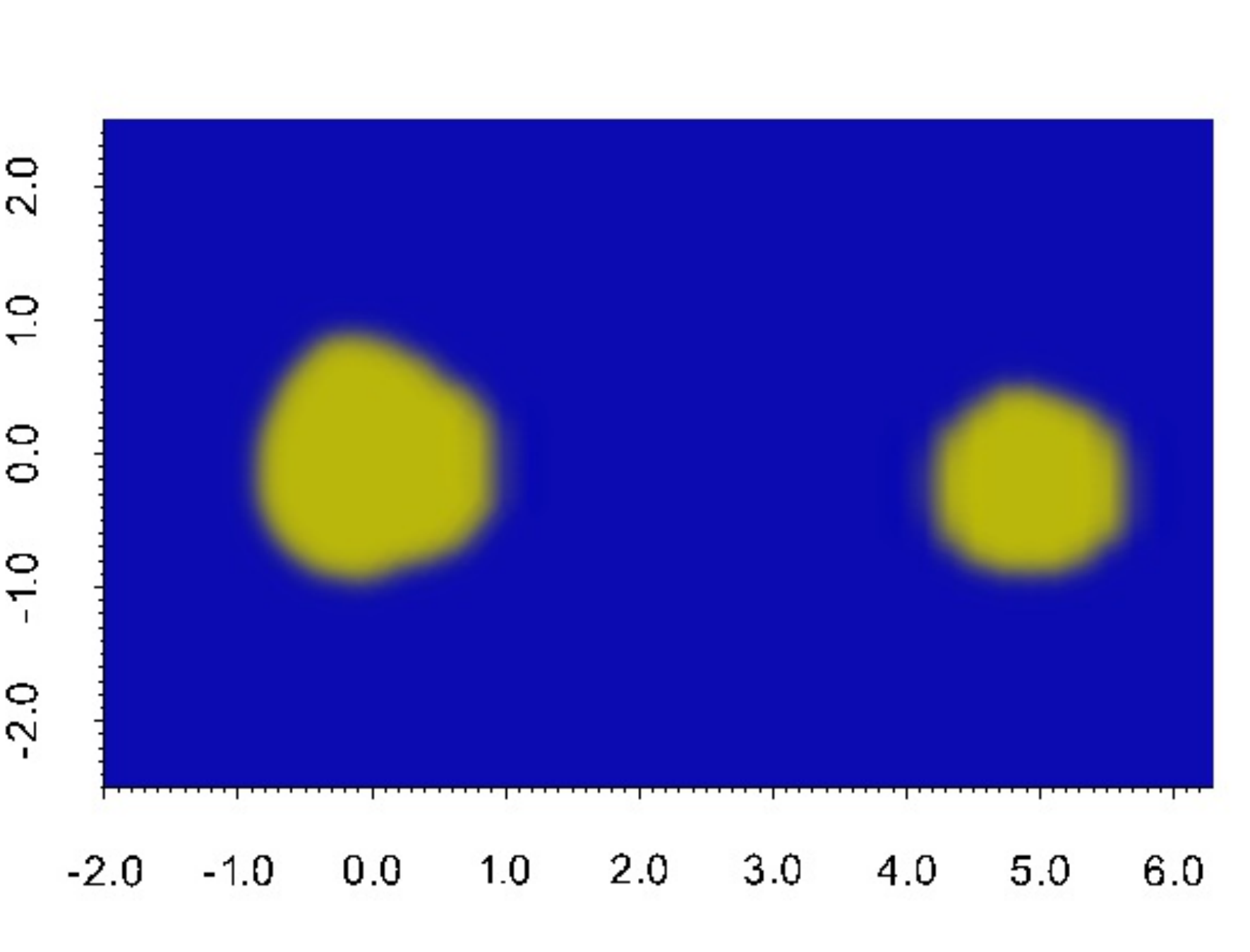}}}
\subfigure[Target data ($\vpo$).]{\label{fig:synthetic_data_PF_target}
{\includegraphics[trim = 10mm 140mm 20mm 20mm,clip, width=0.48\linewidth]{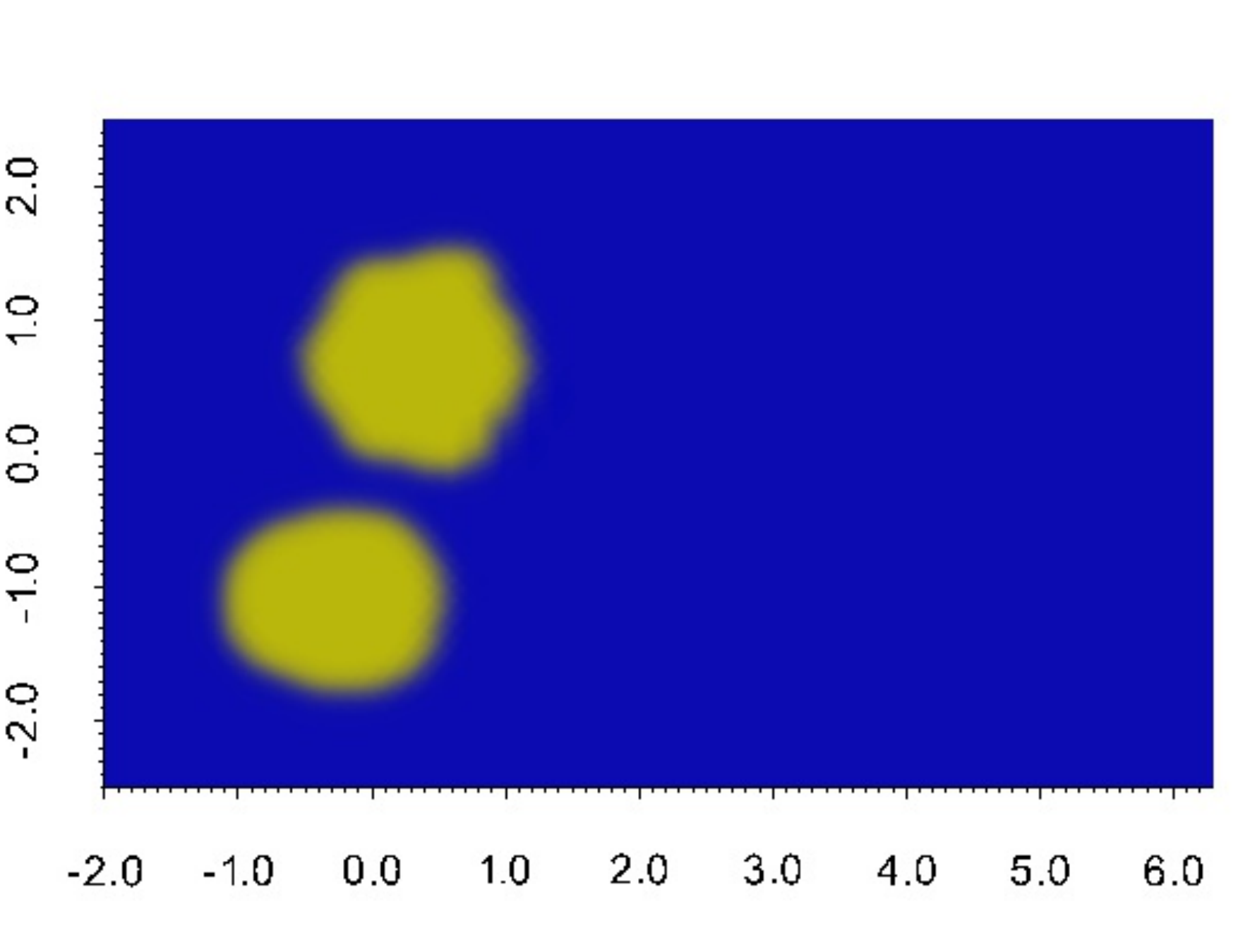}}}
\end{center}
 \caption{Initial and target data for the examples of \S \ref{subsec:results-multi-splits}.}
  \label{fig:splits_initial_target}
   \end{figure}
     \begin{figure}[h!]
\begin{center}
\subfigure[Without the volume constraint]
{{\includegraphics[trim = 10mm 0mm 20mm 0mm,clip,width=0.49\linewidth]{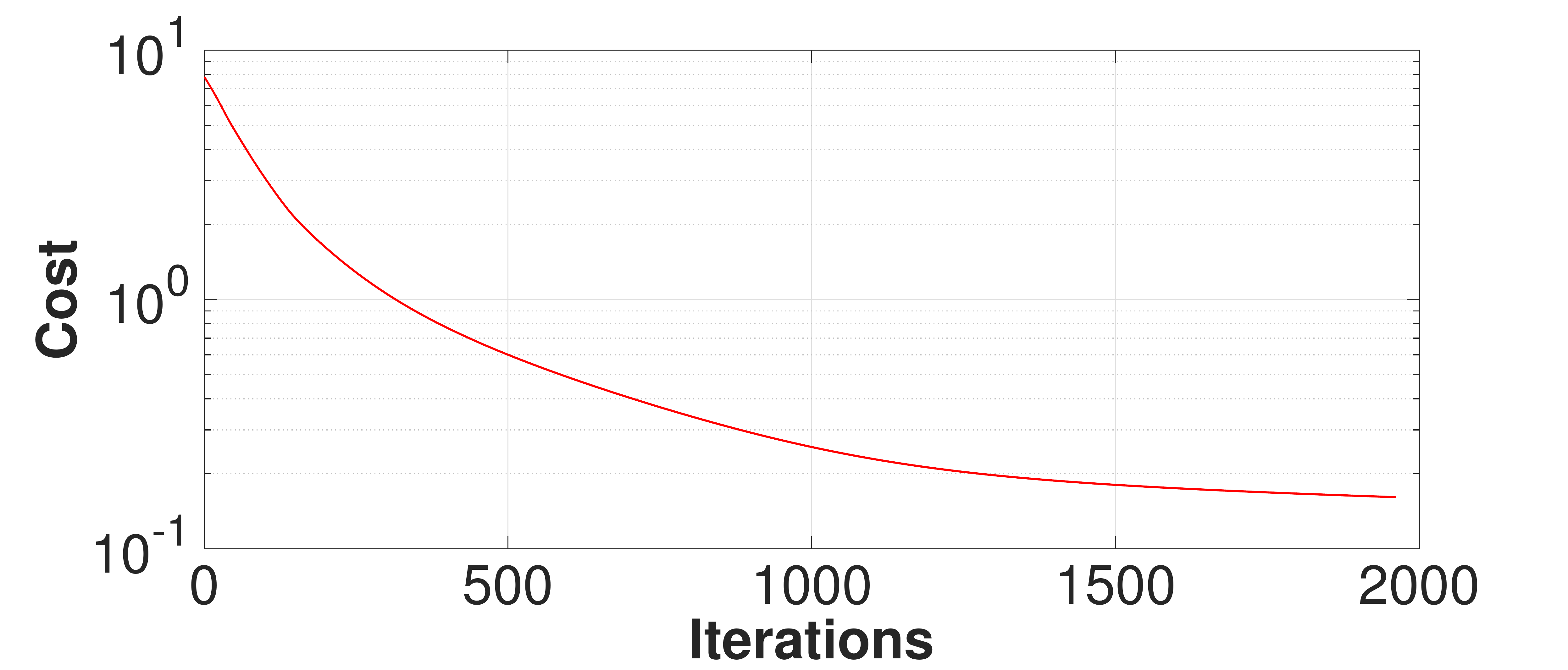}}}
\subfigure[With the volume constraint]
{{\includegraphics[trim = 10mm 0mm 20mm 0mm,clip,width=0.49\linewidth]{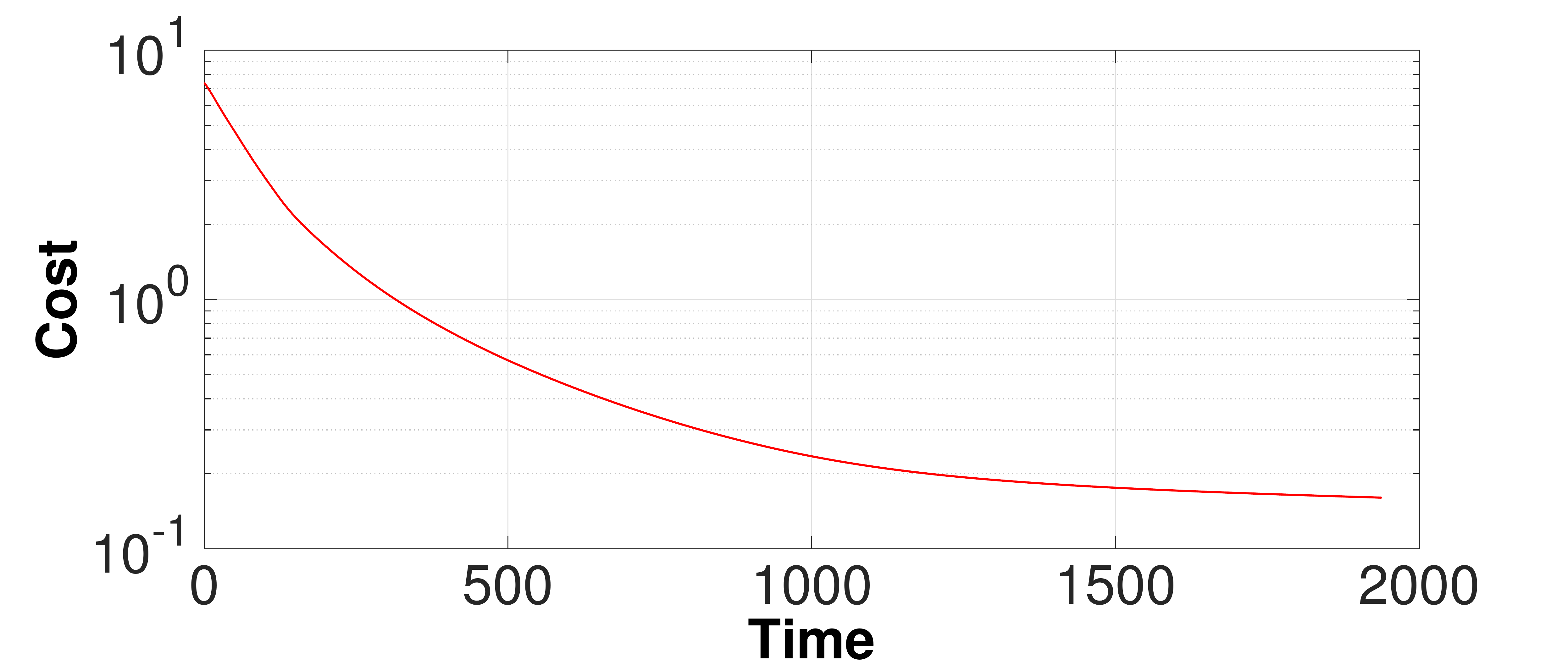}}}
\end{center}
 \caption{The value of the cost functional versus the number of iterations for the examples of \S \ref{subsec:results-multi-splits} with and without the volume constraint. We observe a rapid decrease in the cost initially followed by a much more gradual  decrease as we approach the minimum, this is as expected since the steepest descent algorithm is used for the update of the control.}\label{fig:cost_splits}
  \end{figure}
  \begin{figure}[h!]
\begin{center}
\subfigure[Without the volume constraint]
{{\includegraphics[trim = 10mm 140mm 40mm 10mm,clip,  width=0.49\linewidth]{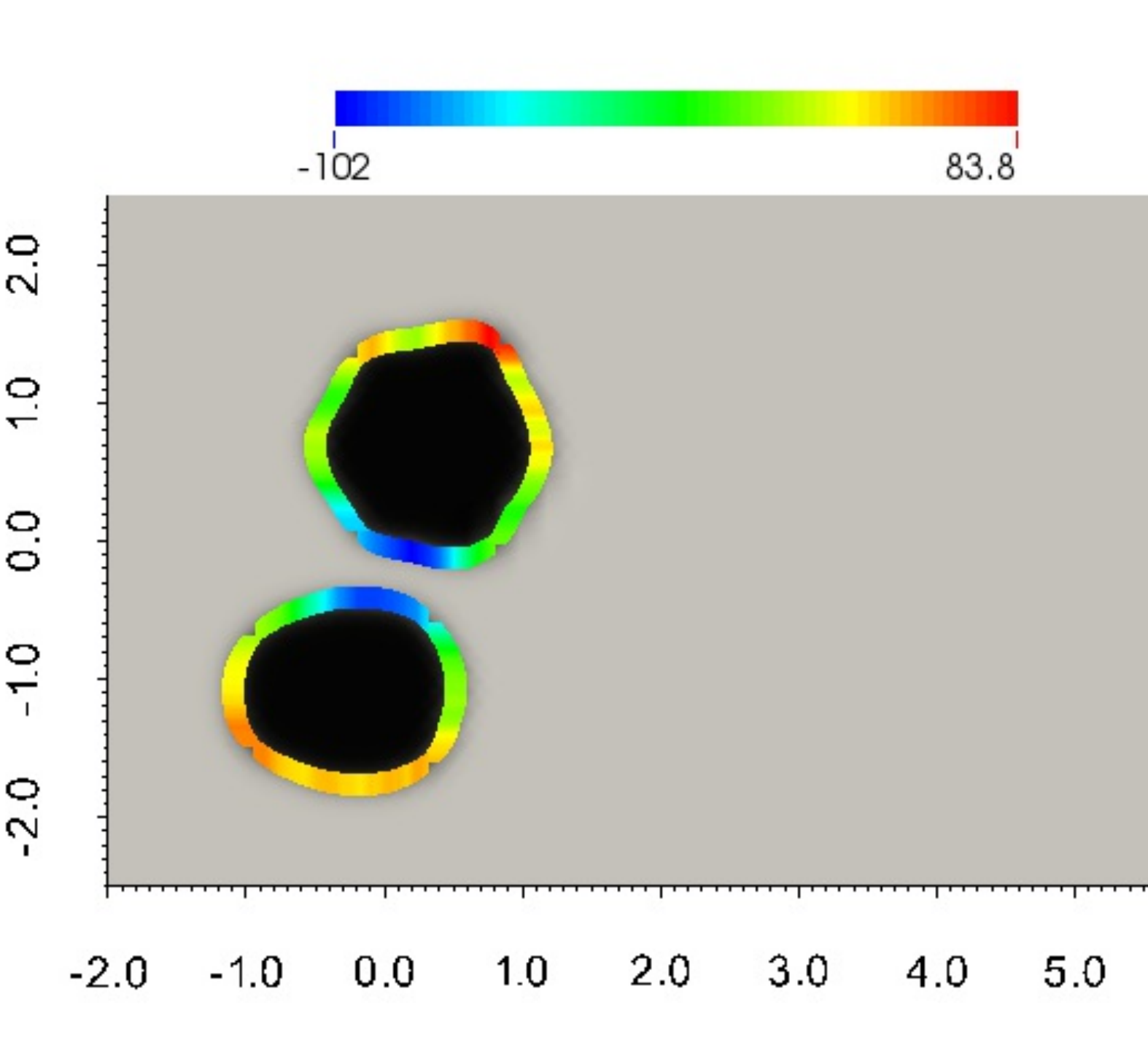}}}
\subfigure[With the volume constraint]
{{\includegraphics[trim = 10mm 140mm 40mm 10mm,clip,  width=0.49\linewidth]{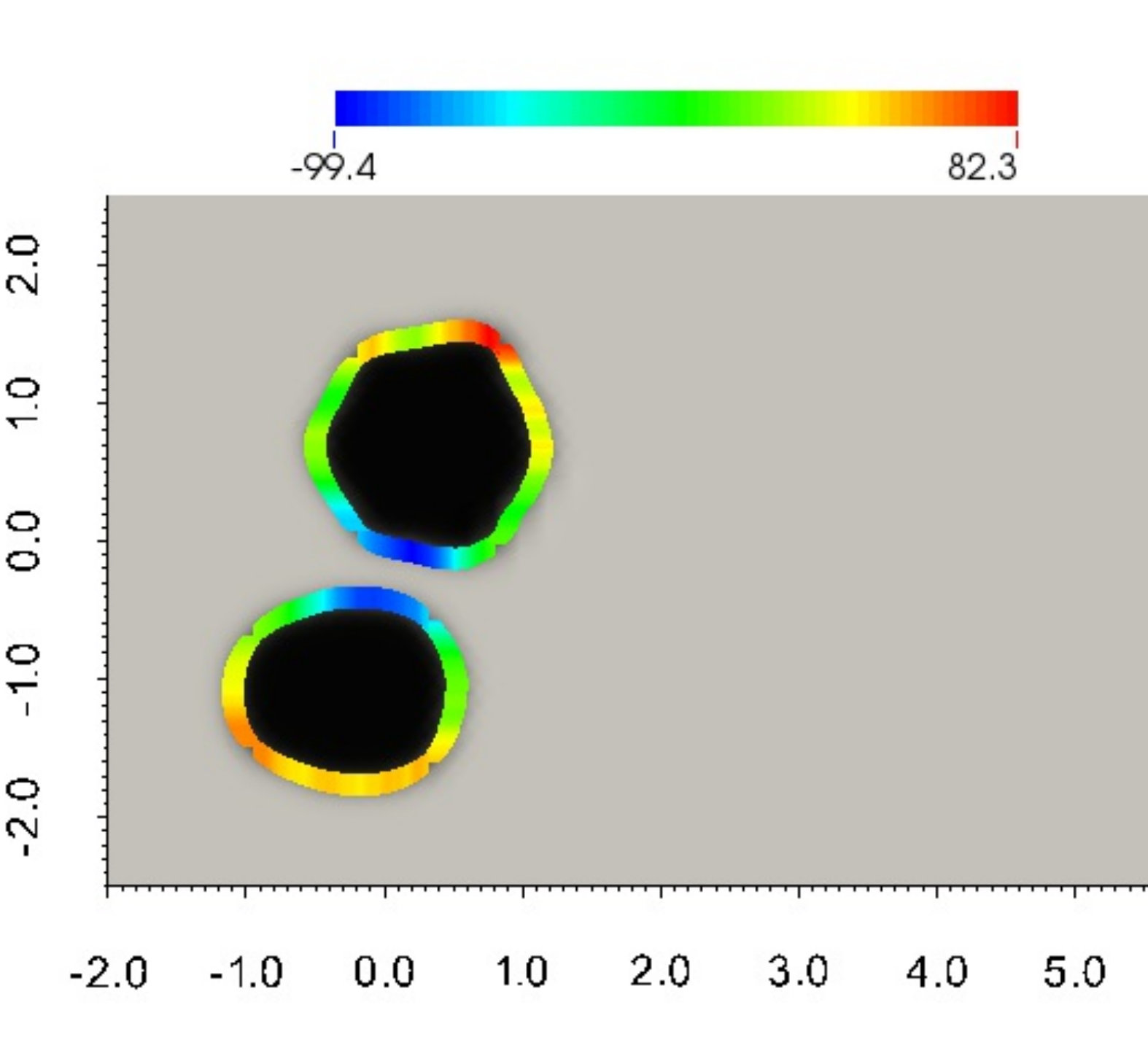}}}
\end{center}
 \caption{
 Zero level-set of the solutions ($\vp(\vec x,T)$) computed using the approximated optimal control ($\eta^*(\vec x,t)$) with and without the volume constraint for the examples of \S \ref{subsec:results-multi-splits}.
  The curve (zero level-set of $\vp(\vec x,T)$)  is shaded by the approximated optimal control ($\eta^*(\vec x,T)$) and the background by the target data ($\vpo(\vec x)$). The color-bar corresponds to the scale for $\eta^*(\vec x,T)$. 
  We see good agreement between the zero level-set of the data computed with the optimal control and the target data in both cases.}\label{fig:splits_optimal_target}
  \end{figure}
    \begin{figure}[h!]
\begin{center}
{{\includegraphics[trim = 40mm 0mm 30mm 0mm,clip,width=\linewidth]{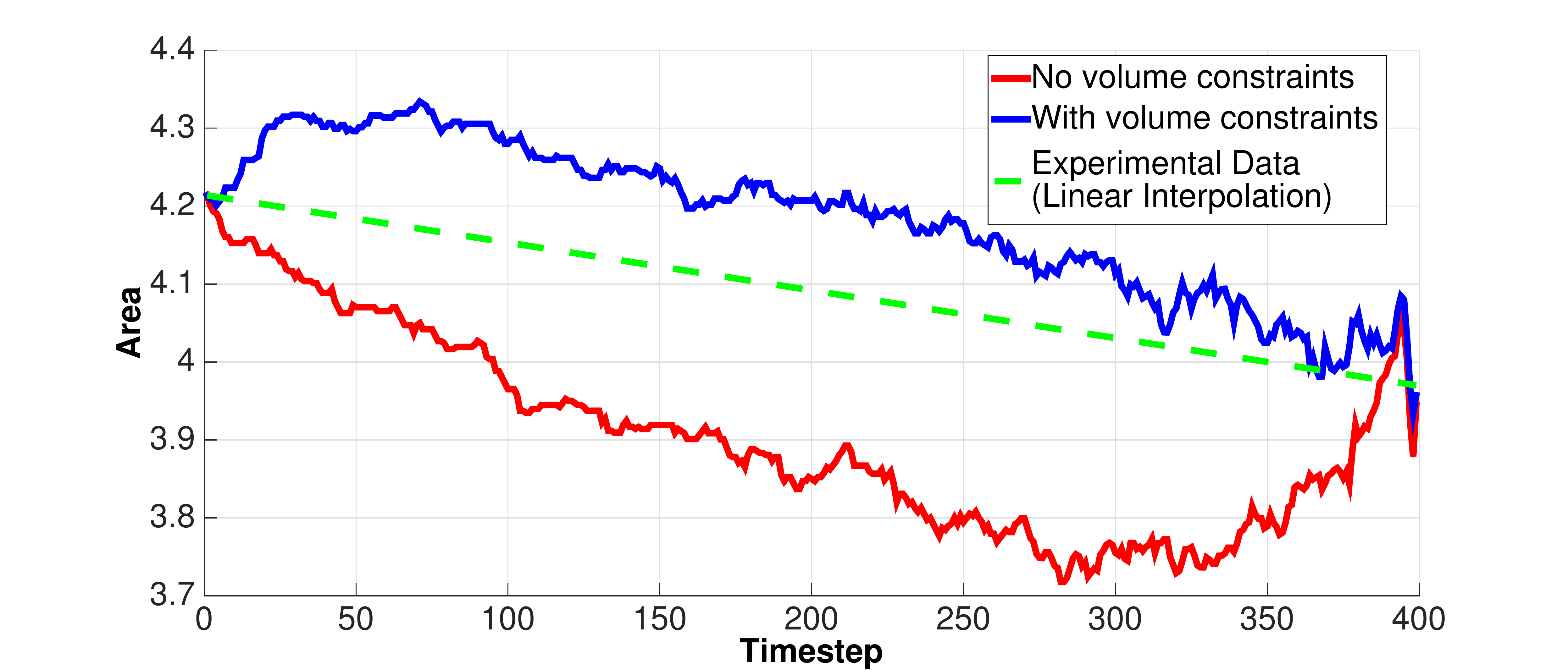}}}
\end{center}
 \caption{Area enclosed by the cell for the experiments of \S \ref{subsec:results-multi-splits} with and without the volume constraint. As with the single cell data, the area (now the sum of the area of the two cells) shrinks considerably during the evolution without the volume constraint, whilst the incorporation of the volume constraint yields a better fit to the linear interpolant of the areas. Unlike the previous examples however, even with the volume constraint the fit to the linear interpolant of the areas of the data is poor.}\label{fig:area_splits}
  \end{figure}
\begin{figure}[h!]
\begin{center}
\subfigure[Without the volume constraint]
{{\includegraphics[trim = 0mm 150mm 0mm 30mm,clip,  width=0.45\linewidth]{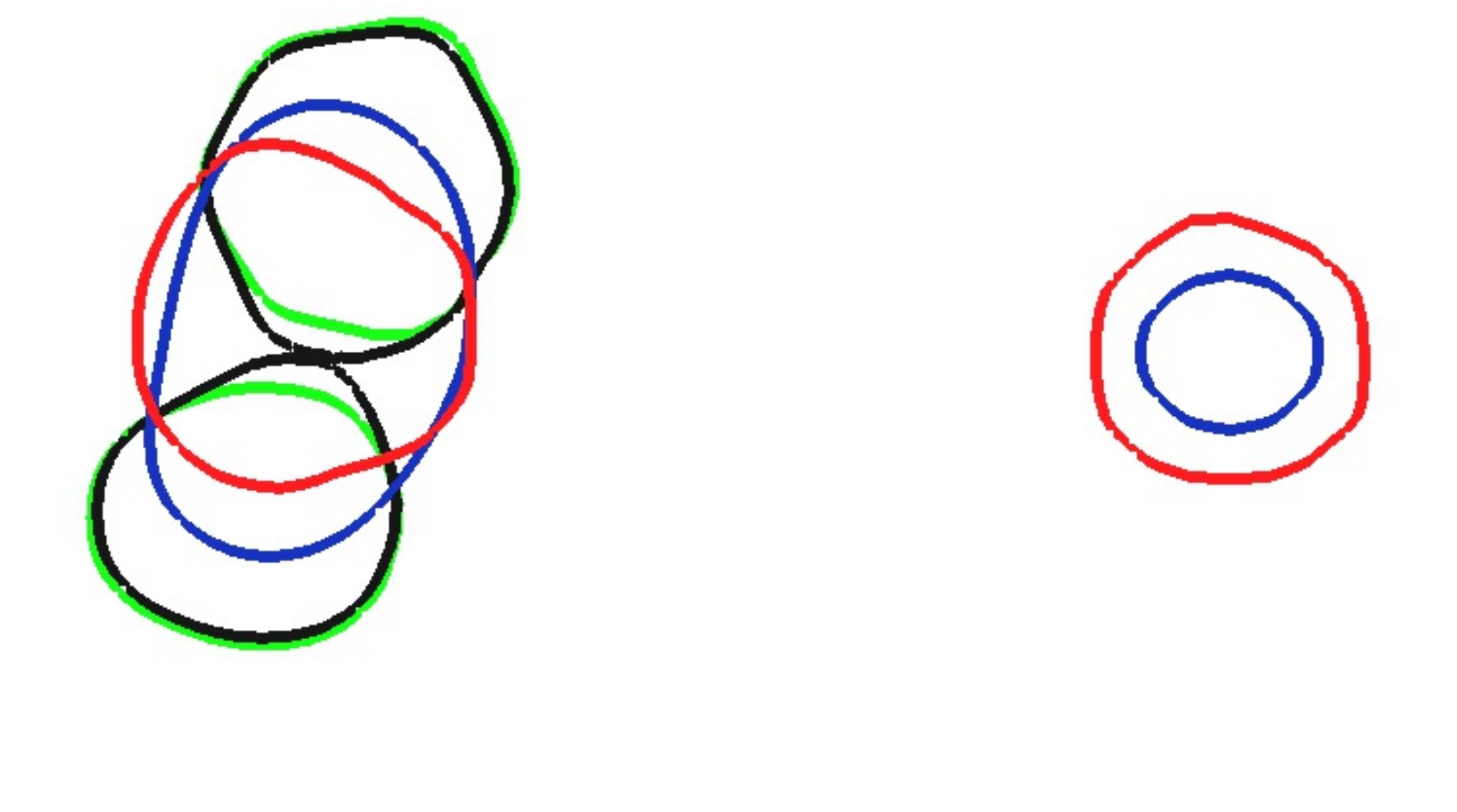}}}
\hskip1em
\subfigure[With the volume constraint]
{{\includegraphics[trim = 0mm 150mm 0mm 30mm,clip,  width=0.45\linewidth]{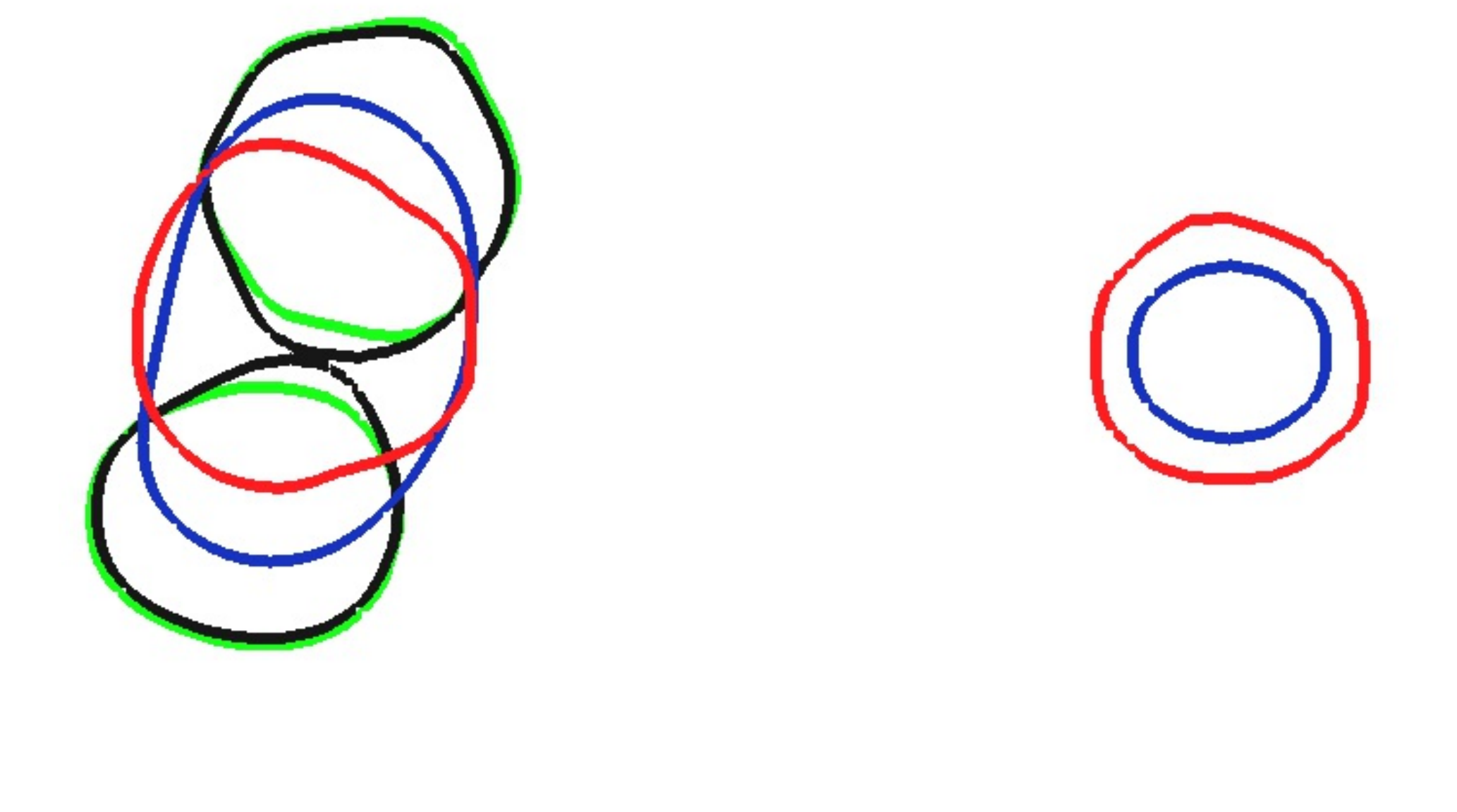}}}
\end{center}
 \caption{Zero level-sets of the solutions computed with the optimal control for the examples of \S \ref{subsec:results-multi-splits} with and without the volume constraint at $t=0$ (red), $t=0.3$ (blue),  $t=0.396$ (black) and $t=0.4$ (green). We observe that in this case the difference between the two schemes is less pronounced. It is also clear that both with and without the volume constraint,  the implicit solution of the matching problem in our algorithm in this case leads to the annihilation of one cell (as it shrinks to a point) while the other cell splits with the zero level-set changing in topology from a single closed curve to two disjoint closed curves.}\label{fig:splits_optimal_curves}
  \end{figure}
\subsection{Comments on the numerical experiments}
 The CPU times for each of the experiments  is of the order of hours. For all the experiments the number of iterations required before the stopping criteria is met are similar, however this leads to simulations with the volume constraint taking approximately four times as long (in terms of CPU time) as those without the volume constraint. This is due to the iterative nature of the algorithm used to compute the Lagrange multiplier c.f., \S \ref{app:FEM} which necessitates multiple solves per timestep. We note that  the CPU times of the algorithms may be too large for many applications. In light of this we mention that the stopping criteria we have used may be too strict for many applications and that a significant decrease in the cost function together with a reasonable fit to the data is observed after as few as 50 iterations (which reduces the CPU time by a factor of around 40) and that for many applications this level of fit may be sufficient hence the stopping criteria could be relaxed. Finally we mention that the current solution procedure based on uniform grids and serial solution of the forward and adjoint problems may be improved and we are currently investigating combining adaptive grids with a parallel solver for the forward and adjoint problems which gives a significant speed up but presents new technical challenges which we wish to avoid in this paper to maintain clarity of exposition.
\section{Conclusion}
\label{sec:conc}
In this study we presented a first step towards the development of cell tracking algorithms based on physical models for cell migration. The presented algorithm seeks to track whole cell morphologies and is applicable to single cell or multi-cell image data sets. Our approach may be regarded as a model fitting procedure in which a physically derived model for the evolution of the cell or cells is fitted to experimental image data sets. The algorithm is based on the theory of  optimal control of PDEs and full details of the derivation and implementation of the algorithm are given.   We also present a number of numerical experiments illustrating the performance of the algorithm with synthetic representative single cell and multi-cell image data sets.

The key novelty of our approach is that the model for the evolution of the cell (or cells), which drives the tracking procedure, is based on a relevant simplification of existing physically derived models for cell motility that reproduce many experimentally observed aspects of cell migration (e.g., \cite{venk11chemotaxis}). Thus this study is an important step towards the development of cell tracking algorithms in which the recovered trajectories are physically meaningful. This is in contrast to the majority of existing algorithms for whole cell tracking in which the models for the evolution of the cell that underly the tracking procedure are purely geometric in nature neglecting completely the physics of cell migration \cite{bosgraaf2009analysis,tyson2010high,kadirkamanathan2012neutrophil}.
One significant advantage of the approach to cell tracking we propose, is that the physics of the model driving the evolution of the cell are reflected in the recovered dynamic data. Thus it is possible to encode physical features of cell migration into the tracking procedure. We illustrate this fact by including volume conservation in the model.
 Comparing the results of the tracking algorithm with and without volume conservation, we observe, that in  a number of simulations neglecting volume conservation leads to physically unrealistic cell morphologies 
with a significant reduction of the cell volume in the recovered morphologies whilst this undesirable effect is no longer evident if volume conservation is included. 
 We note that volume conservation is more physically relevant in three dimensions.  This is since large changes in volume in two-dimensional imaging data of cells 
(i.e., the area of the projections of the cell on to a two-dimensional plane) are often observed in experimental data despite the cells conserving their enclosed (three-dimensional) volume.

Of course volume conservation is only an example of the kind of biology or physics one may wish to encode in the algorithm. A number of models that fit into our framework have been proposed incorporating more complex biophysics such as spontaneous curvatures \cite{marth2013signaling}, for Helfrich type models, adhesion for the migration of cells on substrates or in the ECM \cite{shao2012coupling}, cell-cell or cell-obstacle interactions \cite{venk11chemotaxis} and chemotaxis \cite{neilson2011chemotaxis}, these models are thus potential candidates for the model driving the evolution in our tracking algorithm.

 We work with  diffuse interface representations of the cell membrane to make use of the mature theory for the optimal control of semilinear PDEs. One attractive aspect of this approach is that, as we do not require sharp interface representations of the cell membrane, it may be possible therefore to work directly with the raw experimental image data set without any need for segmentation.  However the diffuse interface or phase field framework we employ does make the algorithm computationally intensive as evidenced by the relatively large CPU times for our experiments and a key area for future work is to investigate improvements in the  computational efficiency of the algorithm. This need is especially evident if one wishes to track cells in $3d$, as although our theoretical framework applies equally to this setting the computational cost becomes prohibitive. Computational aspects under investigation include
 \begin{itemize}
 \item  Spatial and temporal adaptivity which is challenging in this setting as the solution of the state equation enters the adjoint equation.
 \item  Alternative update schemes for the control to the simple yet robust gradient based update considered in this study.
 \item Parallelisation and the development of fast solvers for the solution of the state and adjoint equations.
 \end{itemize}
 Our initial numerical investigations suggest that with a combination of the techniques outlined above it is possible to efficiently track $3d$ cell migration and we report on this elsewhere. 
 Other potentially attractive directions for future work would be to consider higher order finite element spaces for the discretisation of the forward and adjoint problems or the use of spectral element methods, both of which may allow a more accurate solution of the forward and adjoint problem with fewer degrees of freedom, hence reducing the memory requirements.

Investigating the performance of the algorithm with real biological data for different cell types and in different environments is an important and worthwhile task.  We are currently applying the algorithm to the tracking of in vivo neutrophil migration and intend to report on this elsewhere. As mentioned previously one interpretation of the forcing $\eta^*$ is that it accounts for both protrusive forces generated by polymerisation of actin at the leading edge of the cell together with contractile forces generated by the action of myosin motors at the cell rear. Thus a potential avenue for assessing the plausibility of the  cell tracks computed with our algorithm would be to compare the computed $\eta^*$ with experimental imaging data on the location of polymerised actin and myosin-II on the cell membrane with the expectation being that regions in which the computed forcing $\eta^*$ is positive would correspond to regions rich in polymerised actin and regions in which the computed forcing $\eta^*$ is negative would correspond to regions rich in myosin-II.
There are also many extensions of our approach which are likely to prove useful in applications. Our algorithm could equally be applied to the identification of (possibly time-dependent) parameters in models for cell migration (e.g., a spatially constant forcing or material parameters such as surface tension or bending rigidity) however in this case it is likely that the sharp interface approach we propose in \cite{2013arXiv1311.7602C} will be more efficient.  As observed in some of the experiments we report on, the framework we employ allows changes in topology of the cells.  Whilst this may be desirable for some applications, e.g., tracking cells beyond cell division or cell fusion, in many biological experiments the topology of the cells is fixed. Our experiments suggest that topological changes arise primarily in the case of multi-cell image data sets. In this setting it should be possible to track the evolution of certain topological invariants (or more specifically diffuse interface representations of such invariants) and use these as an indicator for when the computed cells are changing in topology. The user could then manually reduce the multi-cell tracking problem to multiple single cell (or smaller scale multi-cell) tracking problems by specifying the correspondence between cells in different frames, with the hope that changes in topology do not occur for these new problems. The model we propose for the evolution in this study is a simplification of more general physically relevant models in which bulk or surface PDEs for the biochemistry are coupled to a geometric evolution law for the motion. An important area for future work is the extension of the framework to this more general setting. We note that the phase field approach we employ makes it computationally straightforward to couple the geometric evolution law for the motion to bulk PDEs (posed either within the cell or in the extra-cellular matrix) \cite{ziebert2011model,shao2010computational,shao2012coupling}. 
 
We hope that our optimal control-model fitting based framework is a useful first step towards incorporating advances in the modelling of cell migration into cell tracking algorithms. 

\section*{Acknowledgments}
This work (AM, VS and CV) is supported by the Engineering and Physical Sciences Research Council, UK grant (EP/J016780/1) and the Leverhulme Trust  Research Project Grant (RPG-2014-149). KB was partially supported by the Embirikion Foundation Grant (2011-2014)-Greece. The computations were carried out using the computational resources of the School of Mathematics and Physical Sciences at the University of Sussex. 

 \appendix
\section{Numerical solution of the forward and adjoint problems}
\label{app:FEM}
\subsection{Discretisation of the state equation}
We introduce the variational form  for the forward problem (\ref{eqn:pf}) defined as follows. 
Find $(\vp,\lambda)\in\Lp{2}([0,T];\Hil{1}(\O))\times L^2(0,T)$ such that 
\begin{align*}
\int_{\Omega}\partial_t\varphi\psi d{\bf{x}}+\int_{\Omega}\nabla\varphi\cdot\nabla\psi d{\bf{x}}
=\frac{1}{\varepsilon}\int_{\Omega}(c_{G}\eta-\lambda)\psi d{\bf{x}}-\frac{1}{\varepsilon^2}\int_{\Omega} G^{\prime}(\varphi)\psi d{\bf{x}}
\myall \psi\in\Hil{1}(\O).
\end{align*}
Let $\T$ be a decomposition of $\O$ into simplexes (for simplicity we assume $\O$ is polygonal). We define the finite element space
\beq
\V:=\lbrace \psi_h\in \Hil{1}(\O)\cap C^{0}(\O):\psi_h|_k\in\mathbb{P}^1\quad\forall k\in\T\rbrace.
\eeq
For the time discretisation employ an implicit-explicit method where the diffusive term is treated implicitly and the reaction terms explicitly. Introducing the shorthand for a time discrete sequence $f^n:=f(t^n)$ and a uniform timestep $\tau$ with $T=M\tau, M\in\mathbb{N}$, the fully discrete scheme reads, 
 for $n=0,\dots,M-1$, given $\vp_h^{n},\eta_h^{n}\in\V$ find  $(\vp_h^{n+1},\lambda^{n+1})\in\V\times\Reals$ such that 
\begin{align*}
  \frac{1}{\tau}\int_{\Omega}(\vp_h^{n+1}-\vp_h^{n})&\psi_h d{\bf{x}}+\int_{\Omega}\nabla\vp_h^{n+1}\cdot\nabla\psi_h d{\bf{x}}=\\
&\frac{1}{\varepsilon}\int_{\Omega}(c_G\eta_h^n-\lambda^{n+1})\psi_h d{\bf{x}}
-\frac{1}{\varepsilon^2}\int_{\Omega} \Lagrange(G^{\prime}(\vp_h^n))\psi_h d{\bf{x}}
\myall\psi_h\in\V,
\end{align*}
where $\Lagrange:C^{0}(\O)\to\V$ denotes the Lagrange interpolant. 

We solve the above problem using the iterative technique introduced and studied in \cite{blowey1993curvature}, which uses a bisection
method for the Lagrange multiplier. In particular we seek an iterative sequence $\{\phi_h^{n+1,k},\lambda^{n+1,k}\}_{k\geq 1}$ where 
$\phi_h^{n+1,k}$ solves 
\begin{align*}
  \frac{1}{\tau}\int_{\Omega}(\vp_h^{n+1,k}-\vp_h^{n})&\psi_h d{\bf{x}}+\int_{\Omega}\nabla\vp_h^{n+1,k}\cdot\nabla\psi_h d{\bf{x}}=\\
&\frac{1}{\varepsilon}\int_{\Omega}(c_G\eta_h^n-\lambda^{n+1,k})\psi_h d{\bf{x}}
-\frac{1}{\varepsilon^2}\int_{\Omega} \Lagrange(G^{\prime}(\vp_h^n))\psi_h d{\bf{x}}
\myall\psi_h\in\V,
\end{align*}
with $\lambda^{n+1,1}=-\frac{2\varepsilon}{\tau}+1$, $\lambda^{n+1,2}=\frac{2\varepsilon}{\tau}-1$ and $\{\lambda^{n+1,k+1}\}_{k\geq 2}$ satisfying
$$
\lambda^{n+1,k+1}=\lambda^{n+1,k}
+\frac{\left(\lambda^{n+1,k}-\lambda^{n+1,k-1}\right)\left(M_\vp^{n+1}-\int_{\Omega}[\vp_h^{n+1,k}]_{+}\right)}{\left(\int_{\Omega}[\vp_h^{n+1,k}]_{+}-\int_{\Omega}[\vp_h^{n+1,k-1}]_{+}\right)},
$$
where we recall
$$M_{\vp}^{n+1}:=\int_\O[\vp_h^0]_{+}+\frac{(n+1)\tau}{T}\left([\vpo]_{+}-[\vp_h^0]_{+}\right)\diff \vec x.$$
We deem this iteration to have converged when $|\lambda^{n+1,k+1}-\lambda^{n+1,k}|< \mbox{tol}$. 

The discretisation of the forward problem without the volume constraint is as above with $\lambda=0$. The discretisation of the adjoint problem without the volume constraint is the same as above as the Lagrange multiplier $\lambda(t)$ does not enter the adjoint problem (\ref{eqn:adjoint}).

\end{document}